\documentclass[10pt]{elsarticle}

\usepackage{charter} 
\usepackage{etex}
\usepackage{ifpdf}
\usepackage{enumerate}
\usepackage[inner=1.2cm,outer=1.2cm, top=2cm, bottom=2cm]{geometry}
\usepackage{amsbsy,amssymb,amsfonts,amsmath}                      
\usepackage{algorithm,algpseudocode} 
\algnewcommand\algorithmicinput{\textbf{Input:}}
\algnewcommand\Input{\item[\algorithmicinput]}
\algnewcommand\algorithmicoutput{\textbf{Output:}}
\algnewcommand\Output{\item[\algorithmicoutput]}
\usepackage{plain}

\usepackage{mathtools}
\newcommand{\bsb}[1]{\boldsymbol{#1}}
\DeclareMathOperator{\spn}{span}

\usepackage{graphicx,subfigure,epstopdf,epsfig}
\usepackage{rotating}
\usepackage{tabularx}
\usepackage[font=small,labelfont=bf]{caption}  
\usepackage[usenames,dvipsnames]{pstricks}
\usepackage{pst-grad} 
\usepackage{pst-plot} 
\usepackage{tikz}
\usepackage{hyperref}
\usepackage{colortbl}
\usepackage{framed, color, soul}
\definecolor{shadecolor}{gray}{0.8}
\sethlcolor{yellow}    
\usepackage{xcolor}
\usepackage{booktabs}

\begin{document}
\begin{frontmatter}



\title{Efficient Estimation of Cardiac Conductivities via POD-DEIM Model Order Reduction}
\author{Huanhuan Yang}
\ead{huan2yang@outlook.com}
\author{Alessandro Veneziani}
\ead{ale@mathcs.emory.edu}
\address{Department of Mathematics and Computer Science, Emory University, Atlanta (GA) USA}


\begin{abstract}
Clinical oriented applications of computational electrocardiology require efficient and reliable identification of patient-specific parameters of mathematical models based on available measures. 
In particular, the estimation of cardiac conductivities in models of potential propagation is crucial, since they have major quantitative impact on the solution.
Available estimates of cardiac conductivities are significantly diverse in the literature and the definition of experimental/mathematical estimation techniques is an open problem with important practical implications in clinics. 
We have recently proposed a methodology based on a variational procedure, where the reliability is confirmed by numerical experiments. 
In this paper we explore model-order-reduction techniques to fit the estimation procedure into timelines of clinical interest. 
Specifically we consider the Monodomain model and resort to Proper Orthogonal Decomposition (POD) techniques to take advantage of an off-line step when solving iteratively the electrocardiological forward model online. 
In addition, we perform the Discrete Empirical Interpolation Method (DEIM) to tackle the nonlinearity of the model. 
While standard POD techniques usually fail in this kind of problems, due to the wave-front propagation dynamics, an educated novel sampling of the parameter space based on the concept of Domain of Effectiveness introduced here dramatically reduces the computational cost of the inverse solver by at least 95\%. 
\end{abstract}


\end{frontmatter}


\section{Introduction}
In recent years, computational electrocardiology has attracted the attention of mathematicians, engineers and clinicians. Computational methods have been continuously refined to match clinical applications. 
Just to mention one example,the ideal lesion pattern with minimal burn in a cardiac ablation therapy has been studied in \cite{Krueger2013} by electrocardiological modeling. 
However, to provide reliable and efficient simulations of cardiac electrical activity is not easy. Differential models in electrocardiology depend on several parameters typically coming from empirical constitutive laws, so their quantification for a specific patient is difficult. 
In particular, these models 
are strongly sensitive to the cardiac conductivity parameter \cite{johnston2011sensitivity}. While experimental data generally disagree on the values of conductivities, mathematically sound estimation methods---generally based on the solution of an inverse problem---have been considered only recently in \cite{graham2010estimation} and in our previous work \cite{HH2015IP}. 
In addition, electrocardiological modeling for clinical application is computationally intensive. This is even more true for the inverse conductivity problem as we carried out in \cite{HH2015IP}, since high computational cost arises in many ``queries'' of forward simulations with different conductivity guesses.
Model-order-reduction techniques have been investigated in the literature, but their application to cardiac conductivity estimation is challenging, due to the nonlinearity of the models and the specific features like wave-front propagation of the solutions. In this paper, we apply model reduction techniques to dramatically decrease the computational cost of solving the inverse conductivity problem.

A numerical solution to the forward problem by Galerkin projection can be represented, in general terms, as an expansion ${u_h}(\boldsymbol{x},t)=\sum\limits_{j=1}^{n} u_j(t)\phi_j(\boldsymbol{x})$. In finite elements, the basis $\{ \phi_j \}_{j=1}^n$ is selected to be a set of piecewise polynomials. This basis is of ``general purpose'' as it does not have any specific clue of the problem to solve. 
Model reduction techniques aim at cheaply solving the forward equations in a low-dimensional space still by a Galerkin projection process. 
To this aim they construct a rather small set of basis functions (known as {\it reduced basis}), which we call ``educated'' basis as it includes features of solutions to the forward problem considered. This allows to an efficient low-dimensional yet accurate representation of the solution.

Among various techniques for the reduced basis construction (see the recent review\cite{Benner2015}), typical ones are the {\it Proper Orthogonal Decomposition} (POD) approach \cite{Kunisch2002, frangos2010surrogate}, the {\it Greedy Reduced Basis} (GRB) approach \cite{quarteroni2011certified, maday}, and their combination \cite{Nguyen2010}. 
In these approaches, the reduced basis is constructed from a set of parameter-dependent solutions of the full-order model. 
The POD approach has drawn widespread attention for its optimal ability to approximate the snapshots of solution with minimized error by the selection of the most important modes, and for its easy-to-use feature in practice. The GRB approach selects the snapshots following a greedy process according to a rule controlled by an {\it a posteriori} error estimator. 
The requirement of a rigorous error estimator currently limits its application to electrocardiological models due to the model complexity. 
In this paper we resort to the POD approach as a starting point, the development of an {\it a posteriori} error estimator for the GRB method is a part of ongoing work.

The POD approach has been used in numerous fields of science and engineering such as fluid-structure interaction \cite{BertagnaVeneziani} and aerodynamics \cite{Thanh2004}, but its practical application in electrocardiology only starts from 2011 \cite{Boulakia2011POD, boulakia2012, corrado2015identification}. In these references, the POD method allows reasonable estimation of cardiac ionic model parameters, however, no systematic study on the improvement of efficiency of solving the full nonlinear electrocardiological model is available. 
In fact, one critical aspect when reducing a nonlinear problem by projecting onto a low-dimensional space is to approximate the projected nonlinear terms in a way independent of the full-order model size. This point was not thoroughly addressed in current electrocardiology publications.

Several techniques are available in general to reduce the cost of evaluating nonlinear terms, such as the {\it trajectory piecewise-linear} (TPWL  \cite{TPWL2003}) approach, the {\it Best Point Interpolation Method} (BPIM \cite{NguyenBestPoint}), the {\it Empirical Interpolation Method} (EIM \cite{Barrault2004667})  and its discrete variant Discrete EIM (DEIM \cite{thesisDEIM}). Although TPWL was successfully applied to some practical problems, it may not be effective or efficient for systems with high order of nonlinearity. The BPIM and EIM approaches are similar \cite{Galbally2010}, both select a small set of spatial interpolation points to avoid the expensive calculation of inner products and use the points for nonlinear approximation. BPIM is optimal in the point selection and gains a little improvement on accuracy, but it is more computational expensive.

Here, we resort to EIM in its discrete variant DEIM. 
Precisely, we apply the POD-DEIM for the first time to the conductivity estimation problem.
It is worth mentioning that the conductivity parameter to be estimated considered in this work is more troublesome than other ionic model parameters, since it dominates the speed and direction of fast transient of electrical potential through the cardiac tissue, which is an intrinsic feature of the forward electrocardiology model. This fact thus prevents a successful model reduction via a 
classical POD procedure. 
Model reduction procedures need to be specifically customized for the problem, and in particular the construction of the educated basis is a delicate step. 
Nevertheless, we show here how an appropriate sampling for the basis computation actually leads to significant reduction of the full-order computational costs with a great level of accuracy. We address the sampling required for basis construction based on the novel concept of ``Domain of Effectiveness'' in the parameter space. 
A rather small set of samples is obtained by sampling the parameter space based on polar coordinates, with refinement in the ``small angle--short arc'' zone of the sample space utilizing Gaussian nodes. 
In this way, we manage to use the POD-DEIM reduced-order model with a computational reduction
of at least 95\% of the full-order conductivity estimation.
The present work relies on but largely improved the study presented in \cite{HHthesis}, to which we refer for more details.

We consider specifically the Monodomain problem. Notwithstanding that some authors consider this model reliable enough for many clinical applications, the extension of the present work to the more accurate
Bidomain system is a follow-up of the present work.
The outline is as follows. After a brief statement of the full-order {\it Monodomain inverse conductivity problem} (MICP) in a discrete form (Sec.~\ref{MICPfull}), the POD-DEIM approach is introduced (Sec.~\ref{DEIMintr}) and applied to solve a reduced MICP (Sec.~\ref{DEIM-use}) by derivative-based optimization. The reduced-order model is tested in Sec.~\ref{DEIMresult}: both the efficiency and accuracy of the POD-DEIM approach in conductivity estimation are investigated; we report both pitfalls and successful examples.

\section{The full-order Monodomain inverse conductivity problem}\label{MICPfull}
\subsection*{The cardiac Bidomain and Monodomain models}
The {\it Bidomain model} is considered as the most physiologically founded description for the dynamics of cardiac electric potentials---the transmembrane potential $u$ and the extracellular potential $u_{\rm e}$---at the level of cardiac tissue. 
Its parabolic-elliptic form (see e.g. \cite{pullan2005mathematically}) reads
\begin{equation} \label{stateequations}
\left\{
\begin{array}{ll}
\beta C_{\rm m} \dfrac{\partial u}{\partial t}
-\nabla \cdot \left( \boldsymbol{\boldsymbol{\sigma}}_{\rm i} \nabla u \right)-\nabla \cdot \left( \boldsymbol{\sigma}_{\rm i}\nabla u_{\rm e} \right) + \beta I_{\rm ion}(u,w) =I_{\rm si} &\hspace{0.05cm}\mbox{ in } \Omega\times[0,T]  \\[0.3cm]
- \nabla \cdot \left( \boldsymbol{\sigma}_{\rm i} \nabla u \right)-\nabla \cdot \left(\boldsymbol{\sigma}_{\rm i}+\boldsymbol{\sigma}_{\rm e}\right)\nabla u_{\rm e}  =I_{\rm si}-I_{\rm se} & \hspace{0.05cm}\mbox{ in } \Omega\times[0,T] 
\end{array} \right. 
\end{equation}
with initial condition $u(\boldsymbol{x},0)=u_0(\boldsymbol{x})$. Typically, homogeneous Neumann boundary conditions are prescribed to model an isolated tissue.
Here $\Omega \subset \mathbb{R}^3$ is a spatial domain denoting the portion of cardiac tissue of interest  and $[0, T]$ is a fixed time interval. 
The symbol $C_{\rm m}$ is the membrane capacitance per unit area with $\beta = 2000 \mbox{ cm}^{-1}$ being the surface-to-volume ratio of the membrane; $\boldsymbol{\sigma}_{\rm i}$ ($ \boldsymbol{\sigma}_{\rm e}$) is the intracellular (extracellular) conductivity tensor;  $I_{\rm si}$ ($I_{\rm se}$) represents the intracellular  (extracellular) stimulation current. An explicit form of the total ionic current $I_{\rm ion}$ is described by an ionic model in which the {\it gating variable} $w$ is used to control the depolarization and repolarization phases of the cardiac action potential. The time evolution of the gating variable $w$ is generally modeled in the form 
$$\dfrac{d w}{d t} + g(u, w) = 0 \quad\mbox{ in } \Omega\times[0,T]$$
with initial condition $w(\boldsymbol{x},0)=w_0(\boldsymbol{x})$.

We can represent the conductivity tensors as
$\boldsymbol{\sigma}_{\rm k}(\boldsymbol{x})=
{\sigma}_{\rm kl}\mathbf{a_l}(\boldsymbol{x})\mathbf{a_l}(\boldsymbol{x})^T + {\sigma}_{\rm kt}\mathbf{a_t}(\boldsymbol{x})\mathbf{a_t}(\boldsymbol{x})^T 
+ {\sigma}_{\rm kn}\mathbf{a_n}(\boldsymbol{x})\mathbf{a_n}(\boldsymbol{x})^T,  $
where k stands for i or e, $ (\mathbf{a_l}, \mathbf{a_t}, \mathbf{a_n} )$ are orthonormal vectors related to the structure of the myocardium with $\mathbf{a_l}$ parallel to the fibre direction of the myocardial tissue. We further assume that the tissue is axial isotropic (i.e.~$\sigma_{\rm kn} = \sigma_{\rm kt}$) and postulate the conductivity $[\sigma_{\rm kl}, \sigma_{\rm kt}]$ to be constant, as has been done by several groups \cite{Clerc1976,Roberts2} for conductivity estimation in experiments. 

The Bidomain model has been widely used due to its ability to reproduce cardiac phenomena \cite{Trayanova}. However, its numerical solution for clinical application requires high computational cost, since it is a degenerate system of PDEs and the mesh and time constraints are significant for simulating fast potential variation.
The {\it Monodomain model} as a heuristic approximation of the Bidomain model has been proposed to provide computational improvements. Its derivation \cite{nielsen2007optimal} is based upon a proportionality assumption 
$\boldsymbol{\sigma}_{\rm e} = \lambda \boldsymbol{\sigma}_{\rm i}$, where $\lambda$ is a constant. A formulation of the Monodomain model is then obtained, by denoting $\boldsymbol{\sigma}_{\rm m} = \frac{\lambda}{1+\lambda}\boldsymbol{\sigma}_{\rm i}$ and $I_{\rm app}=\frac{\lambda}{1+\lambda}I_{\rm si}+\frac{1}{1+\lambda}I_{\rm se}$, as
\begin{equation} \label{monoDerive}
\begin{array}{ll}
\beta C_{\rm m} \dfrac{\partial u}{\partial t}-\nabla \cdot (\boldsymbol{\boldsymbol{\sigma}}_{\rm m} \nabla u) + \beta I_{\rm ion} =I_{\rm app} &\quad\mbox{ in } \Omega\times[0,T].
\end{array} 
\end{equation}

Although the assumption on its derivation lacks  physiological  foundation, the Monodomain model has been intensively used in clinic-oriented simulations \cite{Villongco2014305, Barros2015} since it requires significantly less computational efforts than the Bidomain model. More importantly, a comparison between the Bidomain and Monodomain models in \cite{Bourgault2010} concluded that the discrepancy between the models at the continuous level may be quite small: of order 1\% or even below in terms of activation time relative error. 

\begin{figure}
\begin{center}
\includegraphics[scale=0.45]{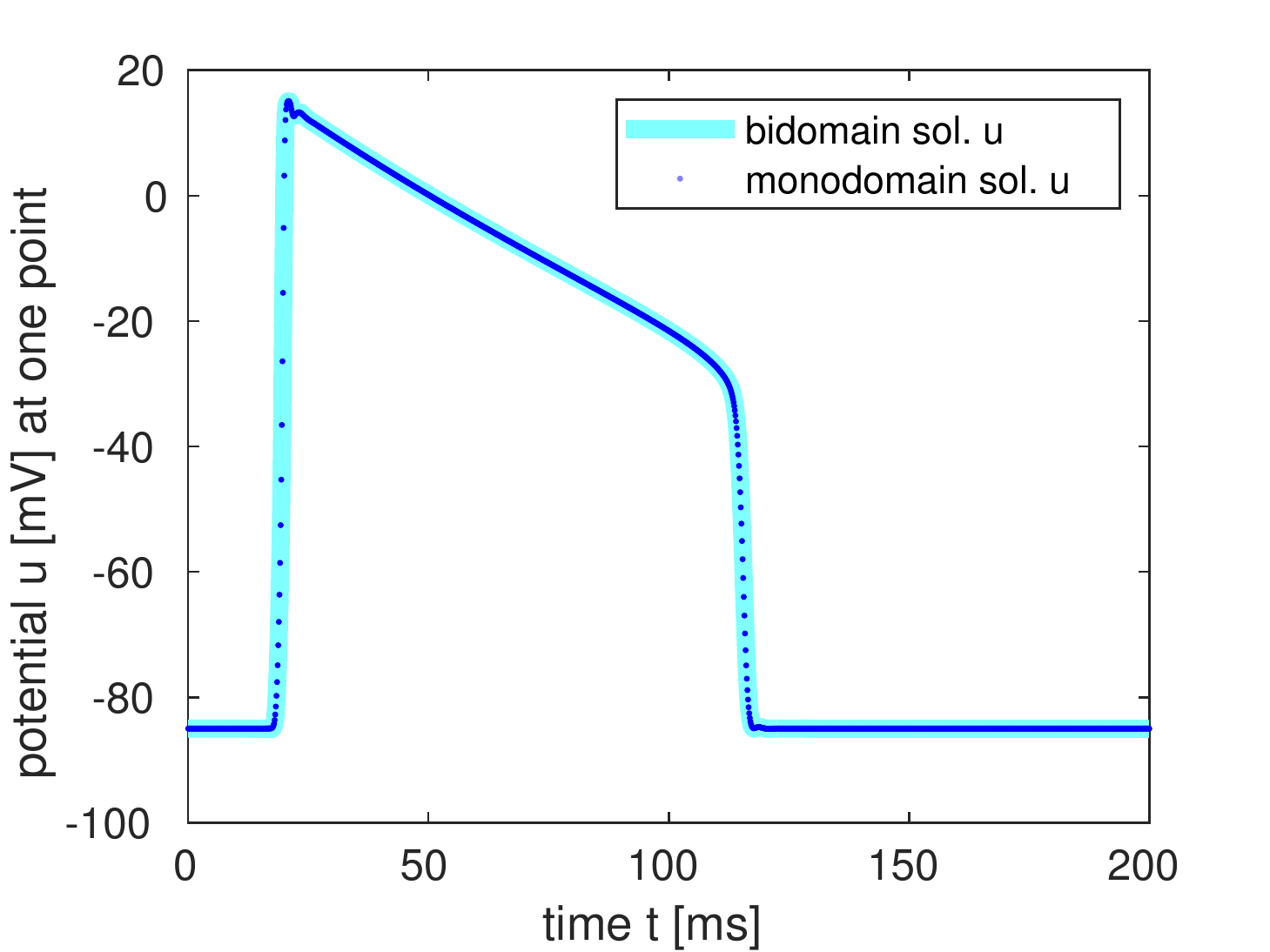}
\caption{The reconstruction of potential by the Monodomain solver, computed with the estimated conductivity $[\sigma_{\rm ml}, \sigma_{\rm mt}] = [1.704, 0.3551]$, gives an excellent matching with the Bidomain solution, which is computed with $[\sigma_{\rm il}, \sigma_{\rm el}, \sigma_{\rm it}, \sigma_{\rm et}] = [3.5, 3, 0.3, 1.8]$.}
\label{uAtPt}
\end{center}
\end{figure}

We resort to the Monodomain model in this paper for the cardiac conductivity estimation, following the ``potential oriented'' line. Namely, we speculate that an appropriate estimate of the conductivity tensor $\bsb{\sigma}_{\rm m}$ based on our variational procedure can still lead to an accurate reconstruction of the potential propagation. This is demonstrated by the numerical result shown in Fig.~\ref{uAtPt}, where the Bidomain solution $u$ computed on a slab mesh with $[\sigma_{\rm il}, \sigma_{\rm el}, \sigma_{\rm it}, \sigma_{\rm et}] = [3.5, 3, 0.3, 1.8]$ was used as synthetic measurement data to estimate the conductivity $\bsb{\sigma}_{\rm m}$ in the Monodomain solver. The reconstruction of potential by the Monodomain solver with the estimated conductivity $[\sigma_{\rm ml}, \sigma_{\rm mt}] = [1.704, 0.3551]$ gives an excellent matching with the Bidomain solution (for more details, see \cite{HHthesis}). 

In the following subsections, we describe the time and space discretization schemes for the full-order Monodomain model, and state the Monodomain inverse conductivity problem in a discrete version. 
As a starting point of studies on reduced-order modeling for conductivity estimation, here we simply use the classical Rogers--McCulloch \cite{rogers} ionic model, given by
\begin{equation}
\begin{array}{l}
I_{\rm ion}(u,w) = C_{\rm m}[\beta_1(u-V_{\rm r})(u-V_{\rm th})(u-V_{\rm p})+c_2(u-V_{\rm r})w ]\label{RMmodel}   \\
g(u,w) = -\beta_2(u-V_{\rm r})+bdw 
\end{array}
\end{equation}
with $\beta_1 = \frac{c_1}{(V_{\rm p}-V_{\rm r})^2}, ~ \beta_2=\frac{b}{V_{\rm p}-V_{\rm r}}.$
The values of ionic model parameters are taken from \cite{colli2004}: $C_{\rm m} = 1~\mu\mbox{Fcm}^{-2}$, $V_{\rm r}=-85$ mV, $V_{\rm th}=-72$ mV,  $V_{\rm p}=15$ mV,  $c_1=11.54 \mbox{ ms}^{-1}$, $c_2=4.4 \mbox{ ms}^{-1}$, $b = 0.012 \mbox{ ms}^{-1}$, $d = 1$.

\subsection*{Time and space discretization} 

To improve the computational efficiency, we split the Monodomain problem into a PDE and an ODE representing the Rogers--McCulloch model. The ODE gating variables are integrated with an implicit backward Euler method. The PDE is solved by a semi-implicit method based on a backward differentiation formula (BDF). The nonlinear terms are tackled by an explicit second-order time extrapolation \cite{HH2015IP}.
Let $\Delta t$ be the time step, we define $L=T/\Delta t, ~t^l=l\Delta t$. Hereafter we use superscripts $l$ and $l+1$ for those variables at time $t^l$ and $t^{l+1}$, respectively.

The finite element method (FEM) is used for space discretization. Let $\{\phi_j\}_{j=1}^n$ be the finite element basis, we denote by $\mathbf{M}$ the mass matrix 
and by $ \mathbf{S_{\rm l}}$ and $ \mathbf{S_{\rm t}}$ the stiffness matrices with entries
$$[\mathbf{S_{\rm l}}]_{jk}=\int_\Omega \mathbf{a}_{\rm l}\mathbf{a}_{\rm l}^T\nabla\phi_k\cdot\nabla\phi_j d\boldsymbol{x}, \quad [\mathbf{S_{\rm t}}]_{jk}=\int_\Omega (\mathbf{I}-\mathbf{a}_{\rm l}\mathbf{a}_{\rm l}^T)\nabla\phi_k\cdot\nabla\phi_j d\boldsymbol{x}.$$
The bold symbol $\mathbf{u}^l$ will denote the vector representation of $u(x,t^l)$ in the finite element space. We adopt a similar notation for the other variables.

At time $t=t^{l+1}$, the gating variable $\mathbf{w}^{l+1}$ in the ionic model is updated by
\begin{equation}\label{ion-algebra}
\dfrac{\mathbf{w}^{l+1}-\mathbf{w}^l}{\Delta t}=-g(\tilde{\mathbf{u}}^{l+1}, \mathbf{w}^{l+1})
\end{equation}
with $\tilde{\mathbf{u}}^{l+1} = 2\mathbf{u}^{l}-\mathbf{u}^{l-1}$ ($\tilde{\mathbf{u}}^{1}=\mathbf{u}^{0}$ in particular) being the second-order time extrapolation of $\mathbf{u}^{l+1}$. Here $g(\cdot)$ is evaluated component-wise, its $i$-th entry is computed as $g([\tilde{\mathbf{u}}^{l+1}]_i, [\mathbf{w}^{l+1}]_i)$.
For the sake of computational efficiency (see Sec.~\ref{pod-DEIM-resultInv}: The measures) and for the convenience of model reduction (Sec.~\ref{DEIM}: (\ref{fappx})), we approximate the nonlinear function $I_{\rm ion}$ and the stimulus $I_{\rm app}$ piecewise linearly. Under this assumption, the discretized Monodomain system reads
\begin{equation}\label{mono-algebra}
\mathbf{A_{\rm m}} \mathbf{u}^{l+1} = \mathbf{b}^{l+1}
\end{equation}
where $\mathbf{A_{\rm m}} = \beta C_{\rm m}\frac{\alpha_0}{\Delta t}\mathbf{M}+\sigma_{\rm ml}\mathbf{S_{\rm l}} + \sigma_{\rm mt}\mathbf{S_{\rm t}}$ and the right-hand side is evaluated as 
\begin{equation}\label{mono-b-ptws}
\mathbf{b}^{l+1} = \mathbf{M}(\mathbf{I}_{\rm app}^{l+1} - \beta I_{\rm ion}(\tilde{\mathbf{u}}^{l+1}, \mathbf{w}^{l+1})) + \beta C_{\rm m}\mathbf{M} \sum\limits_{i=1}^{2}\frac{\alpha_i}{\Delta t}\mathbf{u}^{l+1-i}.
\end{equation}
Here $\alpha_i$'s are the coefficients of the BDF (for a BDF of order two we have $\alpha_0 = 3/2,~\alpha_1 = 2,~\alpha_2 = -1/2$). In (\ref{mono-b-ptws}), $I_{\rm ion}(\tilde{\mathbf{u}}^{l+1}, \mathbf{w}^{l+1})$ is a component-wise evaluation. As we will see in Sec.~\ref{pod-DEIM-resultInv}, this component-wise approximation improves computational efficiency without significant 
loss of accuracy as compared with the exact finite element approximation. 

The linear system (\ref{mono-algebra}) can be solved by an ILU preconditioned conjugate gradient method implemented for instance in the Trilinos package ({\tt{www.trilinos.org}}). 
A mass lumping technique (\cite{Zienkiewicz2005}, Sec.~16.2.4) is employed for the sake of computational efficiency and 
stability. 
Specifically, we apply a diagonal scaling with the factor being the total mass $\mathbf{M}_{\rm tot}$ ($=\sum_{i,j}\mathbf{M}_{ij}$) divided by the trace, i.e.~the diagonal lumped mass matrix $\mathbf{M}_{\rm L}$ satisfies $[\mathbf{M}_{\rm L}]_{ii}  = \dfrac{\mathbf{M}_{\rm tot}}{\mbox{Tr}(\mathbf{M})}\mathbf{M}_{ii}$. 
For simplicity of notation, we will use the same symbol $\mathbf{M}$ for the lumped mass matrix as for the mass.
 We solve the Monodomain and Bidomain models using LifeV, an object oriented C++ finite element library developed by different groups worldwide \cite{lifev}.

\subsection*{The inverse conductivity problem} 
The {\it Monodomain inverse conductivity problem} (MICP) reads: find $\boldsymbol{\sigma} = [\sigma_{\rm ml}, \sigma_{\rm mt}]  \in \mathcal{C}_{\rm ad} \subset \mathbb{R}^2$ that minimizes 
\begin{equation}\label{}
 \mathcal{J}_{\rm m}(\boldsymbol{\sigma}) = \dfrac{1}{2}\sum\limits_{l=1}^L(\mathbf{u}^l-\mathbf{u}_{\rm meas}^l)^T \mathbb{X}_{\rm site} (\mathbf{u}^l-\mathbf{u}_{\rm meas}^l) \chi^l_{\rm snap} + \dfrac{\alpha}{2}\mathcal{R}(\boldsymbol{\sigma})
\end{equation}
subject to the Monodomain system (\ref{mono-algebra}) coupled with (\ref{ion-algebra}). Here $\mathcal{C}_{\rm ad}$ is an admissible domain given by inequality constraints on the conductivity, 
in the general form $\mathbf{h}(\boldsymbol{\sigma})\geq \boldsymbol{0}$.  The term $\mathbf{u}_{\rm meas}^l$ denotes experimental data at time $t^l$.
The data can be obtained in vitro using voltage optical mapping \cite{fentonMeasure}, or in vivo by back-mapping body surface potentials \cite{Cluitmans2013} or possibly by potential reconstruction from
electrocardiogram phase analysis of standard gated SPECT \cite{Chen2008JNC}.
We assume that the measurement sites are always grid points and
 $\mathbb{X}_{\rm site} \in \mathbb{R}^{n\times n}$ is the matrix recording observation sites.
The off-diagonal entries in $\mathbb{X}_{\rm site}$ are zeros; in the diagonal, $[\mathbb{X}_{\rm site}]_{ii} = 1$ if the spatial grid $x_i$ is an observation site and $0$ otherwise. We also introduce the {\it snapshot marker} $\chi^l_{\rm snap}$ which equals 1 if $t^l$ is an observation moment and $0$ otherwise.

The MICP can be solved by the BFGS optimization \cite{nocedal2006} approach as done in \cite{HH2015IP}. However, the computation is intensive since each ``query'' of the forward system or its adjoint counterpart is performed in the full-order $n$ by solving large algebraic systems. 
Moreover, in practice one may find also pathological tissues where scars inside have different (anomalous) values of conductivities. In this case the number of ``queries'' will substantially increase due to the increase of the total number of conductivities (see \cite{HH2015IP, HHthesis}). 
Altogether, in the current scenario, the MICP is too computationally expensive to be applied in clinic. This motivates us to work on a model reduction investigation with the aim of reducing significantly the computational cost.

\section{Model order reduction for nonlinear systems}\label{DEIMintr}
In this section, we describe the method in general by considering the discrete system of a parameterized nonlinear differential equation
\begin{equation}\label{generalSys}
\mathbf{A}\mathbf{y}(\tau) = \mathbf{F}(\mathbf{y}(\tau);\tau)
\end{equation}
where $\tau$ denotes the model parameter of interest in an inverse problem and is contained in a closed bounded domain $\mathcal{D}$. 
We assume that the matrix $\mathbf{A} \in \mathbb{R}^{n\times n}$ has affine dependence on $\tau$ and omit the explicit dependency $\mathbf{A}(\tau)$ for simplicity; $\mathbf{F}(\cdot)$ is a nonlinear function evaluated at the solution $\mathbf{y}(\tau) = [y_1(\tau), \cdots, y_n(\tau)]^T$ component-wise, i.e.  
$\mathbf{F}(\mathbf{y}(\tau);\tau) = [F(y_1(\tau);\tau), \cdots, F(y_n(\tau);\tau)]^T$
 with $F(\cdot)$ being a nonlinear scalar-valued function. In a finite element discretization, this is equivalent to a piece-wise linear approximation.

Our goal is to construct a small set of basis functions $\{\bsb{\varphi}_i \}_{i=1}^N$ in $\mathbb{R}^n$
such that the solution $\mathbf{y}(\tau)$ can be well approximated in the space $\spn\{\bsb{\varphi}_i\}$ (called {\it reduced space}). The functions $\{\bsb{\varphi}_i \}_{i=1}^N$ form the so-called {\it reduced basis} (RB).
Let $\mathbb{Z}_{\rm y} = [\bsb{\varphi}_1, \cdots, \bsb{\varphi}_N] \in \mathbb{R}^{n\times N}$, we represent the solution in the reduced space as
$
\mathbf{y}(\tau) = \mathbb{Z}_{\rm y}\mathbf{y_r}(\tau)
$
with $\mathbf{y_r}(\tau)$ being the vector of coordinates in the reduced space.
A {\it reduced-order model} (ROM) of (\ref{generalSys}) is then obtained by Galerkin projection
\begin{equation}\label{generalReduced}
\underbrace{\mathbb{Z}_{\rm y}^T\mathbf{A}\mathbb{Z}_{\rm y}}_{\mathbf{A_r}}\mathbf{y_r}(\tau) = \mathbb{Z}_{\rm y}^T\mathbf{F}(\mathbb{Z}_{\rm y}\mathbf{y_r}(\tau); \tau).
\end{equation}
Notice that $\mathbf{A_r} = \mathbb{Z}_{\rm y}^T\mathbf{A}\mathbb{Z}_{\rm y} \in \mathbb{R}^{N\times N}$ is a dense matrix but in general it features a very small size, hence the linear system (\ref{generalReduced}) can be tackled with a direct solver. 

The RB is constructed from the {\it full-order model} (FOM) (\ref{generalSys}) which is of large scale, thus the computation is usually expensive and performed offline. 
In the online phase, the ROM (\ref{generalReduced}) is solved many times for different parameter values with remarkably lower computational costs than the FOM.
Techniques for RB construction may rely on a sampling on the parameter $\tau$. We introduce here a sample $S = \{\tau_1, \cdots, \tau_s \}$
consisting of $s$ distinct parameter points in $\mathcal{D}$. 
A RB is constructed so to guarantee that for each $\tau_i \in S$ the error of approximating $\mathbf{y}(\tau_i)$ in the reduced space is bounded by a desired tolerance. Hereafter we follow the POD approach, which usually constructs an ``optimal" reduced basis as specified below.

\subsection{Proper Orthogonal Decomposition (POD)}\label{pod}
For the sake of completeness we briefly recall basic features of POD. More details can be found e.g.~in \cite{Kunisch2002}.
Given the parameter sample $S$, we solve the FOM (\ref{generalSys}) for each parameter value in $S$. The solutions are called {\it snapshots} and denoted by $\{\mathbf{y}_i \}_{i=1}^m$ ($m=s$ in current setting). The Proper Orthogonal Decomposition (POD) approach seeks an orthonormal {\it POD basis} $\{\boldsymbol{\varphi}_1, \cdots, \boldsymbol{\varphi}_N \}$ (also known as a set of {\it POD modes}) in $\mathbb{R}^n$ of a given rank $N~(N\ll m)$ that can best approximate the training space $X^{trn} = \spn\{\mathbf{y}_i \}_{i=1}^m$. Here ``best'' means the POD basis solves 
\begin{equation}\label{minPOD}
\min_{\{\boldsymbol{\psi}_i\}} \sum_{j=1}^m ||\mathbf{y}_j-\sum_{i=1}^N\langle \mathbf{y}_j,\boldsymbol{\psi}_i \rangle\boldsymbol{\psi}_i ||^2
 \qquad\mbox{ s.t. } \langle\boldsymbol{\psi}_i,\boldsymbol{\psi}_j\rangle = \delta_{ij}.
\end{equation}
We gather the snapshots into the so called {\it snapshot matrix} $\mathbf{Y}=[\mathbf{y}_1, \cdots, \mathbf{y}_m] \in \mathbb{R}^{n\times m}$.  The POD modes are given by the $N$ left singular vectors of $\mathbf{Y}$ associated with the $N$ largest singular values \cite{Kunisch2002}. Without loss of generality, we assume that the snapshot mean is zero. In fact we can replace each $\mathbf{y}_i$ with $\mathbf{y}_i-\frac{1}{m}\sum\limits_{j=1}^m\mathbf{y}_j$.

An efficient way for computing the POD modes through snapshot matrix $\mathbf{Y}$
is to first compute the {\it thin QR factorization} of $\mathbf{Y}$ as $\mathbf{Y} = \mathbf{QR}$ , and then compute the singular value decomposition of matrix $\mathbf{R}\in\mathbb{R}^{m\times m}$ as $\mathbf{R}=\mathbf{U}_{\rm R}\mathbf{S}_{\rm R}\mathbf{V}_{\rm R}^T$. The POD modes can be extracted in order as the first $N$ columns of $\mathbf{QU}_{\rm R}$.

\subsubsection*{Snapshots selection / Sampling}
The effectiveness of model reduction is clearly related to the representativity of the snapshots.
Standard schemes of sampling in the parameter space, which determines snapshots selection, include uniform sampling, random sampling, the Latin hypercube sampling (LHS, \cite{McKay1979}), and the centroidal voronoi tessellation (CVT) sampling \cite{Du1999}. When sampling a high-dimensional parameter space, the greedy sampling method could be used \cite{quarteroni2011certified, Nguyen2010}. The key feature of greedy sampling is to adaptively select a parameter at which the estimate of the solution error in the ROM is maximal. In this way,
we select the most effective parameter for controlling the error of the ROM.
However, the application of this method is limited to problems where sharp error estimators are available.

For optimal control applications, several online adaptive sampling procedures have been proposed to let the sampling procedure take into account the optimization trajectory. The Trust Region POD \cite{Arian00trust} approach constructs successively improved POD bases according to parameter values updated during optimization. The updating procedure was embedded with the trust region method which determines whether after an optimization step the POD basis should be updated. 
The Compact POD \cite{Carlberg2008} uses snapshots of the full-order model solutions as well as their derivatives (known as {\it sensitivities}) with respect to the model parameters of interest.

\subsection{Discrete Empirical Interpolation Method (DEIM)}\label{DEIM}

In the reduced system (\ref{generalReduced}) obtained by the POD projection, we have to evaluate the nonlinear term
\begin{equation}\label{nonlinear-n}
\mathbf{n}(\tau) = \underbrace{\mathbb{Z}_{\rm y}^T}_{N\times n}\underbrace{\mathbf{F}(\mathbb{Z}_{\rm y}\mathbf{y_r}(\tau);\tau)}_{n\times 1}.
\end{equation}
This evaluation has computational complexity depending on the size $n$ of the FOM (\ref{generalSys}), which is possibly in the magnitude of hundred thousand. Therefore, solving the ROM (\ref{generalReduced}) without an appropriate methodology may be as expensive as solving the full one. 

The complexity of evaluating (\ref{nonlinear-n}) can be made independent of the full order $n$ by using the Discrete Empirical Interpolation Method (DEIM \cite{thesisDEIM}). The DEIM provides an interpolation approximation for the nonlinear term $\mathbf{F}(\mathbb{Z}_{\rm y}\mathbf{y_r}(\tau);\tau)$ (simply denoted as $\mathbf{f}(\tau)$ in the sequel) by a projection onto a low-dimensional subspace. For this purpose, we introduce an $M$ dimensional space ($M \ll n$) where we look for an approximation of $\mathbf{f}(\tau)$ for values $\tau$ of interest. In particular, we can sample $\tau$, take snapshots of $\mathbf{f}(\tau)$ computed from the FOM with those samples, and then apply the POD on the snapshots to extract a projection basis $\{\mathbf{z}_1,\cdots,\mathbf{z}_M \}$.

Let $\mathbb{Z}_{\rm f} = [\mathbf{z}_1,\cdots,\mathbf{z}_M] \in \mathbb{R}^{n\times M}$, the DEIM approximation of $\mathbf{f}$ is in the form 
$\mathbf{f} \approx \hat{\mathbf{f}} =\mathbb{Z}_{\rm f} \mathbf{c}.$
To determine the coefficient vector $\mathbf{c}$, the DEIM optimally extracts $M$ distinct rows from the over-determined system $\mathbf{f}=\mathbb{Z}_{\rm f} \mathbf{c}$. Specifically, the DEIM selects row indices $p_1,\cdots, p_M$ in $\{1, ... ,n\}$ and requires:
$[\mathbf{f}]_{p_i} = [\hat{\mathbf{f}}]_{p_i}.$
If we denote ${\rm P} = [\mathbf{e}_{p_1}, \cdots, \mathbf{e}_{p_M}] \in \mathbb{R}^{n\times M} $ with $\mathbf{e}_{p_i}$ being the $p_i\mbox{-th}$ unit vector in $\mathbb{R}^n$, the  coefficient vector $\mathbf{c}$ is solved 
from ${\rm P}^T\mathbf{f}=({\rm P}^T\mathbb{Z}_{\rm f})\mathbf{c}$.
Finally the approximation of $\mathbf{f}$ writes
\begin{equation}\label{fappx}
\mathbf{f} \approx \hat{\mathbf{f}} = \mathbb{Z}_{\rm f} \mathbf{c} = \mathbb{Z}_{\rm f}({\rm P}^T \mathbb{Z}_{\rm f})^{-1} {\rm P}^T \mathbf{f} = \mathbb{Z}_{\rm f}({\rm P}^T \mathbb{Z}_{\rm f})^{-1}\mathbf{F}({\rm P}^T\mathbb{Z}_{\rm y}\mathbf{y_r}(\tau);\tau).
\end{equation}
The last equality in (\ref{fappx}) follows from the assumption that the function $\mathbf{F}(\cdot)$ evaluates component-wise  at its input vector. The nonlinear term (\ref{nonlinear-n}) can then be efficiently computed through
\begin{equation}
\mathbf{n}(\tau) \approx \underbrace{\mathbb{Z}_{\rm y}^T\mathbb{Z}_{\rm f}({\rm P}^T \mathbb{Z}_{\rm f})^{-1}}_{N\times M}\underbrace{\mathbf{F}({\rm P}^T\mathbb{Z}_{\rm y}\mathbf{y_r}(\tau);\tau)}_{M \times 1}.
\end{equation}
Notice that the matrices $\mathbb{Z}_{\rm y}^T\mathbb{Z}_{\rm f}({\rm P}^T \mathbb{Z}_{\rm f})^{-1} \in \mathbb{R}^{N\times M}$  and ${\rm P}^T\mathbb{Z}_{\rm y} \in \mathbb{R}^{M\times N}$ 
can be precomputed so that the computational complexity of evaluating $\mathbf{n}(\tau)$ is only $\mathcal{O}(MN)$.

The interpolation indices ${p_1}, \cdots, {p_M}$ are selected inductively from the projection basis $\{\mathbf{z}_1 , \cdots, \mathbf{z}_M \}$ by the DEIM algorithm. For the sake of completeness, 
we recall in Algorithm \ref{DEIM_alg} the DEIM described in \cite{thesisDEIM}. At each iteration, an interpolation index is selected to limit growth of the error bound of the approximation $\hat{\mathbf{f}}$. In particular, the first index $p_1$ is the index on which $\mathbf{z}_1$ has the largest magnitude; each of the remaining indices $p_l$ is the index on which the residual of approximating $\mathbf{z}_l$ by the first $l-1$ basis vectors $\{\mathbf{z}_1,\cdots,\mathbf{z}_{l-1} \}$ has the largest magnitude.
It is demonstrated that the DEIM algorithm is well-defined (\cite{thesisDEIM}, lem.~2.2.2). In fact, ${\rm P}^T\mathbf{Z}_{\rm f}$ is non-singular in each iteration of Algorithm \ref{DEIM_alg} and the interpolation indices are not repeated. 

\begin{algorithm}[tp]
\caption{DEIM \cite{thesisDEIM}}\label{DEIM_alg}
\begin{algorithmic}[1]
\Input $\{\mathbf{z}_l\}_{l=1}^M \subseteq \mathbb{R}^n $ linear independent
\Output indices $\mathbf{p} = [p_1, \cdots, p_M]^T \in \mathbb{R}^M$,  $\mathbf{M}$ as the inverse of ${\rm P}^T\mathbb{Z}_{\rm f}$
\State $[|\rho|, p_1] = \max\{|\mathbf{z}_1|\} $
\State $\mathbb{Z}_{\rm f}\leftarrow[\mathbf{z}_1],~ {\rm P}\leftarrow[\mathbf{e}_{p_1}], ~\mathbf{p} \leftarrow [p_1]$
\For{$l=2, \cdots, M$}
	\State Solve $ ({\rm P}^T\mathbb{Z}_{\rm f})\mathbf{c} = {\rm P}^T\mathbf{z}_l $ 
	\State $\mathbf{r}=\mathbf{z}_l - \mathbb{Z}_{\rm f}\mathbf{c}$
	\State $[|\rho|, p_l] = \max\{|\mathbf{r}|\} $
	\State $\mathbf{a}^T = \mathbf{e}_{p_l}^T \mathbb{Z}_{\rm f}$
	\State $$\mathbf{M} \leftarrow 
\left[ 
\begin{array}{cc}
\mathbf{I} & -\mathbf{c} \\
\boldsymbol{0} & 1
\end{array}
\right]
\left[ 
\begin{array}{cc}
 \mathbf{M} & \boldsymbol{0} \\
-\rho^{-1}\mathbf{a}^T \mathbf{M}  & \rho^{-1}
\end{array}
\right] $$
	\State $\mathbb{Z}_{\rm f}\leftarrow [\mathbb{Z}_{\rm f} ~\mathbf{z}_l],~ {\rm P}\leftarrow [{\rm P} ~\mathbf{e}_{p_l}], ~\mathbf{p} \leftarrow [\mathbf{p}^T ~p_l]^T$
\EndFor
\end{algorithmic}
\end{algorithm}

\section{The reduced-order Monodomain inverse conductivity problem}\label{DEIM-use}
The combination of POD and DEIM as described can be applied to the MICP. However, the real accuracy (and also efficiency) of the procedure depends on the nature of the problem and ultimately on the snapshots selection. In this section we present the standard tools to apply POD-DEIM to the MICP. The performance and the specific customization needed to make the procedure effective are discussed in Sec. \ref{DEIMresult}.

\subsection*{The POD basis for the Monodomain model}
In constructing a reduced basis for the parameterized Monodomain model, two ways can be followed as mentioned in \cite{Boulakia2011POD}. One is to store a set of solutions of the full-order model computed at different instants with a set of different parameter values in a given sample space, then collect all these solutions to build a unique snapshot matrix, upon which the POD basis is finally built. In other words, we treat both the conductivity tensor $\boldsymbol{\sigma}$ and also the time $t$ as parameters of the model. Another way for the reduced basis construction is to build multiple POD bases instead of a unique one. Each POD basis is constructed from the snapshots of solutions computed with a particular conductivity parameter, which we call {\it generating parameter} and denote by $\bsb{\sigma}_{\rm gen}$. 
The idea is that for a given value of the tensor $\bsb{\sigma}$ we select the POD basis obtained with the closest generating parameter (among all available) to approximate the Monodomain problem.

\begin{figure}[!tp]
\begin{center}
\begin{minipage}{0.49\textwidth}
\includegraphics[scale=0.35]{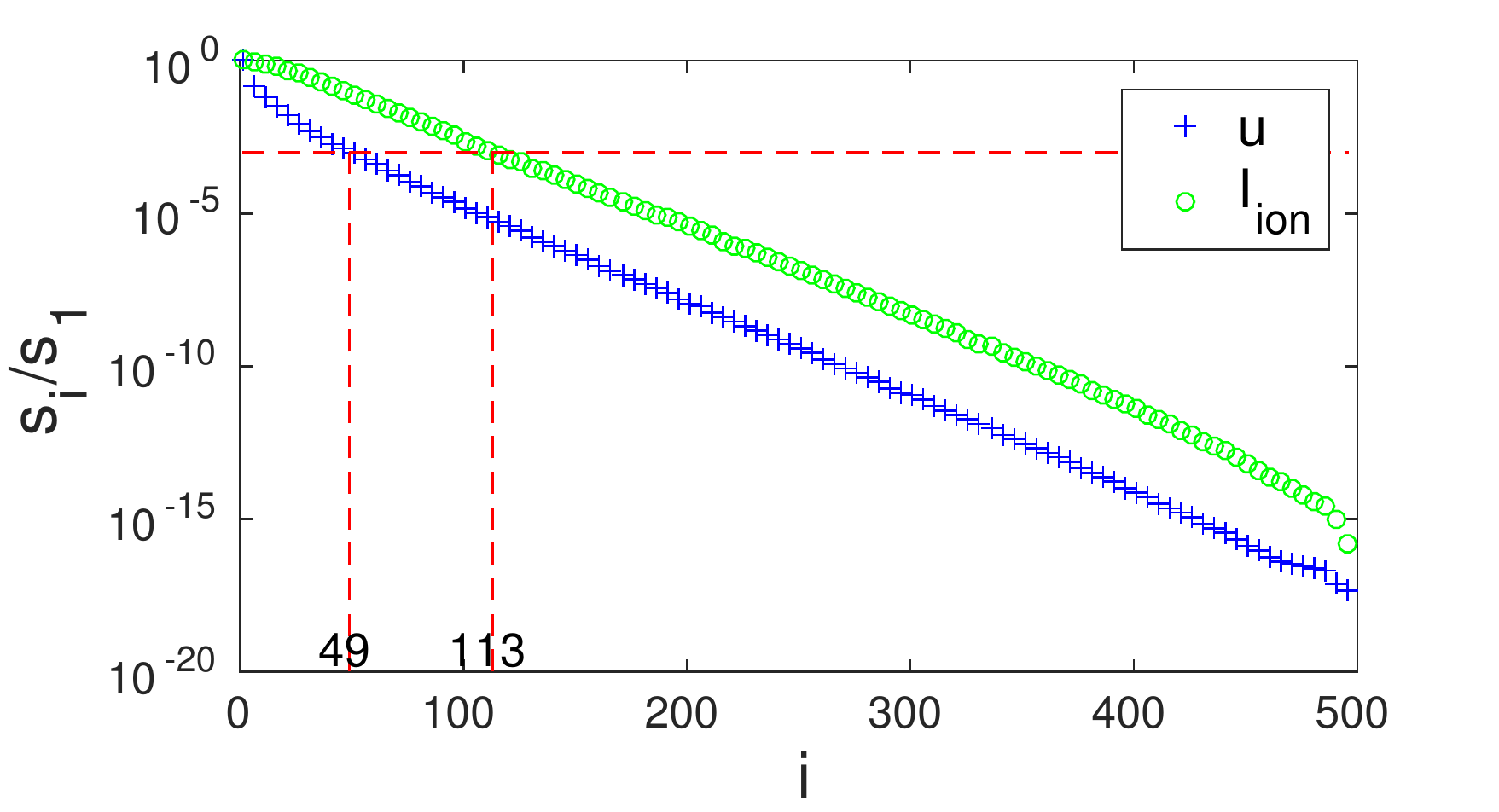}
\end{minipage}
\begin{minipage}{0.49\textwidth}
\includegraphics[scale=0.35]{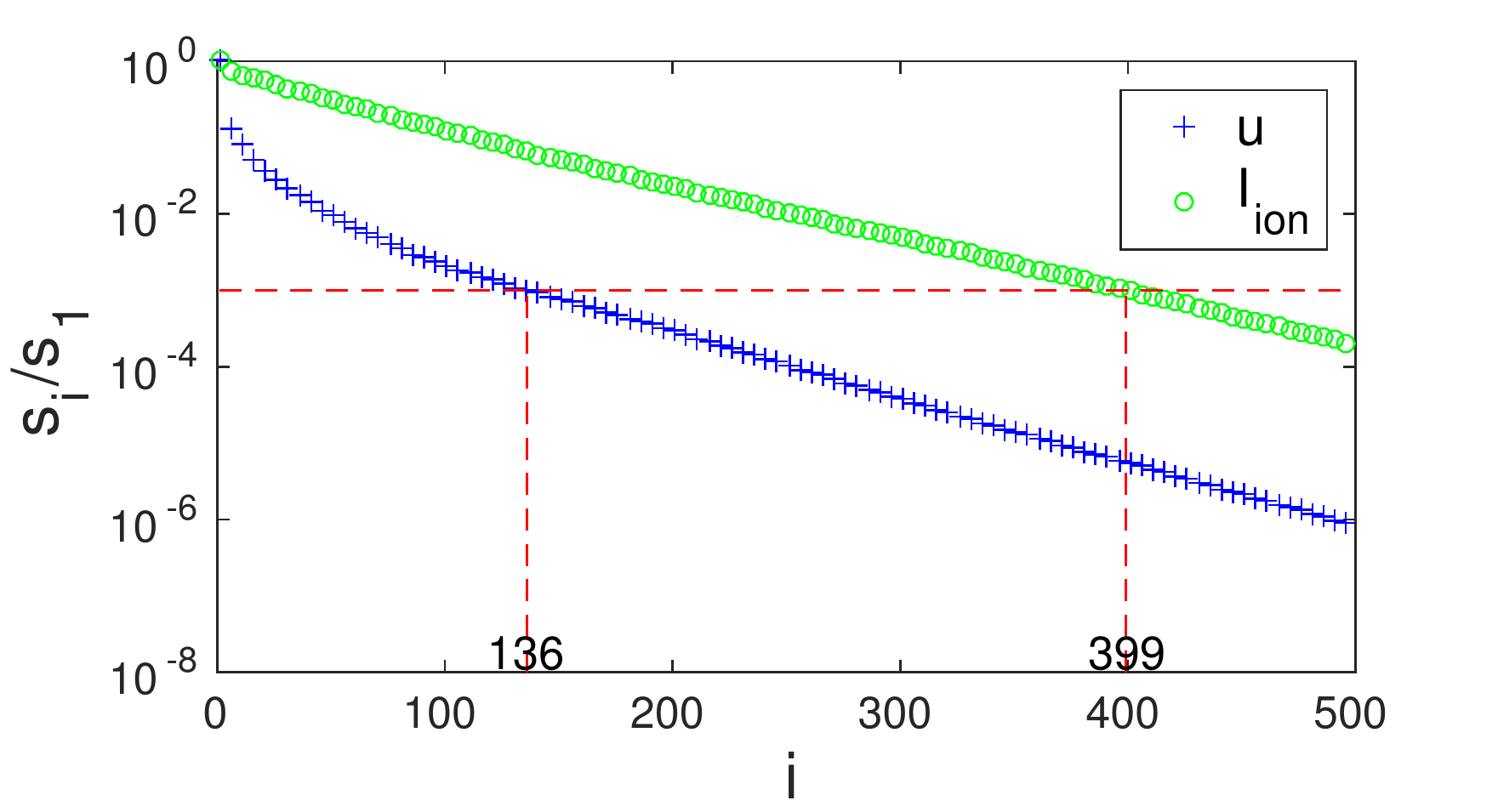}
\end{minipage}
\caption{Singular values of snapshot matrices of $u$ and $I_{\rm ion}$. Left: snapshots generated by one parameter value. Right: snapshots generated by four parameter values.}
\label{svd-onePara-multi}
\end{center}
\end{figure}

The dynamic of the transmembrane potential $u$ is not quite smooth in time, due to the wavefront propagation or the upstroke (depolarization) spreading. This has an important impact on POD procedures. In fact, the singular values of the snapshot matrix of $u$ do not decay fast in general. We confirm this by a numerical experiment shown in Fig.~\ref{svd-onePara-multi}~(left). In this test, the time step $\Delta t = 0.05$~ms and the full-order model dimension is 24272. We collected snapshots of the transmembrane potential $u$ (and the ionic current $I_{\rm ion}$) computed with a fixed conductivity parameter $\boldsymbol{\sigma}=[3, 1]$ for 500 time steps. Slow decay shown in Fig.~\ref{svd-onePara-multi}~(left) is apparent especially when compared with other problems, in which singular values of a snapshot matrix usually decrease fast hence few POD modes are enough to give an accurate approximation of the solution considered (see e.g.~\cite{BertagnaVeneziani} for a FLuid-Structure Interaction problem). 
The slow decay is detrimental to the actual model reduction as many modes need to be considered in the ROM construction.

Fig.~\ref{svd-onePara-multi} (left) shows that the singular value ${s_i}$ of snapshots of the transmembrane potential relative to the leading singular value ${s_1}$ decays to $10^{-3}$ (cross markers 
in Fig.~\ref{svd-onePara-multi}) when $i$ is 49 (113 for the ionic current). 
It suggests that $u$ ($I_{\rm ion}$) may still be well approximated by a POD basis of dimension less than 50 (120), provided that the model parameter is close enough to the generating parameter of the POD basis. We also infer that the ionic current $I_{\rm ion}$ features even more complex nonlinearity than the transmembrane potential $u$, since more POD modes are needed for approximating the ionic current with the same accuracy of the transmembrane potential (circular markers in Fig.~\ref{svd-onePara-multi})

A consequence of this nature of the problem is that building a unique POD basis actually provides worthless results. This can be identified from Fig.~\ref{svd-onePara-multi} right, in which the leading 500 singular values of snapshots generated by four different conductivity parameters are reported.
The relative singular value $s_i/s_1$ of snapshots of the transmembrane $u$ decays to $10^{-3}$  when $i$ is 136 (399 for the ionic current). Compared with the one-parameter case shown in Fig.~\ref{svd-onePara-multi} left, a combination of snapshots from different parameter values does not reduce the number of POD modes necessary for accurate POD approximation.
This feature demands for specific customization of the POD procedure based on an educated selection of the snapshots as detailed in Sec.~\ref{DEIMresult}.

\subsection*{Sensitivity equations}
One possible way to improve our model reduction is to add snapshots of the sensitivity of the Monodomain solution to the conductivity tensor. To this aim, we explicitly report the sensitivity equations.

Applying differentiation with respect to $\sigma_{\rm mk}$ on the Monodomain model ((\ref{mono-algebra}) coupled with (\ref{ion-algebra})), at time $t^{l+1}$ we can solve $\frac{\partial \mathbf{u}^{l+1}}{\partial \sigma_{\rm mk}}$ from
\begin{equation}
\frac{\partial\mathbf{A_{\rm m}}}{\partial \sigma_{\rm mk}} \mathbf{u}^{l+1} + \mathbf{A_{\rm m}}\frac{\partial \mathbf{u}^{l+1}}{\partial \sigma_{\rm mk}} = \frac{\partial \mathbf{b}^{l+1}}{\partial \sigma_{\rm mk}}.
\end{equation}
The sensitivity of the ionic current, stated in the continuous form without loss of generality, can be evaluated through
$$ \frac{\partial I_{\rm ion}}{\partial \sigma_{\rm mk}} = \partial_u I_{\rm ion}\frac{\partial u}{\partial \sigma_{\rm mk}} + \partial_w I_{\rm ion}\frac{\partial w}{\partial \sigma_{\rm mk}}$$
where ${\rm k}$ stands for either ${\rm l}$ or ${\rm t}$.

\subsection*{Application of DEIM on the Monodomain solver}
Let us denote the reduced bases for the transmembrane potential $u$ and the ionic current $I_{\rm ion}$ (evaluated in the degrees of freedom) by $\mathbb{Z}_{\rm u}$ and $\mathbb{Z}_{\rm ion}$ respectively. By projecting the discrete Monodomain system (\ref{mono-algebra}, \ref{mono-b-ptws}) onto the reduced space $\mathbb{Z}_{\rm u}$ , we obtain
\begin{equation}\label{reducedMono}
\left\{
\begin{array}{l}
\mathbb{Z}_{\rm u}^T \mathbf{A_{\rm m}}\mathbb{Z}_{\rm u}\mathbf{u}_\mathbf{r}^{l+1} = \mathbf{b}_\mathbf{r}^{l+1} 
\quad\quad l\in\{0, 1, \cdots, L-1 \} \cr
\mbox{with }\quad \mathbf{b}_\mathbf{r}^{l+1}  = \mathbb{Z}_{\rm u}^T\mathbf{b}^{l+1}
= \underbrace{\mathbb{Z}_{\rm u}^T\mathbf{M}\mathbf{I}_{\rm app}^{l+1}}_{N\times 1}
+ \beta C_{\rm m}\underbrace{\mathbb{Z}_{\rm u}^T\mathbf{M}\mathbb{Z}_{\rm u}}_{N\times N}\sum\limits_{i=1}^{2}\frac{\alpha_i}{\Delta t}\underbrace{\mathbf{u}_\mathbf{r}^{l+1-i}}_{N\times 1} \cr
\hspace{4cm} - \beta\underbrace{\mathbb{Z}_{\rm u}^T\mathbf{M}}_{N\times n}\underbrace{I_{\rm ion}(\mathbb{Z}_{\rm u}\widetilde{\mathbf{u}}_\mathbf{r}^{l+1}, \mathbf{w}^{l+1})}_{n\times 1}  
\end{array}
\right.
\end{equation}
 where $\mathbf{u}_\mathbf{r}^{l+1}$ is the solution at time $t=t^{l+1}$ in the reduced space and
 $\widetilde{\mathbf{u}}_\mathbf{r}^{l+1} = 2\mathbf{u}_\mathbf{r}^{l} - \mathbf{u}_\mathbf{r}^{l-1}$ is its second-order extrapolation.
 As discussed in Sec.~\ref{DEIM}, the complexity in computing the nonlinear term $I_{\rm ion}$ in the right hand side is $\mathcal{O}(n)$. This mainly results from the fact that $I_{\rm ion}(\mathbb{Z}_{\rm u}\widetilde{\mathbf{u}}_\mathbf{r}^{l+1}, \mathbf{w}^{l+1})$ cannot be precomputed, since it depends nonlinearly on the full-order vector $\mathbb{Z}_{\rm u}\widetilde{\mathbf{u}}_\mathbf{r}^{l+1}$ and $\mathbf{w}^{l+1}$. We will apply the DEIM approximation to $I_{\rm ion}(\mathbb{Z}_{\rm u}\widetilde{\mathbf{u}}_\mathbf{r}^{l+1}, \mathbf{w}^{l+1})$.
We use the POD basis $\mathbb{Z}_{\rm ion} \in \mathbb{R}^{n\times M}$ of snapshots of $I_{\rm ion}$ as an input basis for the DEIM algorithm, where $M$ is the number of POD modes. The DEIM algorithm generates interpolation indices $\mathbf{p} = [p_1, \cdots, p_M]^T $ for constructing the extraction matrix ${\rm P}$. The DEIM approximation reads then
\begin{equation}
I_{\rm ion}(\mathbb{Z}_{\rm u}\widetilde{\mathbf{u}}_\mathbf{r}^{l+1}, \mathbf{w}^{l+1})
\approx \underbrace{\mathbb{Z}_{\rm ion}({\rm P}^T\mathbb{Z}_{\rm ion})^{-1}}_{n\times M}I_{\rm ion}(\underbrace{{\rm P}^T\mathbb{Z}_{u}\tilde{\mathbf{u}}_\mathbf{r}^{l+1}}_{M\times 1},~ \underbrace{{\rm P}^T\mathbf{w}^{l+1}}_{M \times 1}).
\end{equation}
If we set
\begin{align*}
&\mathbf{A_r} = \mathbb{Z}_{\rm u}^T\mathbf{A_{\rm m}}\mathbb{Z}_{\rm u} \in \mathbb{R}^{N\times N},\\
&\mathbf{M}_{\rm u} = \mathbb{Z}_{\rm u}^T\mathbf{M}\mathbb{Z}_{\rm u} \in \mathbb{R}^{N\times N},
~\mathbf{M}_{\rm iu} = \underbrace{\mathbb{Z}_{\rm u}^T\mathbf{M}}_{N\times n}\underbrace{\mathbb{Z}_{\rm ion}({\rm P}^T\mathbb{Z}_{\rm ion})^{-1}}_{n\times M} \in \mathbb{R}^{N\times M}, \\
&\mathbf{I}_\mathbf{r}^{l+1} = \mathbb{Z}_{\rm u}^T\mathbf{M}\mathbf{I}_{\rm app}^{l+1} \in \mathbb{R}^{N\times 1} ,
~\mathbf{w}_\mathbf{r}^{l+1} = {\rm P}^T\mathbf{w}^{l+1} \in \mathbb{R}^{M\times 1},
~{\rm U} = {\rm P}^T\mathbb{Z}_{\rm u} \in \mathbb{R}^{M\times N} ,
\end{align*}
the reduced Monodomain system is then formulated as
\begin{equation}\label{rbDEIM_mono}
\left\{
\begin{array}{l}
\mathbf{A_r}\mathbf{u}_\mathbf{r}^{l+1} 
= \mathbf{I}_\mathbf{r}^{l+1}
+ \beta C_{\rm m}\mathbf{M}_{\rm u}\sum\limits_{i=1}^{2}\frac{\alpha_i}{\Delta t}\mathbf{u}_\mathbf{r}^{l+1-i}
- \beta \mathbf{M}_{\rm iu} I_{\rm ion}({\rm U}\widetilde{\mathbf{u}}_\mathbf{r}^{l+1},~\mathbf{w}_\mathbf{r}^{l+1}) \\
\dfrac{\mathbf{w}_\mathbf{r}^{l+1}-\mathbf{w}_\mathbf{r}^{l}}{\Delta t} = -g({\rm U}\widetilde{\mathbf{u}}_\mathbf{r}^{l+1}, \mathbf{w}_\mathbf{r}^{l+1})
\end{array}
\right..
\end{equation}

\subsection*{The reduced Monodomain inverse conductivity problem}
We are in position of formulating the reduced MICP: find $\boldsymbol{\sigma} = [\sigma_{\rm ml}, \sigma_{\rm mt}]$ that minimizes 
\begin{equation}\label{}
 \mathcal{J}_{\rm r}(\boldsymbol{\sigma}) = \dfrac{1}{2}\sum\limits_{l=1}^L(\mathbb{Z}_{\rm u}\mathbf{u}_\mathbf{r}^l-\mathbf{u}_{\rm meas}^l)^T \mathbb{X}_{\rm site} (\mathbb{Z}_{\rm u}\mathbf{u}_\mathbf{r}^l-\mathbf{u}_{\rm meas}^l) \chi^l_{\rm snap} + \dfrac{\alpha}{2}\mathcal{R}(\boldsymbol{\sigma})
\end{equation}
subject to the reduced Monodomain system (\ref{rbDEIM_mono}) and the inequality constraint $\mathbf{h}(\boldsymbol{\sigma})\geq \boldsymbol{0}$. 
In particular, we restrict the admissible domain of the conductivities as
$$\big\{ \sigma_{\rm ml}/\sigma_{\rm mt}\geq 1,~\sigma_{\rm ml}/\sigma_{\rm mt}\leq 100, ~\sigma_{\rm mt}\geq 0.05, ~\sigma_{\rm ml}\leq 7 \big\},$$
inspired from the range of conductivity measures listed in Table 1 of \cite{HH2015IP}. We will denote the components of $\mathbf{h}(\boldsymbol{\sigma})$ as $h_i(\boldsymbol{\sigma})$ in 
the algorithm description.
In order to obtain linear algebra operations of complexity $\mathcal{O}(N)$ only, in computing the cost function $\mathcal{J}_{\rm r}$, we can rewrite
\begin{equation}\label{}
\mathcal{J}_{\rm r}(\boldsymbol{\sigma}) = \dfrac{1}{2}\sum\limits_{l=1}^L
( (\mathbf{u}_\mathbf{r}^l)^T  \mathbf{X}_{\rm u}    \mathbf{u}_\mathbf{r}^l
-2(\widehat{\mathbf{u}}_{\rm meas}^l)^T\mathbf{u}_\mathbf{r}^l
+ ||\widehat{\vphantom{\rule{1pt}{6.5pt}}\smash{\widehat{\mathbf{u}}}}_{\rm meas}^l||^2
) \chi^l_{\rm snap} + \dfrac{\alpha}{2}\mathcal{R}(\boldsymbol{\sigma})
\end{equation}
where
$\mathbf{X}_{\rm u} =  \mathbb{Z}_{\rm u}^T\mathbb{X}_{\rm site}\mathbb{Z}_{\rm u} \in \mathbb{R}^{N\times N}$.
The projected measurements
$\widehat{\mathbf{u}}_{\rm meas}^l= \mathbb{Z}_{\rm u}^T\mathbb{X}_{\rm site}\mathbf{u}_{\rm meas}^l  \in \mathbb{R}^{N\times 1} \mbox{ and } \widehat{\vphantom{\rule{1pt}{6.5pt}}\smash{\widehat{\mathbf{u}}}}_{\rm meas}^l= \mathbb{X}_{\rm site}\mathbf{u}_{\rm meas}^l \in \mathbb{R}^{n\times 1}$
can be precomputed.

As we did in \cite{HH2015IP} for the full-order problem, we introduce the Lagrange multipliers $\{\mathbf{q}_\mathbf{r}^1,\cdots,\mathbf{q}_\mathbf{r}^L, \mathbf{r}_\mathbf{r}^1,\cdots,\mathbf{r}_\mathbf{r}^L\}$
and the Lagrangian functional
\begin{equation}
\begin{split}
&\mathcal{L}_{\rm r}(\mathbf{u}_\mathbf{r}^1,\cdots,\mathbf{u}_\mathbf{r}^L, \mathbf{w}_\mathbf{r}^1,\cdots,\mathbf{w}_\mathbf{r}^L, \boldsymbol{\sigma}, \mathbf{q}_\mathbf{r}^1,\cdots,\mathbf{q}_\mathbf{r}^L, \mathbf{r}_\mathbf{r}^1,\cdots,\mathbf{r}_\mathbf{r}^L) = \mathcal{J}_{\rm r}(\boldsymbol{\sigma})  \\
&\qquad\qquad-\sum\limits_{l=1}^L (\mathbf{q}_\mathbf{r}^{l})^T(\mathbf{A_r}\mathbf{u}_\mathbf{r}^{l} - \mathbf{I}_\mathbf{r}^{l} - \beta C_{\rm m}\mathbf{M}_{\rm u}\sum\limits_{i=1}^{2}\frac{\alpha_i}{\Delta t}\mathbf{u}_\mathbf{r}^{l-i} + \beta \mathbf{M}_{\rm iu} I_{\rm ion}({\rm U}\widetilde{\mathbf{u}}_\mathbf{r}^{l},~\mathbf{w}_\mathbf{r}^{l}))  \\
&\qquad\qquad-\sum\limits_{l=1}^L (\mathbf{r}_\mathbf{r}^{l})^T(\dfrac{\mathbf{w}_\mathbf{r}^{l}-\mathbf{w}_\mathbf{r}^{l-1}}{\Delta t} + g({\rm U}\widetilde{\mathbf{u}}_\mathbf{r}^{l}, \mathbf{w}_\mathbf{r}^{l}) ) .
\end{split}
\end{equation}
The adjoint form of the reduced discretized Monodomain system can be constructed by setting $\dfrac{\partial \mathcal{L}_{\rm r}}{\partial \mathbf{u}_\mathbf{r}^{l}} = 0$ for $l = 1, \cdots, L$. It reads
\begin{equation}\label{}
\mathbb{Z}_{\rm u}^T \mathbf{A_{\rm m}}\mathbb{Z}_{\rm u} \mathbf{q}_\mathbf{r}^{l}  =
 \mathbf{d}_\mathbf{r}^{l} \quad\quad l\in\{1, \cdots, L \} 
\end{equation}
with
\begin{align*}
& \mathbf{d}_\mathbf{r}^{l} = 
\beta{\rm U}^T \{\mathbf{M}_{\rm iu}^T\mathbf{q}_\mathbf{r}^{l+2} \}\circ \partial_u I_{\rm ion}({\rm U}\tilde{\mathbf{u}}_\mathbf{r}^{l+2}, \mathbf{w}_\mathbf{r}^{l+2}) 
-2 \beta{\rm U}^T \{ \mathbf{M}_{\rm iu}^T\mathbf{q}_\mathbf{r}^{l+1} \}\circ \partial_u I_{\rm ion}({\rm U}\tilde{\mathbf{u}}_\mathbf{r}^{l+1}, \mathbf{w}_\mathbf{r}^{l+1}) \\ 
&\qquad  -\partial_ug{\rm U}^T\tilde{\mathbf{r}}_\mathbf{r}^l
+(\mathbf{X}_{\rm u}\mathbf{u}_\mathbf{r}^l - \hat{\mathbf{u}}_{\rm meas}^l ) \chi_{\rm snap}^l 
+ \beta C_{\rm m}\mathbf{M}_{\rm u}\sum\limits_{i=1}^{2}\frac{\alpha_i}{\Delta t}\mathbf{q}_\mathbf{r}^{l+i}. 
\end{align*} 
Here the operation $\circ$ means entry-wise product.
The dual gating variable $\mathbf{r}_\mathbf{r}^l$ is updated by the equation below derived from setting $\dfrac{\partial \mathcal{L}_{\rm r}}{\partial \mathbf{w}_\mathbf{r}^{l}} = 0$ for $l = 1, \cdots, L$:
\begin{equation}
\frac{\mathbf{r}_\mathbf{r}^{l+1}-\mathbf{r}_\mathbf{r}^{l}}{\Delta t} =
\partial_wg\mathbf{r}_\mathbf{r}^l + \beta \{\mathbf{M}_{\rm iu}^T\mathbf{q}_\mathbf{r}^{l} \}\circ \partial_{w}I_{\rm ion}({\rm U}\tilde{\mathbf{u}}_\mathbf{r}^l, \mathbf{w}_\mathbf{r}^{l}).
\end{equation}
For superscripts exceeding $L$, we set $\mathbf{q}_\mathbf{r}^{L+1} = \bsb{0} = \mathbf{q}_\mathbf{r}^{L+2}$ and $\mathbf{r}_\mathbf{r}^{L+1} = \bsb{0} = \mathbf{r}_\mathbf{r}^{L+2}$.

Based on the adjoint equations we then get the derivatives of $\mathcal{J}_{\rm r}$  
\begin{eqnarray}
& \dfrac{{\mathcal D} \mathcal{J}_{\rm r}}{{\mathcal D} \sigma_{\rm mk}}  =
-\sum\limits_{l=1}^L(\mathbf{q}_\mathbf{r}^{l})^T\mathbf{S}_{\rm ku} \mathbf{u}_\mathbf{r}^{l} + 
\dfrac{\alpha}{2}\dfrac{\partial\mathcal{R}}{\partial\sigma_{\rm mk}} \label{DJ-rb1}
\end{eqnarray}
where 
$\mathbf{S}_{\rm ku} = \mathbb{Z}_{\rm u}^T \mathbf{S}_{\rm k}\mathbb{Z}_{\rm u} \in \mathbb{R}^{N\times N}$ and ${\rm k}$ stands for ${\rm l}$ and ${\rm t}$. Notice that 
$\mathbf{A_r} = \mathbb{Z}_{\rm u}^T \mathbf{A_{\rm m}}\mathbb{Z}_{\rm u}  =  \beta C_{\rm m}\frac{\alpha_0}{\Delta t}\mathbf{M}_{\rm u} + \sigma_{\rm ml}\mathbf{S}_{\rm lu} + \sigma_{\rm mt}\mathbf{S}_{\rm tu}.$

Eventually, the reduced inverse conductivity problem can be solved by a line search or trust region interior-point method \cite{nocedal2006}, for instance the primal-dual method or the logarithmic barrier method with a quasi-Newton update on the Hessian computation. For simplicity, we describe in Algorithm \ref{optRB-alg} only the line search barrier method. In particular, the norm $||\cdot||$ in step \ref{searchNorm} of Algorithm \ref{optRB-alg} could be customized. For instance, we can use the Euclidean norm of the polar coordinates to introduce a priori information weighting different values as we propose in Sec.~\ref{pod-DEIM-resultInv}. The parameters $\mu$ and $\varsigma$ are related to the enforcement of the unilateral constraints on the solution.
The optimization iteration updating in our simulation is based on the Optizelle package\footnote{http://www.optimojoe.com/products/optizelle/}, an open source library for nonlinear optimization.

If one wants to avoid off-line reduced bases construction, as an alternative the conductivity parameter can be estimated by an online POD-DEIM framework with adaptivity. This is explained by Algorithm \ref{adaptivePOD-alg}.
Specifically, one can solve the full-order Monodomain system at the initial guess, build the POD bases from
snapshots and solve the reduced inverse solver, then use the reduced-order solution as a new initial guess for the next computing cycle. In particular, for the inverse conductivity problem, instead of storing a unique basis by merging previous POD modes we typically store multiple reduced bases online. In the reduced-order computation, we still use polar coordinates to choose a single appropriate basis (see step \ref{searchNorm} of Algorithm \ref{optRB-alg}).

\begin{algorithm}[tp]
\caption{Optimization-POD-DEIM:~ \texttt{opt-reduced($\mathbf{\boldsymbol{\sigma}}^0$, $\{\mathbb{Z}_{\rm u}^i\}_{i=1}^s$, $\{\mathbb{Z}_{\rm ion}^i\}_{i=1}^s$, $\{\boldsymbol{\sigma}_{\rm gen}^i \}_{i=1}^s$, $\varsigma$)}}\label{optRB-alg}
\begin{algorithmic}[1]
\Input initial guess $\mathbf{\boldsymbol{\sigma}}^0\in\mathcal{C}_{\rm ad}$, POD basis $\{\mathbb{Z}_{\rm u}^i\}_{i=1}^s$ and $\{\mathbb{Z}_{\rm ion}^i\}_{i=1}^s$, POD bases generating parameters $\{\boldsymbol{\sigma}_{\rm gen}^i \}_{i=1}^s$, factor $\varsigma\in(0, 1)$
\Output estimated conductivity values $\mathbf{\boldsymbol{\sigma}}$
\State $\boldsymbol{\sigma}\leftarrow \boldsymbol{\sigma}^0, \mu \leftarrow 1, k\leftarrow 0, i_0\leftarrow -1$ 
\While{$k<k_{\rm max}$ and not converged}
\While{stopping tolerance not reached}
	\State Search $i_* = arg\min\limits_{1\leq j\leq s}||\boldsymbol{\sigma}-\boldsymbol{\sigma}_{\rm gen}^j || $ \label{searchNorm}
	\If{($i_*\neq i_0 )$} 
		\State Import the POD basis $\mathbb{Z}_{\rm u}^{i_*}$ and $\mathbb{Z}_{\rm ion}^{i_*}$
		\State $i_0 \leftarrow i_*$
	\EndIf
	\State Solve $\mathbf{u_r}^{1\cdots L}$ with $\boldsymbol{\sigma}$ from $t^1$ to $t^L$, using bases $\mathbb{Z}_{\rm u}^{i_*}$ and $\mathbb{Z}_{\rm ion}^{i_*}$  
	\State Compute the perturbed cost function value $\mathcal{J}_{\rm r}(\boldsymbol{\sigma}) - \mu \sum\limits_i\log \big(h_i(\boldsymbol{\sigma})\big)$
	\State Solve $\mathbf{q_r}^{1\cdots L}$ with $\boldsymbol{\sigma}$ from $t^L$ to $t^1$, using bases $\mathbb{Z}_{\rm u}^{i_*}$ and $\mathbb{Z}_{\rm ion}^{i_*}$
	\State Compute the gradient $\nabla\mathcal{J}_{\rm r}(\boldsymbol{\sigma}) - \mu\sum\limits_i\frac{\nabla h_i(\boldsymbol{\sigma})}{h_i(\boldsymbol{\sigma})}$, using (\ref{DJ-rb1})
	\State Update the inverse Hessian approximation and compute the search direction	 $\boldsymbol{v}^k$ (BFGS \cite{nocedal2006}) 
	\State Set $\boldsymbol{\sigma} = \boldsymbol{\sigma} + \gamma_k \boldsymbol{v}^k$ where $\gamma_k \in (0, \infty)$ is computed from a line search
	\State $k\leftarrow k+1$
\EndWhile
\State $\mu \leftarrow \varsigma\mu$
\EndWhile
\end{algorithmic}
\end{algorithm}

\begin{algorithm}[tp]
\caption{Adaptive POD-DEIM Optimization}\label{adaptivePOD-alg}
\begin{algorithmic}[1]
\Input initial guess $\boldsymbol{\sigma}^0$, factor $\varsigma\in(0, 1)$
\Output estimated conductivity value $\boldsymbol{\sigma}^0$
\State $R_{\rm u} = [~]$, $R_{\rm ion} = [~]$, $\Sigma_{\rm gen} = [~]$
\While{stopping criterion not fulfilled}
	\State Create snapshots: solve full-order $u(\boldsymbol{\sigma}^0)$ and $I_{\rm ion}(\boldsymbol{\sigma}^0)$\label{creSnap}
	\State Build POD bases $\mathbb{Z}_{\rm u}(\boldsymbol{\sigma}^0)$ and $\mathbb{Z}_{\rm ion}(\boldsymbol{\sigma}^0)$ using above snapshots
	\State Include new POD bases: $R_{\rm u} = [R_{\rm u}, \mathbb{Z}_{\rm u}(\boldsymbol{\sigma}^0)]$, $R_{\rm ion} = [R_{\rm ion}, \mathbb{Z}_{\rm ion}(\boldsymbol{\sigma}^0)]$, $\Sigma_{\rm gen} = [\Sigma_{\rm gen}, \boldsymbol{\sigma}^0]$
	\State Solve the reduced-order problem $\min\limits_{\boldsymbol{\sigma}}\mathcal{J}_{\rm r}(\boldsymbol{\sigma})$:~ \texttt{opt-reduced($\boldsymbol{\sigma}^0$, $R_{\rm u}$, $R_{\rm ion}$, $\Sigma_{\rm gen}$, $\varsigma$)}
	\State Set $\boldsymbol{\sigma}^0 \leftarrow \mbox{arg}\min\limits_{\boldsymbol{\sigma}}\mathcal{J}_{\rm r}(\boldsymbol{\sigma})$
\EndWhile
\end{algorithmic}
\end{algorithm}

\section{Model reduction in action: pitfalls and success}\label{DEIMresult}
In the beginning of this section, we will focus only on the offline-online procedure of POD-DEIM approximation. That is, we would like to obtain a set of reduced bases offline so that it can be use at any moment and for many times to solve an inverse conductivity problem. The observation on adaptive POD-DEIM will be provided at the last.

\subsection{POD-DEIM on the forward problem}\label{resultFWD}
To investigate the performance of POD-DEIM model reduction technique on the Monodomain system, we preliminarily solved the ROM on a slab $\Omega=[0, 5]\times[0, 5]\times[0, 0.5] \subseteq \mathbb{R}^3$ with 24272 mesh nodes. In each simulation, five stimuli of $ I_{\rm app} = 10^5 ~\mu \mbox{A/cm}^{3}$ were applied with four at the corners and one at the center of the domain for a duration of $1$~ms. We set the myocardium fibers to be constantly along the $x$-axis. The snapshots for POD basis construction were taken every 0.05 ms with a duration of 25 ms.
\begin{figure}[!tp]
\begin{minipage}{0.30\textwidth}
\centering{POD mode 1}

\includegraphics[scale=0.2]{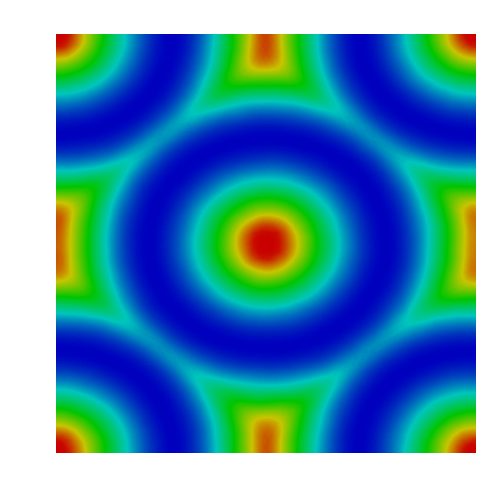}
\end{minipage}
\begin{minipage}{0.30\textwidth}
\centering{POD mode 2}

\includegraphics[scale=0.2]{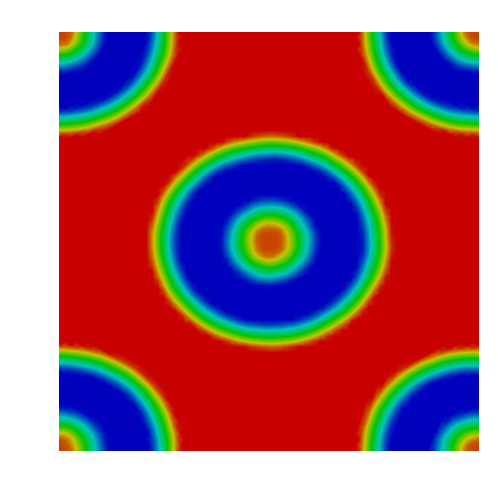}
\end{minipage}
\begin{minipage}{0.35\textwidth}
\centering{POD mode 3\hspace{1cm} \par}
\includegraphics[scale=0.35]{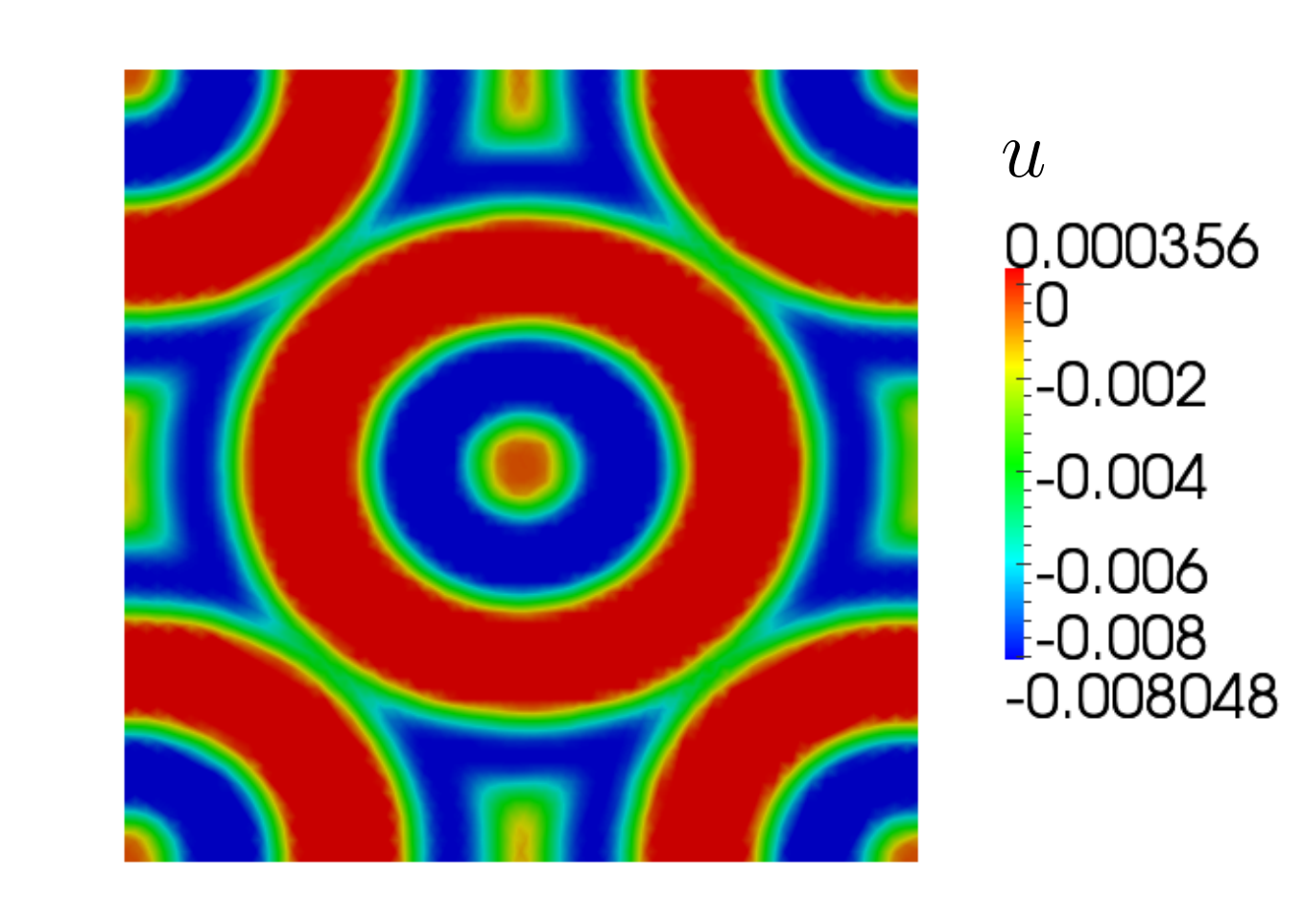}
\end{minipage}

\begin{minipage}{0.30\textwidth}
\centering{POD mode 1}

\includegraphics[scale=0.2]{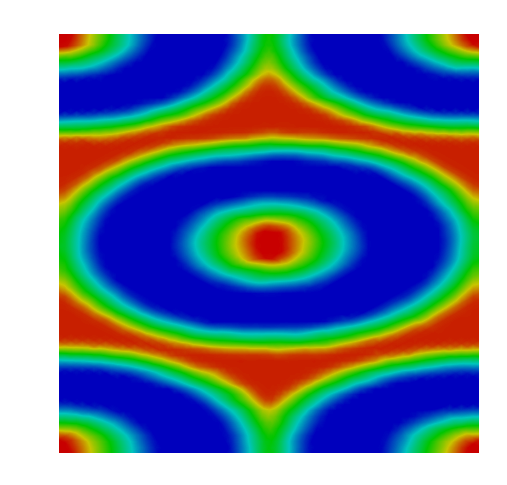}
\end{minipage}
\begin{minipage}{0.30\textwidth}
\centering{POD mode 2}

\includegraphics[scale=0.2]{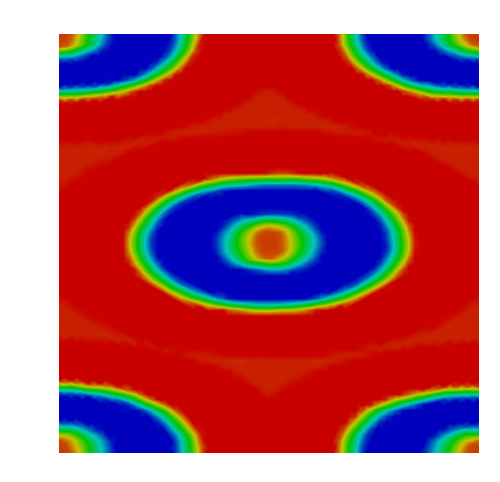}
\end{minipage}
\begin{minipage}{0.35\textwidth}
\centering{POD mode 3\hspace{1cm} \par}
\includegraphics[scale=0.35]{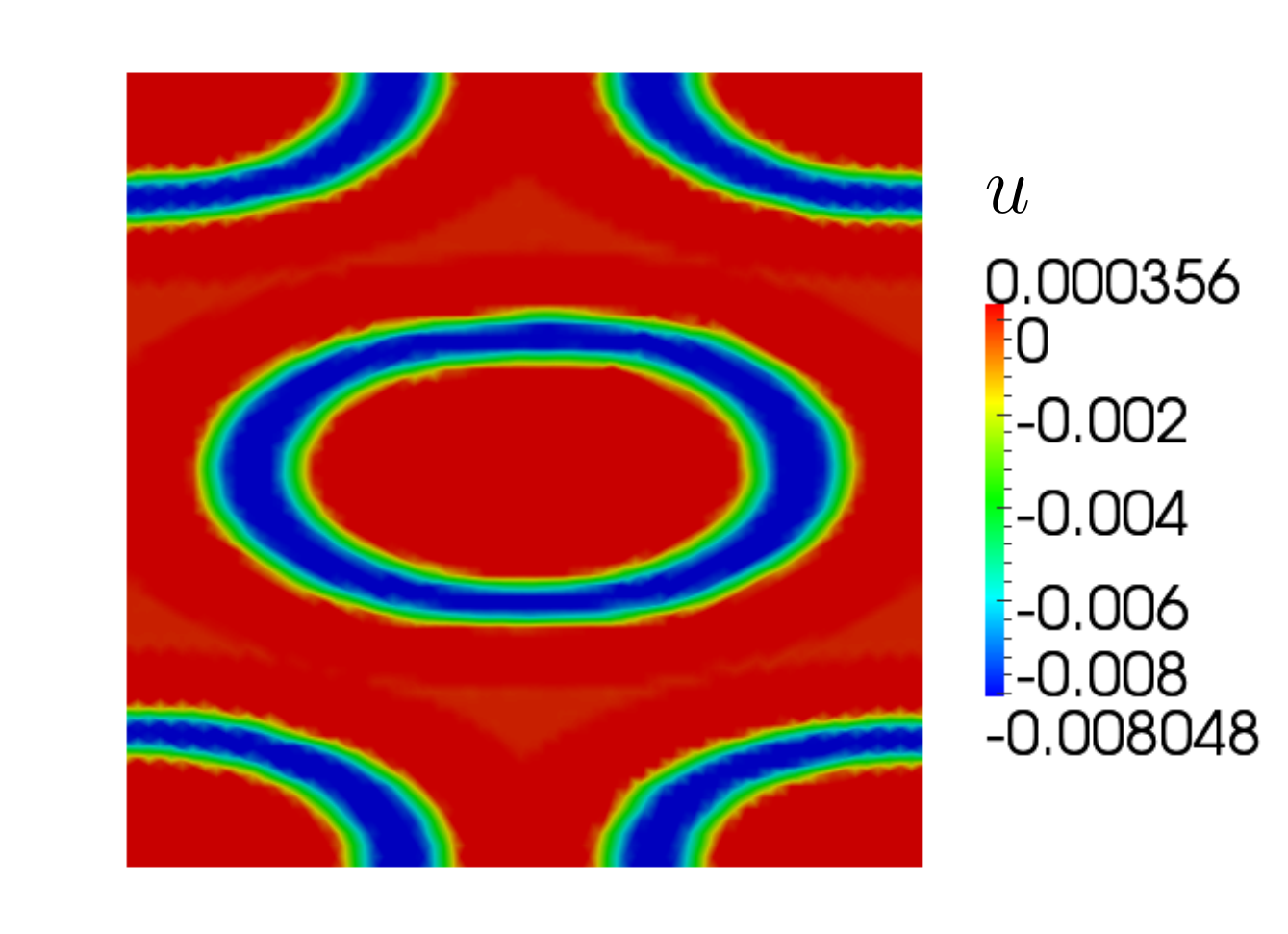}
\end{minipage}
\caption{POD modes of the transmembrane potential constructed with $\boldsymbol{\sigma}_{\rm gen} = [3, 2]$ (first row) and $\boldsymbol{\sigma}_{\rm gen} = [4, 0.1]$ (second row).}
\label{pod-u1-u2}
\end{figure}

In Fig.~\ref{pod-u1-u2} (first row), we plot the first three POD modes of the transmembrane potential constructed from snapshots of $u$ that were computed with conductivity parameter $\boldsymbol{\sigma}_{\rm gen} = [3, 2]$. The parameter chosen has $\sigma_{\rm ml}/\sigma_{\rm mt}$ close to one, which makes the tissue almost isotropic. This explains why the wave front propagated from the center of the slab tissue is almost circular. We also plot in the second row the leading POD modes of the transmembrane potential constructed with $\boldsymbol{\sigma}_{\rm gen} = [4, 0.1]$. The wave front in this case is elliptic since the ratio $\sigma_{\rm ml}/\sigma_{\rm mt} \gg 1$.
This comparison suggests that the POD modes constructed with different conductivity values are very weakly correlated. Therefore, the approach of extracting POD modes from combined snapshots computed with different parameter values can not achieve enough dimension reduction, as already pointed out in Fig.~\ref{svd-onePara-multi} (right).

\begin{figure}[!tp]
\begin{minipage}{0.49\textwidth}
\includegraphics[scale=0.25]{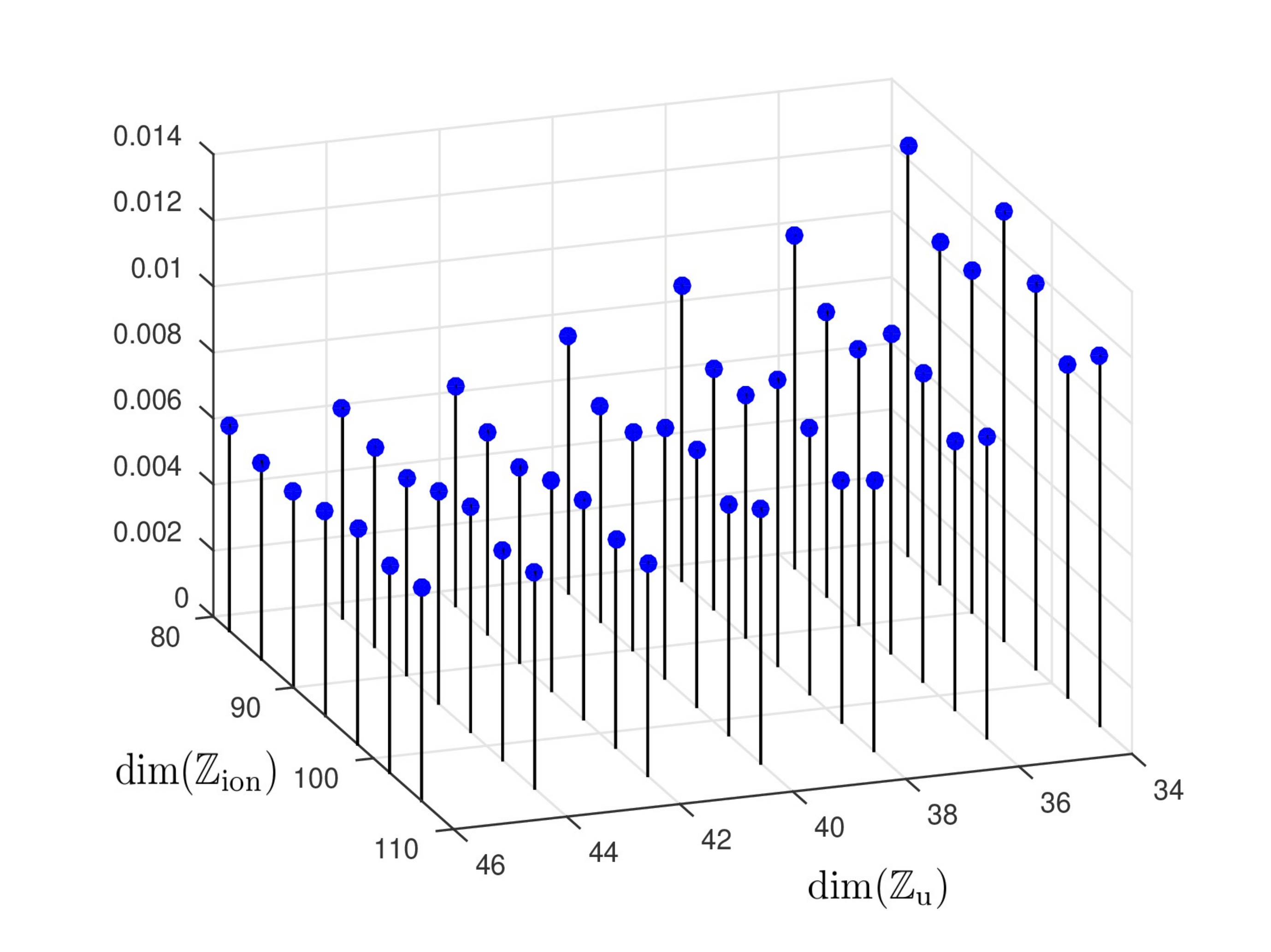}
\end{minipage}
\begin{minipage}{0.49\textwidth}
\includegraphics[scale=0.258]{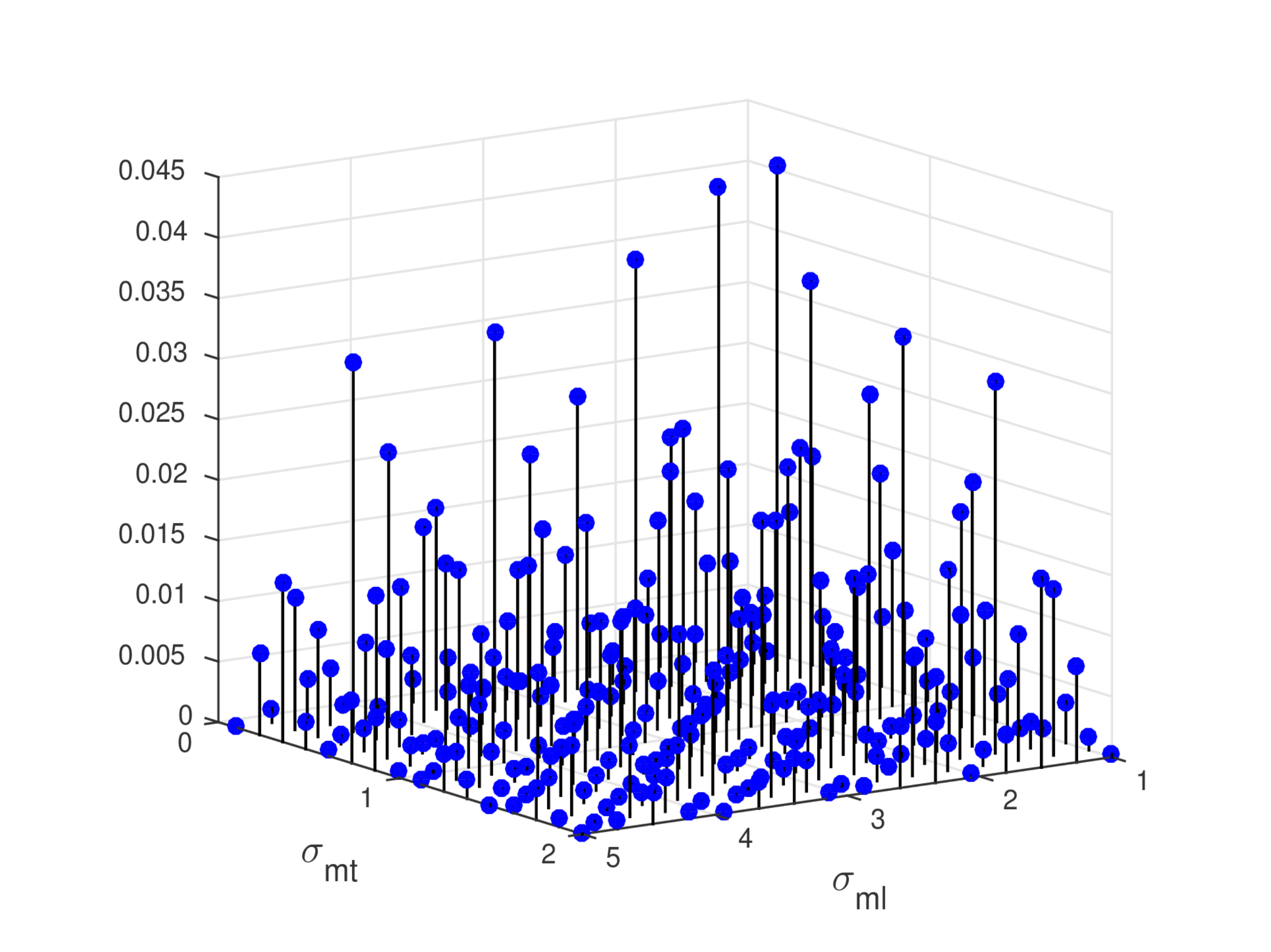}
\end{minipage}
\caption{Errors of $u$ by POD-DEIM approximation. Left: errors w.r.t.~the dimension of the POD basis, with $\dim(\mathbb{Z}_{\rm u})$ ranging from 34 to 46 and $\dim(\mathbb{Z}_{\rm ion})$ ranging from 82 to 106 ; Right: errors on $16\times 16$ different conductivity parameters, with fixed POD-basis dimensions: $\dim(\mathbb{Z}_{\rm u}) = 45$ and $\dim(\mathbb{Z}_{\rm ion}) = 100$.}
\label{err-pod-DEIM}
\end{figure}

To check the performance of reduced-order modeling, we took a uniform sampling on the conductivity parameter: 25 samples were generated over the domain $[1, 5]\times [0, 2]$. For each value of the parameter, snapshots were computed to construct offline a POD basis of the transmembrane potential (ionic current). In this way, 25 POD bases were available for importing during online computation of the reduced-order Monodomain model (Algorithm \ref{optRB-alg}).

As a preliminary assessment of the impact of model reduction, we computed the transmembrane potential with the POD-DEIM procedure. Fig.~\ref{err-pod-DEIM} reports the corresponding relative error. The errors in Fig.~\ref{err-pod-DEIM} (left) are with respect to the dimension of the POD basis, with $\dim(\mathbb{Z}_{\rm u})$ ranging from 34 to 46 and $\dim(\mathbb{Z}_{\rm ion})$ ranging from 82 to 106. The test conductivity parameter is $\boldsymbol{\sigma} = [2.6, 1.05]$. We observe that the POD-DEIM method provides stable and accurate approximation when the POD basis dimension is around 40 for $u$ and 90 for $I_{\rm ion}$.

In Fig.~\ref{err-pod-DEIM} (right) we plot the errors on $16\times 16$ different conductivity parameters, with fixed POD-basis dimensions: $\dim(\mathbb{Z}_{\rm u}) = 45$ and $\dim(\mathbb{Z}_{\rm ion}) = 100$. We observe that 60.55\% of the parameters in this test feature errors below 0.005, 21.48\% of them lead errors greater than 0.01 and they are mainly associated with conductivities at the boundary of the sample space. These errors indicate that a {\it uniform sampling} with 25 parameter values would be generally inappropriate for its use in the inverse problem of conductivity estimation.

\subsection{Domain of effectiveness (DOE) of the reduced basis}\label{DOEresult}
From previous results we argue that a single reduced basis may be effective only for problems with parameter values close to the generator one. To quantify the role of this basis, we introduce the concept of {\it ``domain of effectiveness''} (DOE). This is the region of the parameter space where a reduced-order solution provide accurate results as we specify later on.
We carry out a supportive study on the domains of effectiveness of different reduced bases, when each basis has a unique generating parameter $\boldsymbol{\sigma}_{\rm gen}$ that differs from others in the conductivity ratio.

\begin{figure}[!tp]
\begin{center}
\begin{minipage}{0.49\textwidth}
\includegraphics[scale=0.45]{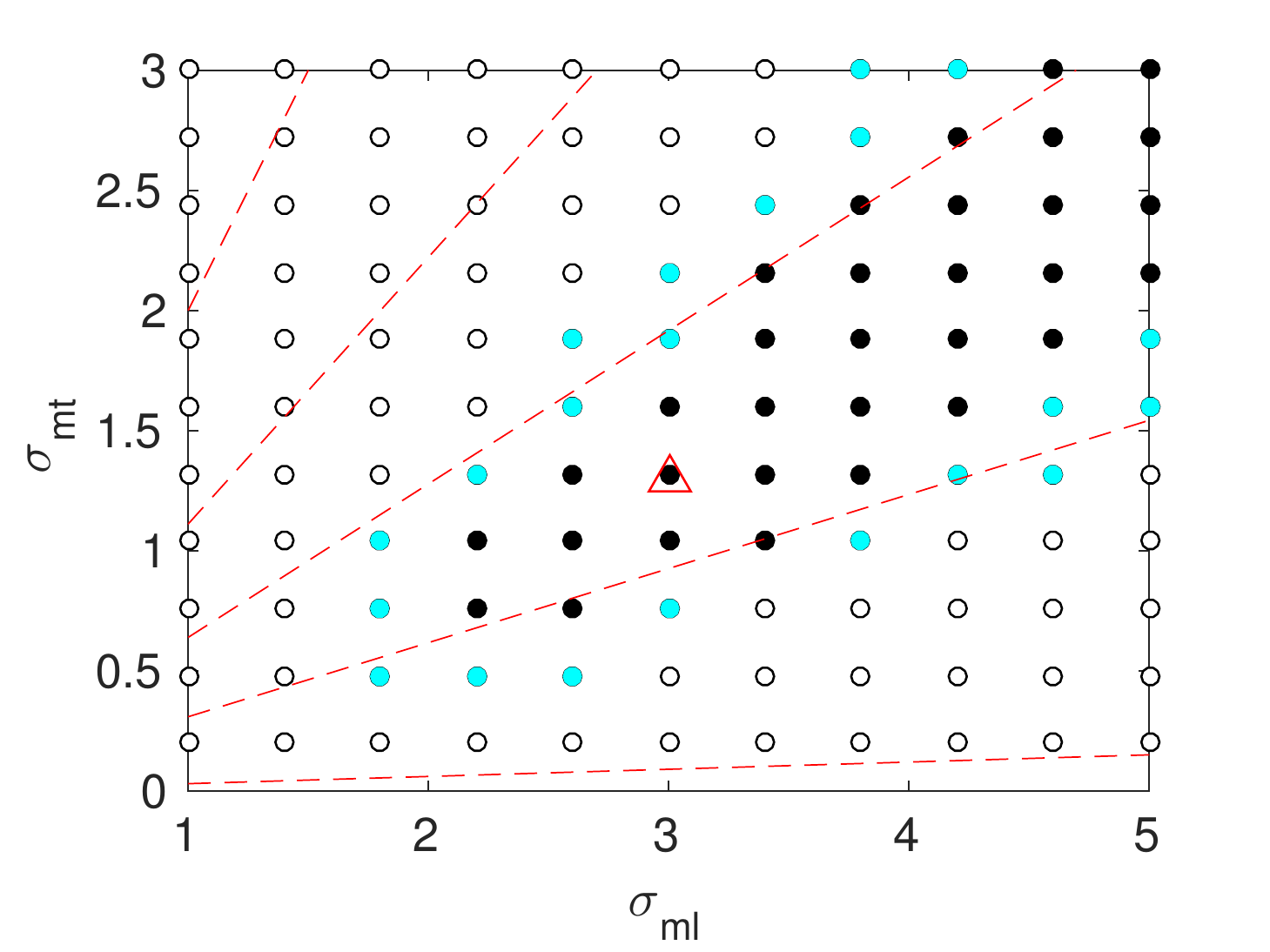} 
\end{minipage}
\begin{minipage}{0.49\textwidth}
\includegraphics[scale=0.45]{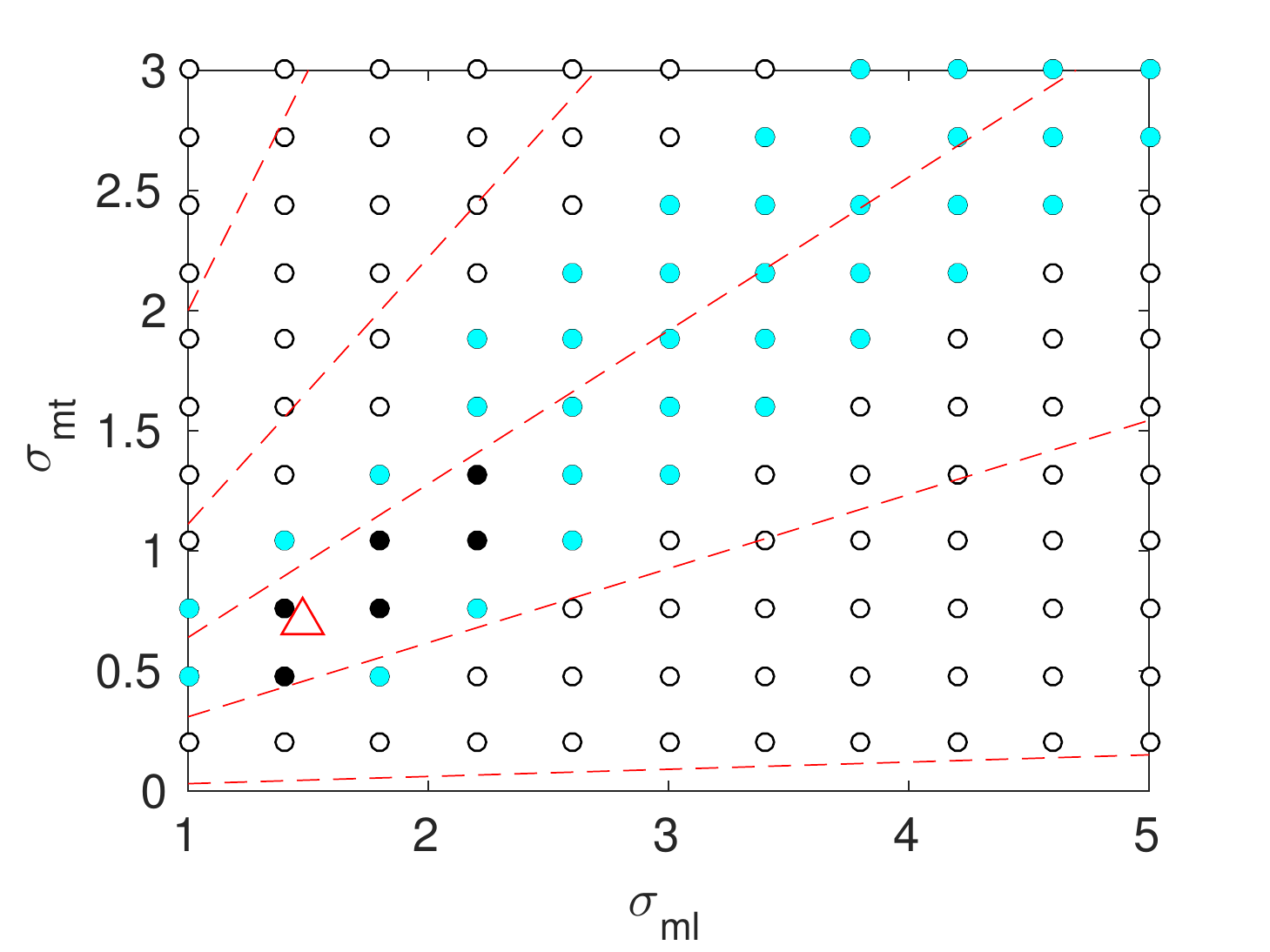}
\end{minipage}

\begin{minipage}{0.49\textwidth}
\includegraphics[scale=0.45]{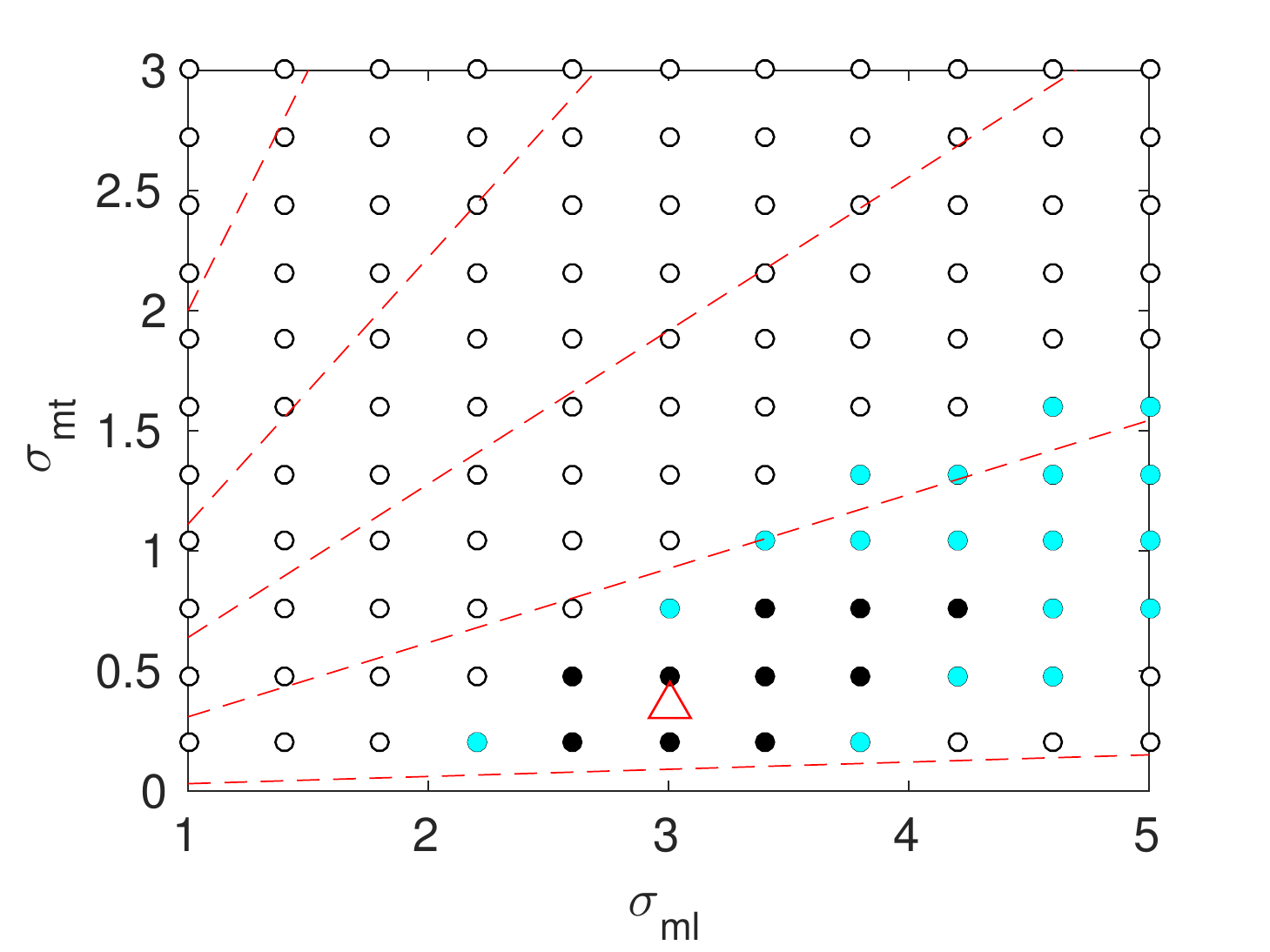}
\end{minipage}
\begin{minipage}{0.49\textwidth}
\includegraphics[scale=0.45]{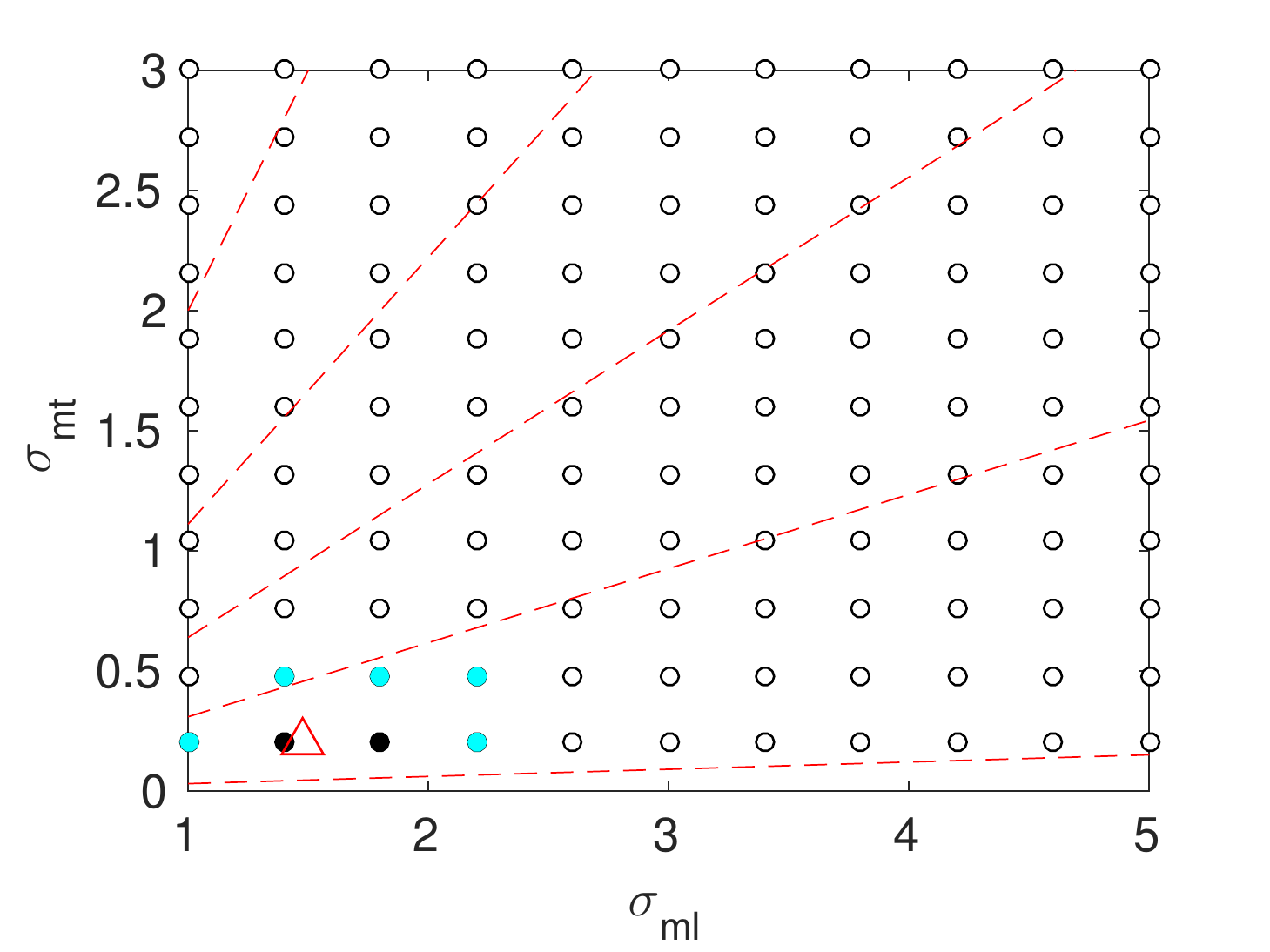}
\end{minipage}
\caption{Domains of effectiveness of different reduced bases. The relative errors (denoted by $e$) of the ROM solutions are indicated in different colors. Black: $e \leq 0.002$, Cyan(Gray): $0.002 < e \leq 0.005$, White: $e > 0.005$. The parameter space is partitioned by the red dash lines using an equi-spaced partition on $\arctan(\frac{\sigma_{\rm mt}}{\sigma_{\rm ml}})$. The POD basis generating parameter $\boldsymbol{\sigma}_{\rm gen}$ in each picture is indicated by a red triangle. }
\label{err-oneRB-noSensi}
\end{center}
\end{figure}

\begin{figure}[!tp]
\begin{minipage}{0.49\textwidth}
\includegraphics[scale=0.45]{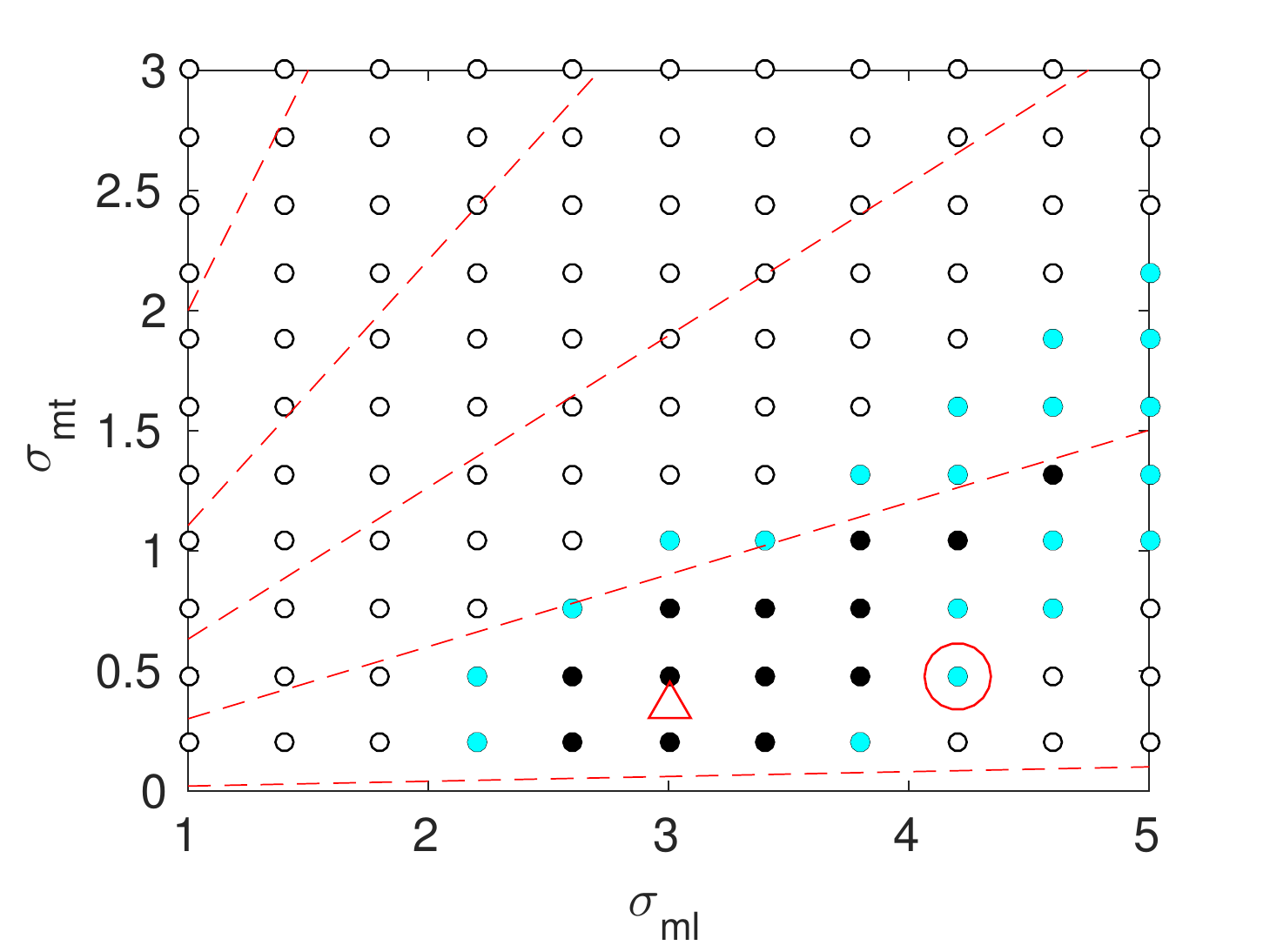} 
\end{minipage}
\begin{minipage}{0.49\textwidth}
\includegraphics[scale=0.45]{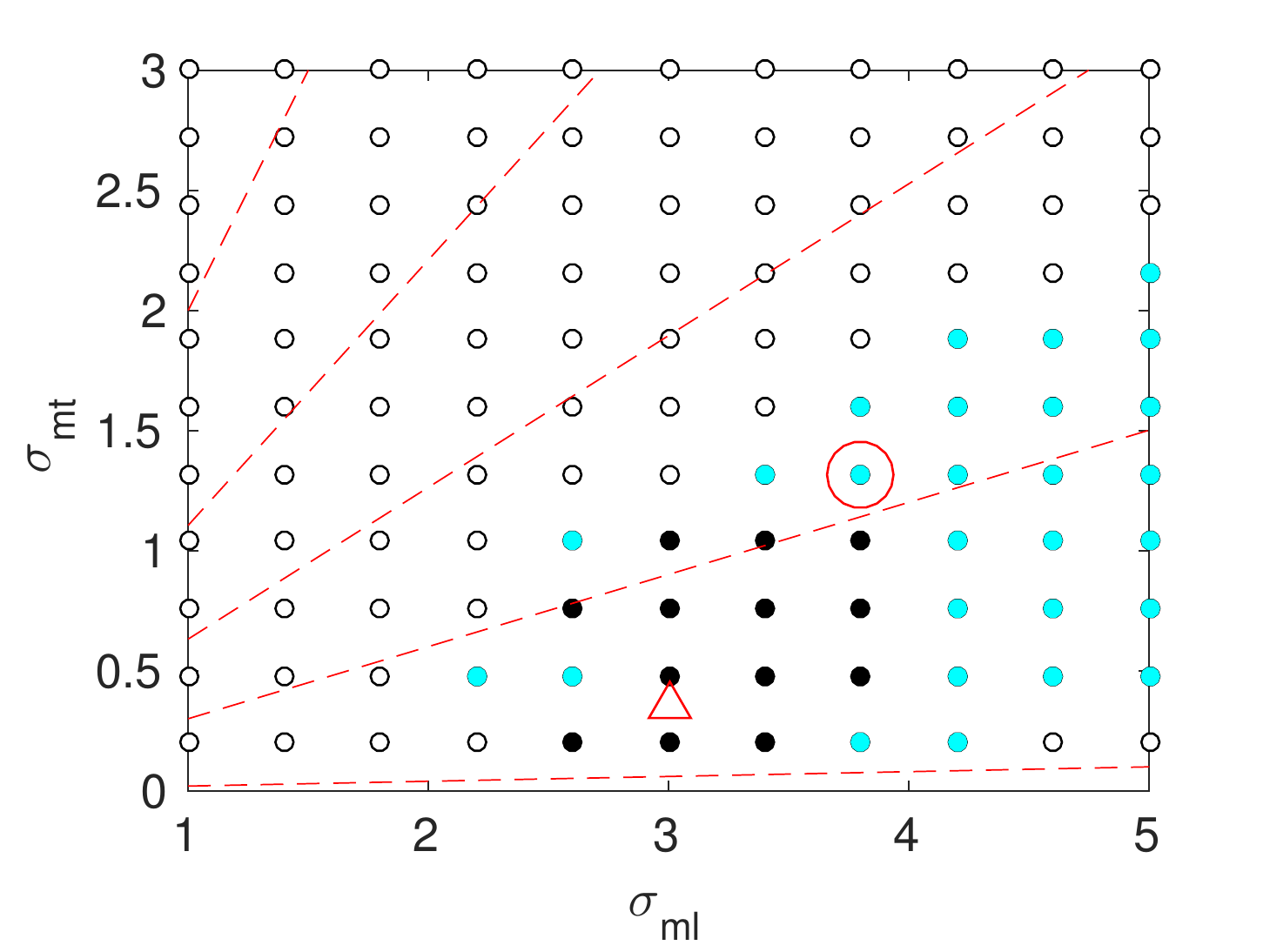}
\end{minipage}
\caption{Domains of effectiveness of sensitivity-based reduced bases both generated by $\bsb{\sigma}_{\rm gen} = [3, 0.35]$. The left one is based on $\bsb{\sigma}_{\rm gen} + [\delta_{\rm l}, \delta_{\rm t}] = [3.8, 1.32]$ (marked by the red circle) and the right one takes $\bsb{\sigma}_{\rm gen} + [\delta_{\rm l}, \delta_{\rm t}] = [4.2, 0.48]$ (marked by the red circle).}
\label{err-oneRB-Sensi}
\end{figure}

Given a POD basis $\mathbb{Z}_{\rm u}$ and a test parameter $\boldsymbol{\sigma}$, we measure the effectiveness of $\mathbb{Z}_{\rm u}$ at $\boldsymbol{\sigma}$ by the relative error of reduced-order solution
$e(\boldsymbol{\sigma}) = {\sum\limits_{l=1}^L || \mathbb{Z}_{\rm u}\mathbf{u}^l_\mathbf{r}-\mathbf{u}^l ||^2}/{\sum\limits_{l=1}^L ||\mathbf{u}^l ||^2}$
where $\mathbf{u}^l$ is the full-order solution at time $t^l$ given $\boldsymbol{\sigma}$ and $\mathbf{u}^l_\mathbf{r}$ is the corresponding reduced-order solution solved with the reduced basis $\mathbb{Z}_{\rm u}$. We study the domains of effectiveness of different reduced bases and plot them in Fig.~\ref{err-oneRB-noSensi}. In each picture, a unique generating parameter $\boldsymbol{\sigma}_{\rm gen}$ (indicated by the red triangle) was used for the construction of the reduced basis. The reduced-order simulation errors are less than 0.002 at black points, between 0.002 and 0.005 at cyan/gray points, and greater than 0.005 at white points. 
We define the DOE associated to the reduced basis as the region collecting black and cyan points, where the reconstruction error is less than 0.005.

In Fig.~\ref{err-oneRB-noSensi}, the parameter space is partitioned by the red dash lines using an equi-spaced partition on the range of $\arctan ({\sigma_{\rm mt}}/{\sigma_{\rm ml}})$.
It is interesting to see that the DOE of a reduced basis is apparently confined to an angular region of its generating parameter. The region is wide when the arc and the angle of $\boldsymbol{\sigma}_{\rm gen}$ are large (as in Fig.~\ref{err-oneRB-noSensi} upper left), and is relatively narrow while the arc or the angle is relatively small (as in Fig.~\ref{err-oneRB-noSensi} upper right, lower left and lower right). 
This study confirms that the transmembrane potentials solved with the Monodomain model are strongly sensitive to the conductivity ratio and amplitude \cite{johnston2011sensitivity}. 

As we mentioned earlier, a natural idea to enlarge the domain of effectiveness of a reduced basis is to include extra sensitivity snapshots in the POD basis construction. Following the idea in Compact POD \cite{Carlberg2008},  we took snapshots for the transmembrane potential $u$ as 
\begin{equation*}
 \bigg[ \mathbf{u}^1, \delta_{\rm l}\frac{\partial\mathbf{u}^1}{\partial \sigma_{\rm ml}}, \delta_{\rm t}\frac{\partial\mathbf{u}^1}{\partial \sigma_{\rm mt}}, \cdots , \mathbf{u}^L, \delta_{\rm l}\frac{\partial\mathbf{u}^L}{\partial \sigma_{\rm ml}}, \delta_{\rm t}\frac{\partial\mathbf{u}^L}{\partial \sigma_{\rm mt}} \bigg],
\end{equation*}
where $\delta_{\rm l}$ and $\delta_{\rm t}$ are scaling factors applied on the sensitivity snapshots.
We also took snapshots for the ionic current $I_{\rm ion}$ in a similar way. To investigate the feasibility of this concept, we chose a generating parameter $\bsb{\sigma}_{\rm gen} = [3, 0.35]$ and compare the DOE of the sensitivity-based RB with the older one shown in Fig.~\ref{err-oneRB-noSensi} lower left. Two sets of scaling factors were chosen in a way such that
 $$\bsb{\sigma}_{\rm gen} + [\delta_{\rm l}, \delta_{\rm t}] = [3.8, 1.32]\qquad\bsb{\sigma}_{\rm gen} + [\delta_{\rm l}, \delta_{\rm t}] = [4.2, 0.48].$$
Two reduced bases were generated accordingly and their DOE are plotted in Fig.~\ref{err-oneRB-Sensi}. The ROM in this test took 100 POD modes for $u$ and 250 for $I_{\rm ion}$ which are appropriate numbers based on our experience. From this test, we notice that adding sensitivity snapshots is not significantly effective for widening the DOE. In addition, the largely increased necessary number of POD modes makes this idea not feasible for application in inverse problem such as conductivity estimation.

\subsection{POD-DEIM for the inverse problem}\label{pod-DEIM-resultInv}
\subsubsection*{The measures}
Considering the loss of accuracy in using the reduced-order model, we may need to acquire enough measures to maintain the stability of the inverse conductivity solver. In the following tests, we used $100\times 100$ measurement sites on the tissue surface, which is achievable in experiments\footnote{F. Fenton (GA Tech), personal communication.}. 
We point out that when reducing the number of sites (down to 1000), we experienced just modest (expected)  slowing of the iterative procedure and slightly worse estimation in most test cases.
We took snapshots every $dt_{\rm snap} = 2$ ms for a duration of $T = 30$ ms.

\begin{figure}[!tp]
\begin{center}
\includegraphics[scale=0.2]{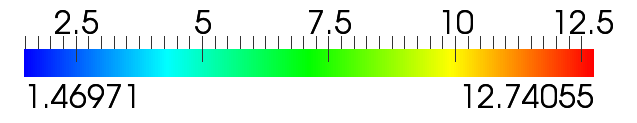}

\begin{minipage}{0.32\textwidth}
\includegraphics[scale=0.2]{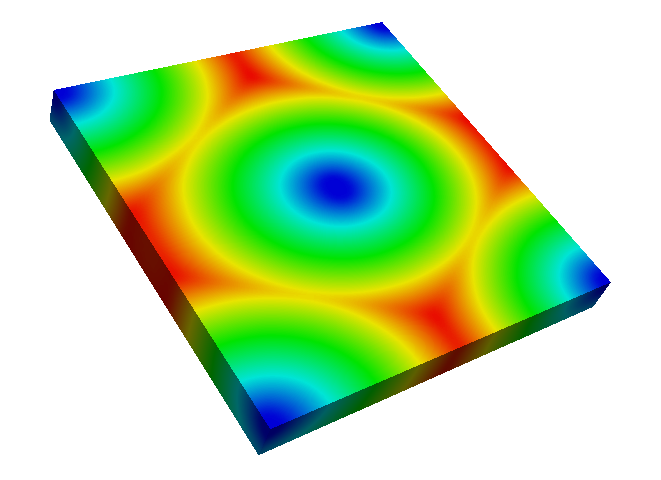}
\end{minipage}
\begin{minipage}{0.32\textwidth}
\includegraphics[scale=0.2]{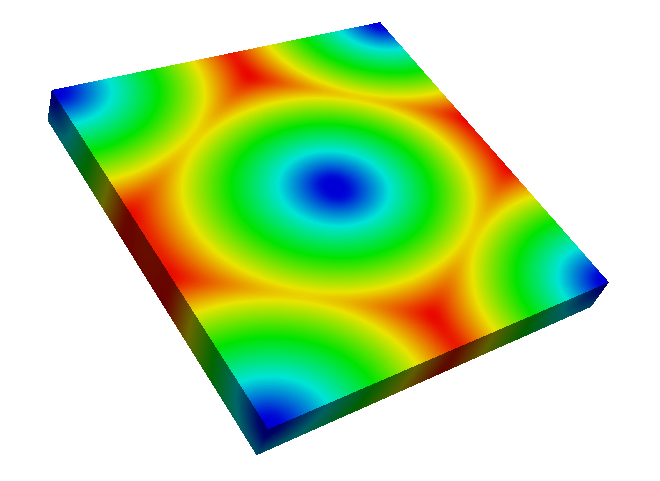}
\end{minipage}
\begin{minipage}{0.32\textwidth}
\includegraphics[scale=0.25]{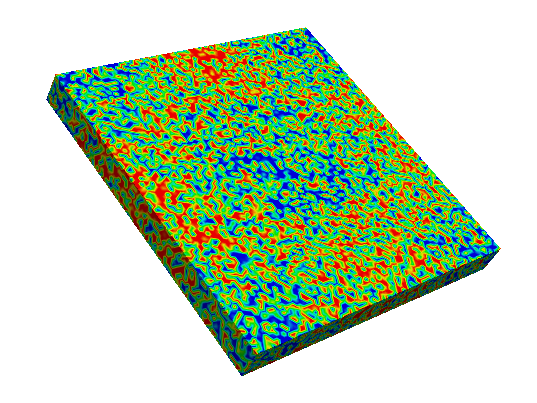}	
\end{minipage}
\caption{Left: $u$ computed with exact finite element approximation on the nonlinear term $I_{\rm ion}$; middle: $u$ computed with component-wise evaluation on the nonlinear term $I_{\rm ion}$; right: synthetic measure of $u$ generated from the left by adding 15\% noise uniformly.}
\label{uPW+noise}
\end{center}
\end{figure}

The DEIM approximation for nonlinearity was applied under the assumption that the nonlinear term $I_{\rm ion}$ is evaluated component-wise in the FOM. This assumption greatly improves the computational efficiency of the full-order system while keeping enough accuracy when the mesh is fine. This is verified by a numerical experiment shown in Fig.~\ref{uPW+noise}. In this test the mesh contains 76832 nodes. As one can see, the component-wise computation on $I_{\rm ion}$ (Fig.~\ref{uPW+noise} middle) presents small error as compared with the exact finite element approximation (Fig.~\ref{uPW+noise} left). More importantly, the computational cost is reduced from 759.284 seconds to 53.6954 seconds (more than 95\% reduction).
Fig.~\ref{uPW+noise} (right) displays the synthetic measures generated by adding 15\% noise to the potential $u$ computed with exact evaluation of $I_{\rm ion}$. This typically represents the way we generate measures in our subsequent numerical tests.

\subsubsection*{Sampling}

\begin{figure}[tp]
\begin{center}
\includegraphics[scale=0.5]{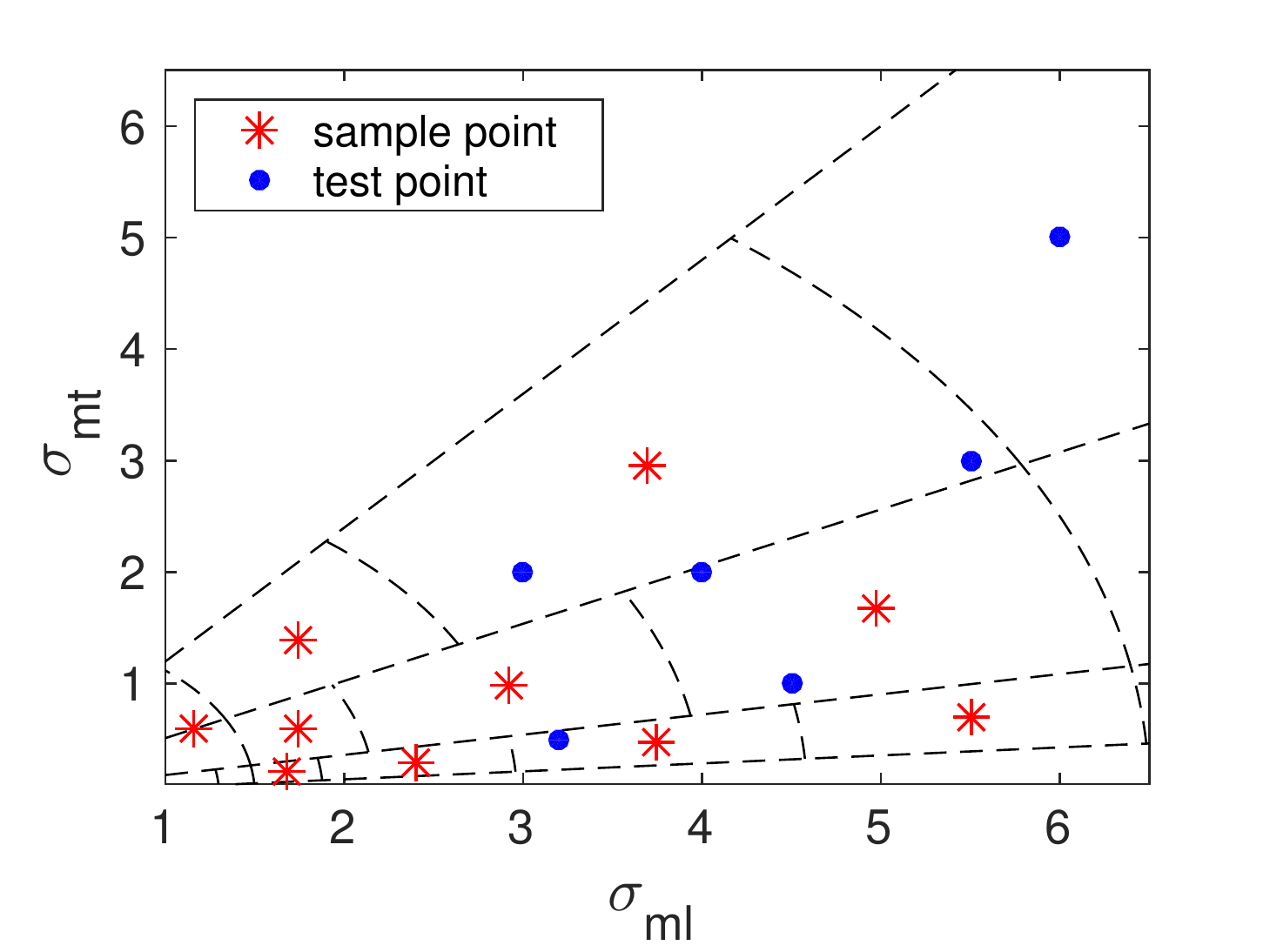}    
\caption{Ten samples (stars in red) generated by a nonuniform sampling on the polar coordinates of conductivity values. The ``small angle--short arc'' zone of the sample space is refined by the use of Gaussian nodes. The six blue dots will be used as test points.}
\label{10smpNewandTst}
\end{center}
\end{figure}

As we have reported in Sec.~\ref{resultFWD}, a uniform sampling on the conductivity parameter to construct POD bases is not a viable approach, since it needs too many sample points and the corresponding ROM approximation still lacks accuracy. However, based on the DOE study in Sec.~\ref{DOEresult}, a nonuniform sampling with refinement in the ``small angle--short arc'' zone could substantially decrease the number of sample points while preserving the accuracy of the ROM approximation. This can be achieved by a sampling on polar coordinates, so to take advantage of the angular pattern of the DOE. 
We illustrate an example in Fig.~\ref{10smpNewandTst}, where nine Gaussian nodes were generated to cover nonuniformly the parameter space and one extra node is added in the left corner to slightly extend the coverage of those samples. Specifically, the Gaussian points were obtained in the polar coordinates $[\rho, \theta]$:
\begin{eqnarray}
&\theta \leftarrow \frac{\theta_i+\theta_{i+1}}{2}, \quad\rho \leftarrow \frac{\rho_j^i+\rho_{j+1}^i}{2}\\
&\theta_i = \theta_{\rm max}-(\theta_{\rm max} - \theta_{\rm min})\cos\big(\frac{(i-1)\pi}{2(n_\theta-1)} \big) , \quad i=1,\cdots n_\theta \\
&\rho_j^i = \rho_{\rm max}-(\rho_{\rm max} - \rho_{\rm min})\cos\big(\frac{(j-1)\pi}{2(n_\rho^i-1)} \big), \quad j=1,\cdots n_\rho^i. 
\end{eqnarray} 
For the example above, we typically take $\theta_{\rm max}=\arctan(1.2)$, ~$\theta_{\rm min}=\arctan(1/14)$, ~$n_\theta=4$,~ $n_{\rho}^i=\max(6-i,3)$, ~$\rho_{\rm max}=6.5$, ~$\rho_{\rm min}=1.5$. The extra point in the left corner was given as $[\rho, \theta] \leftarrow [1.3, \frac{\theta_{\rm min}+\theta_{\rm max}}{2}]$.

\subsubsection*{Conductivity estimation}
We conduct a group of tests on the performance of the POD-DEIM approach applying to the inverse conductivity problem, using the ten samples plotted in Fig.~\ref{10smpNewandTst}. For the sake of reliability of our results, six test points (blue dots in Fig.~\ref{10smpNewandTst}) were carefully chosen not too close to the sample points. With each test point, we solved the full-order Monodomain equation and then added 15\% noise to generate the synthetic measures.
We performed the numerical experiments on a 76832-node mesh with a simulation time step $\Delta t = 0.025$ ms. In each run of the reduced inverse conductivity solver, 35 POD modes were taken for $u$  and 80  for $I_{\rm ion}$, the norm $||\cdot||$ in step \ref{searchNorm} of Algorithm \ref{optRB-alg} was taken as the Euclidean norm of the polar coordinates. The optimization iteration started from a prescribed initial guess $\boldsymbol{\sigma}_{\rm initial}=[1.5, 1]$ and was constrained by a maximum iteration number 40. 
Each simulation was performed on a processor having an Intel(R) Core(TM) i7-3740QM CPU @ 2.70GHz.

\begin{table}[tp]
 \begin{center}
\caption{Conductivity estimation on a slab mesh with DOF = 76832. $T=~30$ ms, $\Delta t=~0.025$ ms, $dt_{\rm snap}=~2$ ms, noise = 15\%, $\boldsymbol{\sigma}_{\rm initial}=[1.5, 1]$.  }
{\footnotesize
\begin{tabular}{|c||c|c|c|c|c|}
\hline\hline
& $\boldsymbol{\sigma}_{\rm exact}$ & $\boldsymbol{\sigma}_{\rm estimated}$ & \# fwd $|$ bwd & Total exe.~time & Time perc\\
\hline
	Full Order & [3.2, 0.5] & [3.237, 0.625] &  51 $|$ 12  & 2720 s& 100\%\\  
	POD+DEIM  & [3.2, 0.5] & [3.084 0.4277] &  103 $|$ 40  & 123.7  s & 4.5\%\\ 
\hline
	Full Order & [4.5, 1] & [4.535, 1.099] &  58 $|$ 17  & 5409 s&  100\%\\  
	POD+DEIM  & [4.5, 1] & [4.482 0.9846] & 161 $|$ 25 &  72.62 s &  1.3\%\\    
\hline
	Full Order & [5.5, 3] & [5.538, 3.059] &  109 $|$ 36   & 8699 s & 100\%\\   
	POD+DEIM  & [5.5, 3] & [5.293 3.099] & 156 $|$ 24 &  58.25 s &  0.6\%\\     
\hline
	Full Order & [4, 2] & [4.045, 2.076] &  71 $|$ 23   & 5388 s & 100\%\\   
	POD+DEIM  & [4, 2] & [3.762 1.965] & 94 $|$ 40 &  97.44 s & 1.8\%\\     
\hline
	Full Order & [3, 2] & [3.055, 2.077] &  59 $|$ 19   & 5626 s & 100\%\\  
	POD+DEIM  & [3, 2] & [2.741, 2.148] & 83 $|$ 11 &  49.52 s & 0.9\%\\ 
\hline
	Full Order & [6, 5] & [6.038, 5.049] &  50 $|$ 16   & 3366 s & 100\%\\   	
	POD+DEIM  & [6, 5] & [6.035, 4.910] & 54 $|$ 16 &  31.31 s & 0.9\%\\ 
\hline
\end{tabular}
}
\label{76k6testTable}
\end{center}
\end{table}

\begin{figure}[!h]
\begin{minipage}{0.49\textwidth}
\includegraphics[scale=0.38]{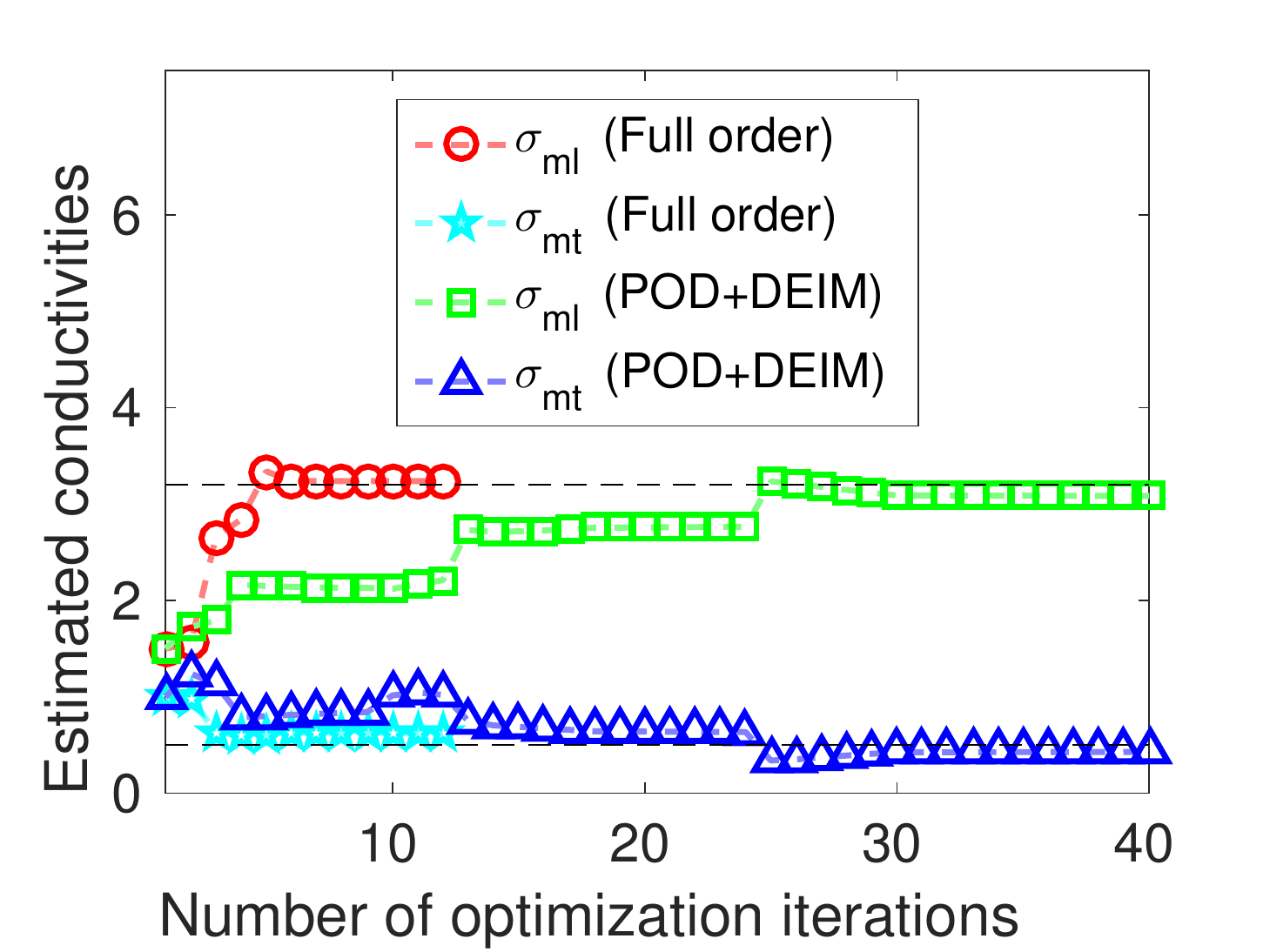} 
\end{minipage}
\begin{minipage}{0.49\textwidth}
\includegraphics[scale=0.38]{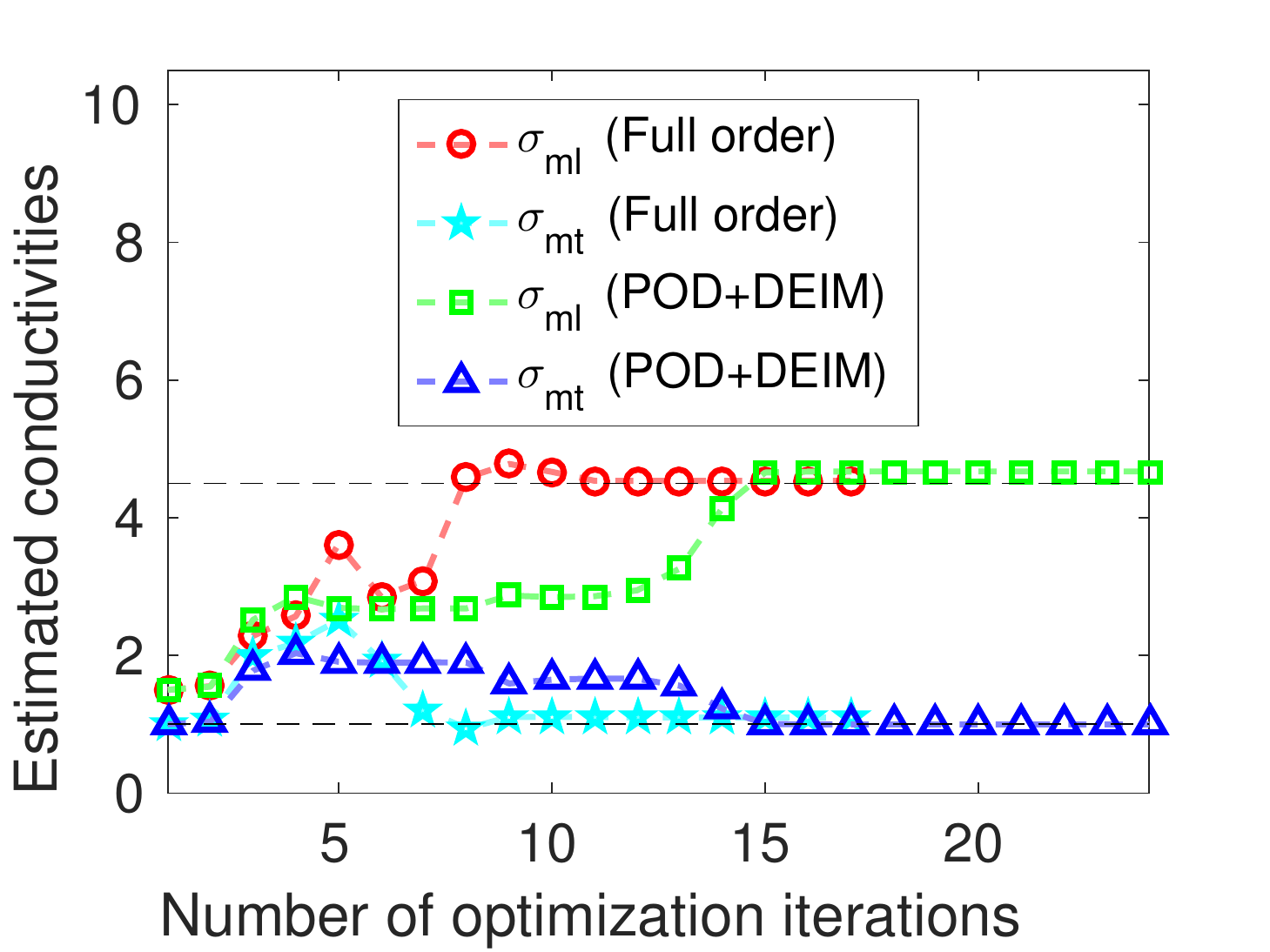}
\end{minipage}

\vspace{1.2cm}

\begin{minipage}{0.49\textwidth}
\includegraphics[scale=0.38]{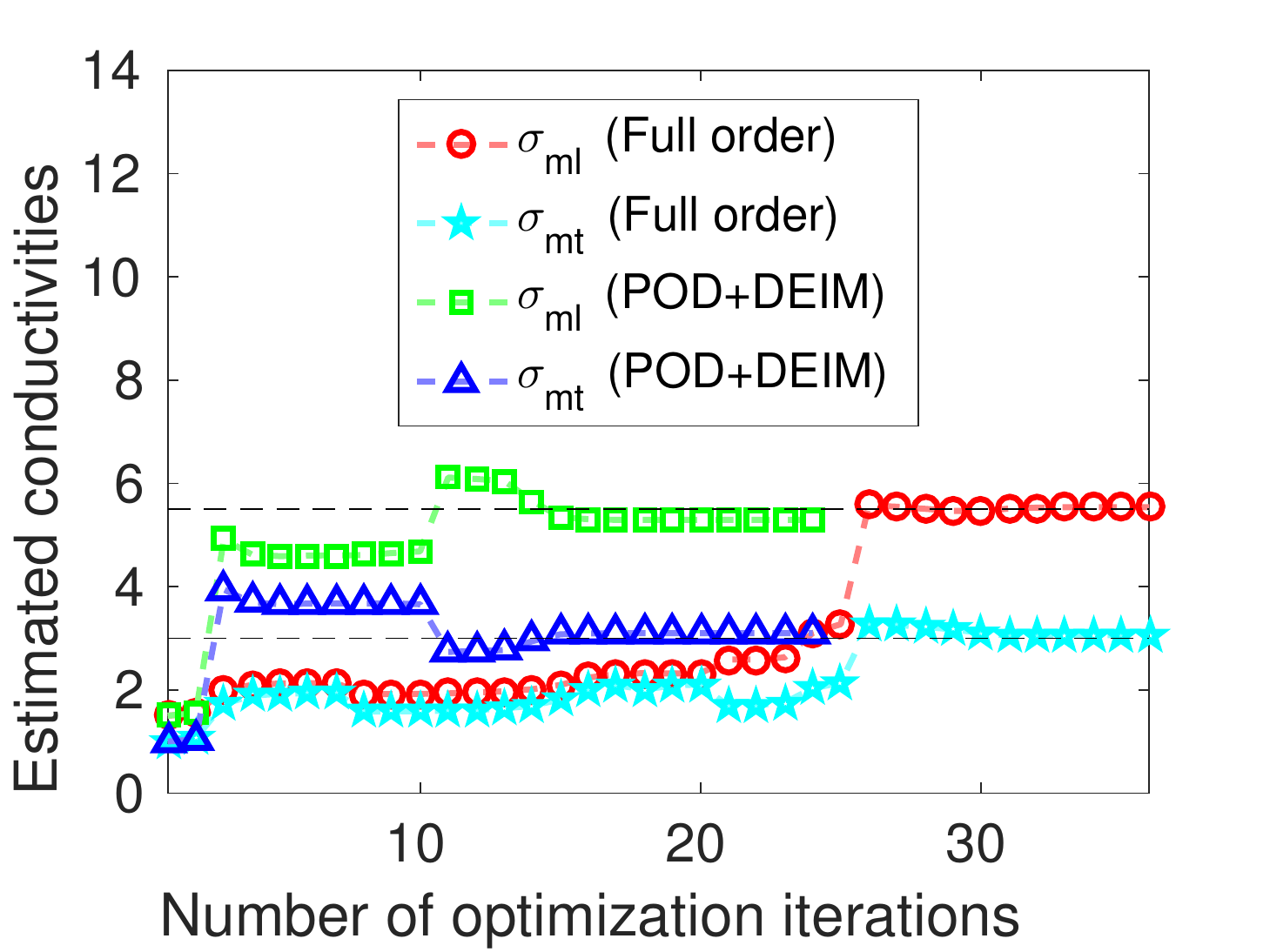}
\end{minipage}
\begin{minipage}{0.49\textwidth}
\includegraphics[scale=0.38]{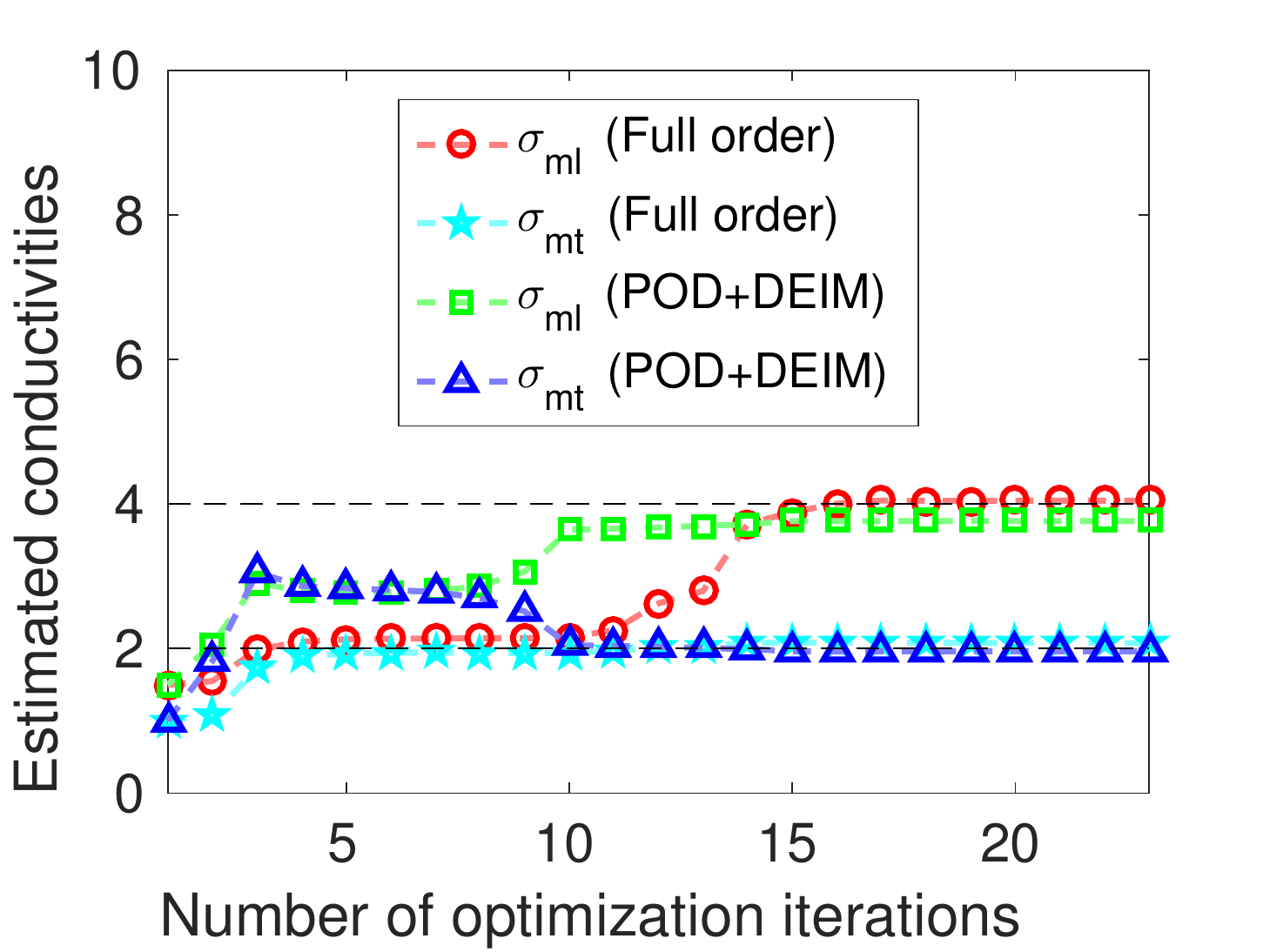}
\end{minipage}

\vspace{1.2cm}

\begin{minipage}{0.49\textwidth}
\includegraphics[scale=0.375]{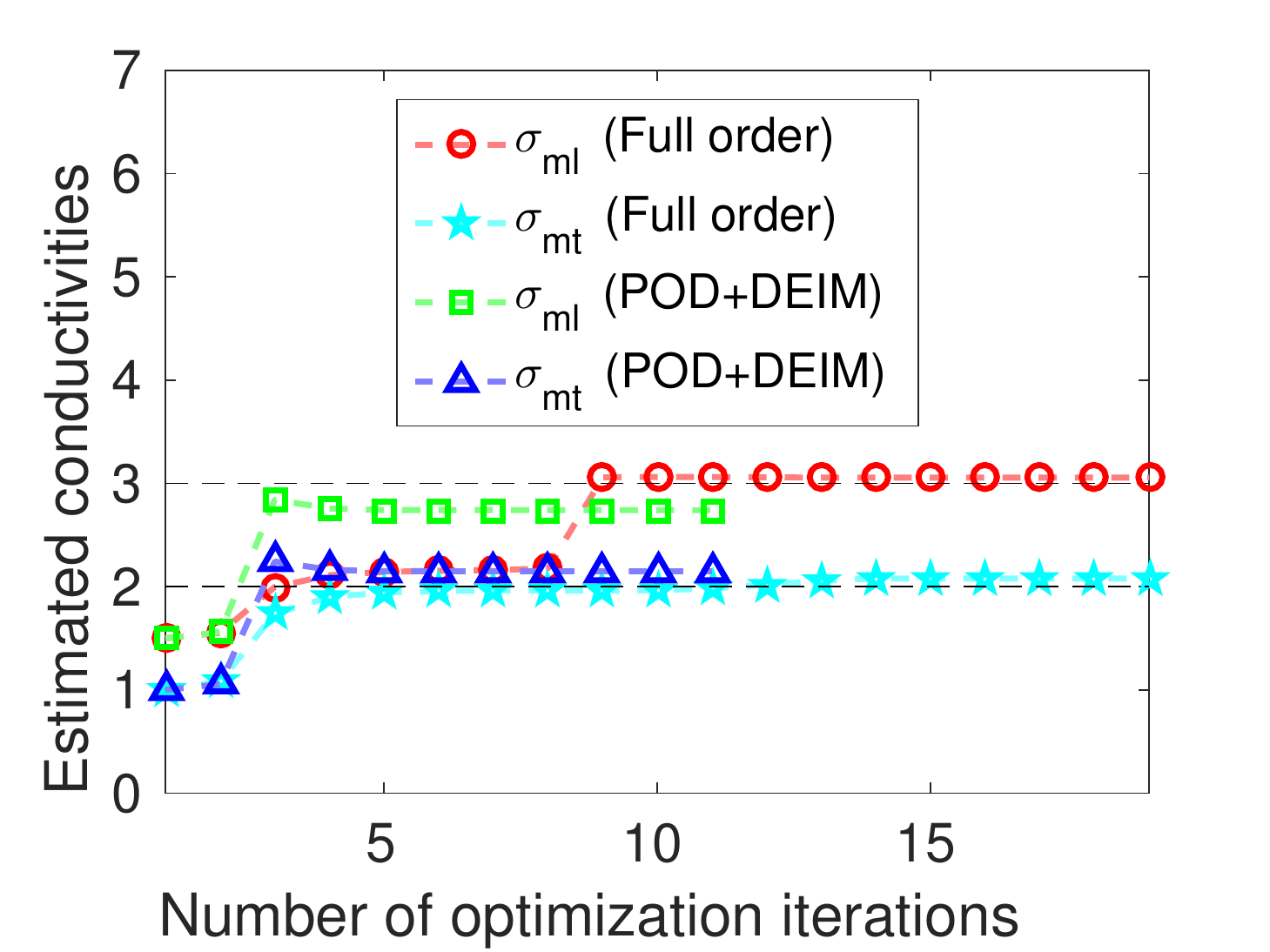}
\end{minipage}
\begin{minipage}{0.49\textwidth}
\includegraphics[scale=0.375]{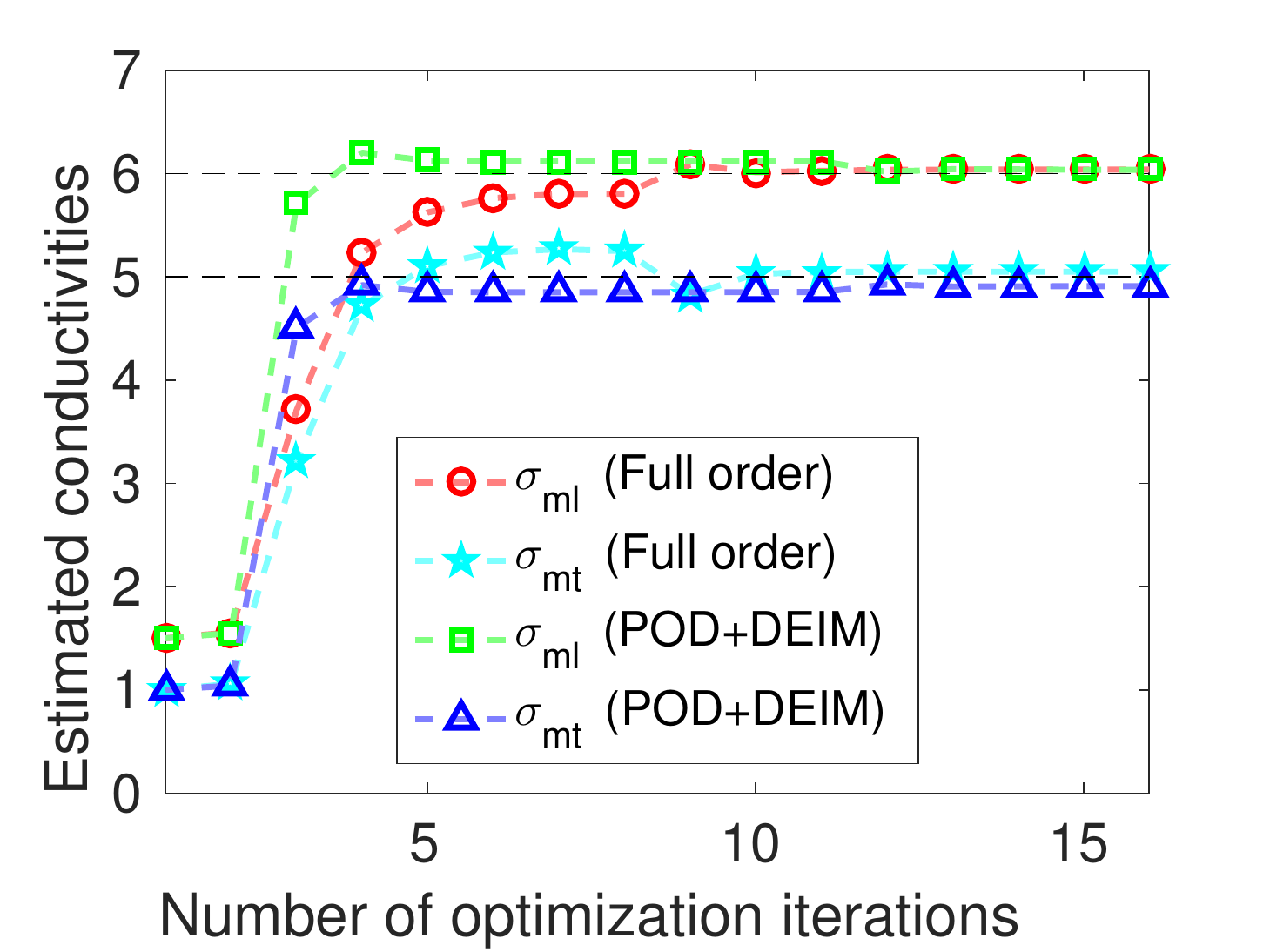}
\end{minipage}
\caption{Optimization iterations corresponding to Table \ref{76k6testTable}}
\label{6testPt-itrFunDJ}
\end{figure}

\begin{figure}[!h]
\begin{center}
\begin{minipage}{0.49\textwidth}
\includegraphics[scale=0.38]{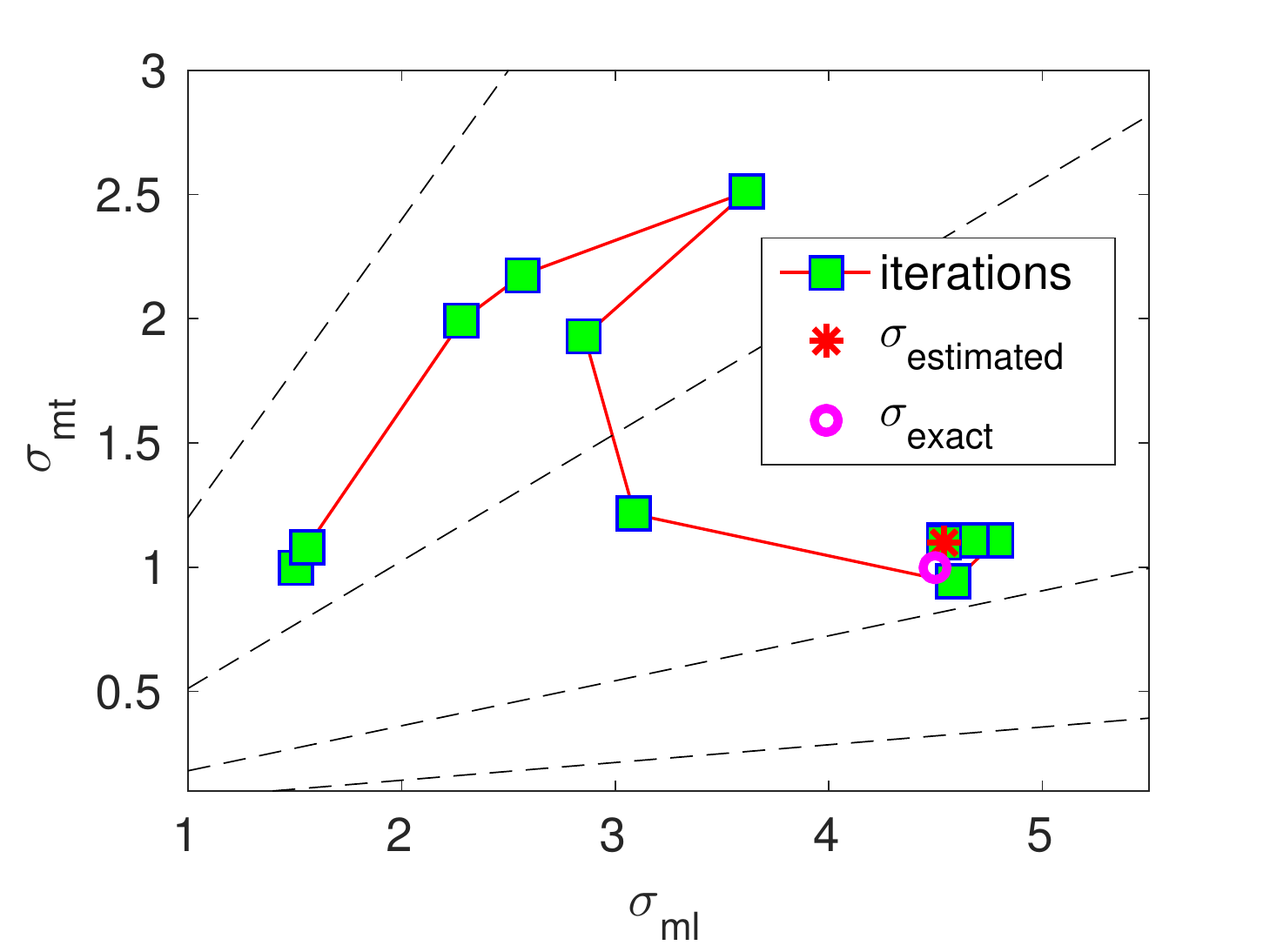}
\end{minipage}
\begin{minipage}{0.49\textwidth}
\includegraphics[scale=0.38]{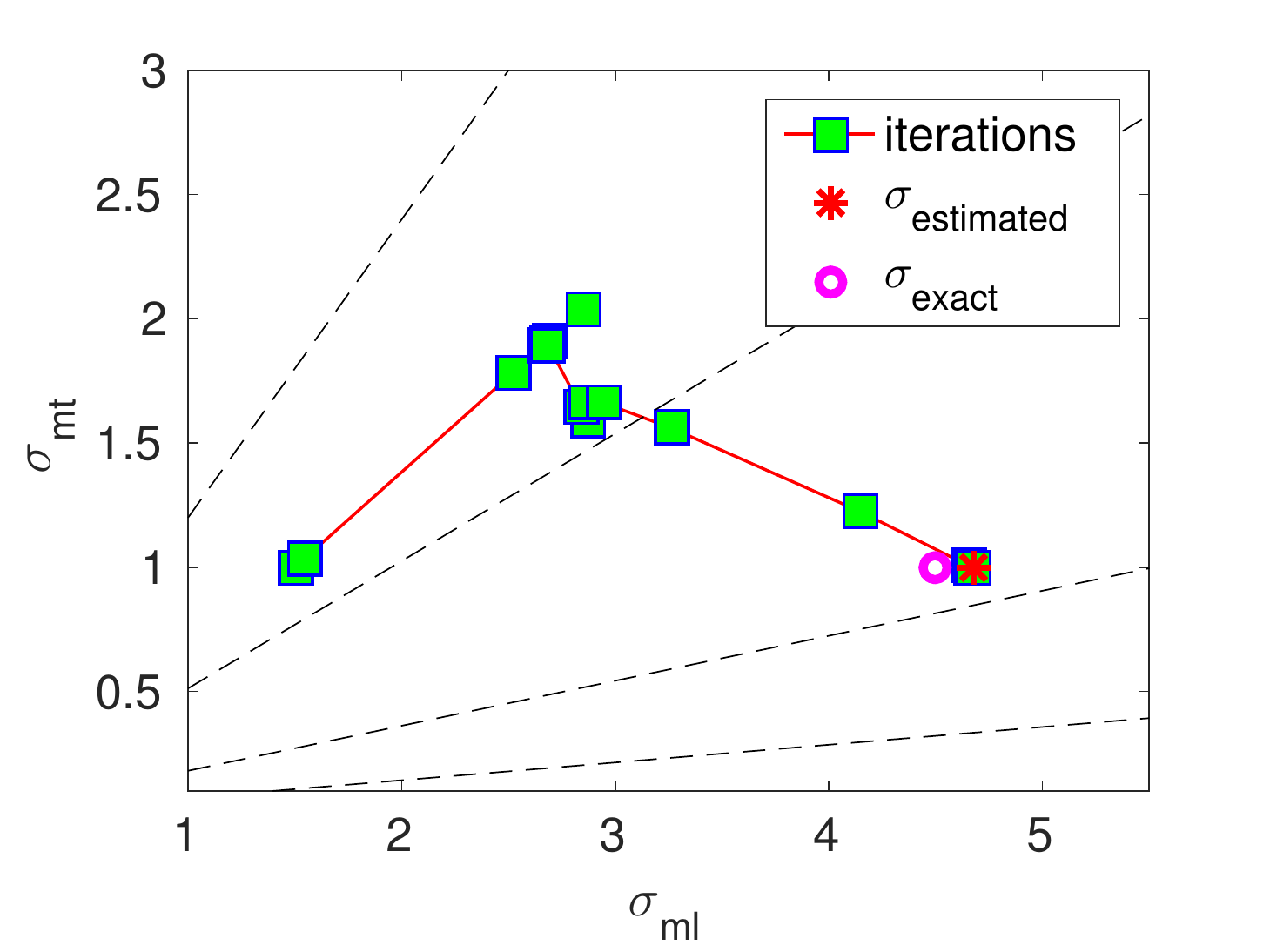}
\end{minipage}

\vspace{1cm}

\begin{minipage}{0.49\textwidth}
\includegraphics[scale=0.38]{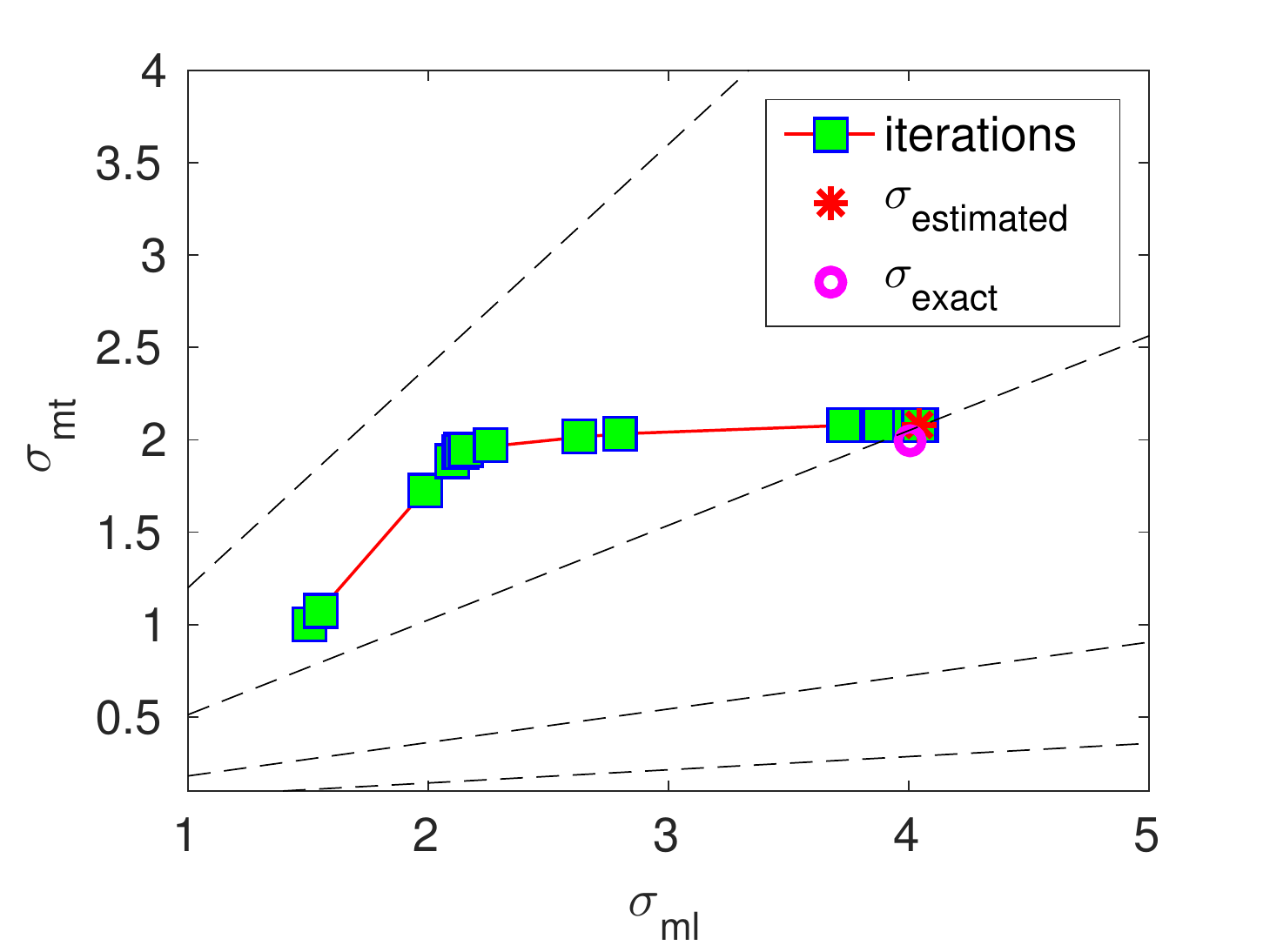}
\end{minipage}
\begin{minipage}{0.49\textwidth}
\includegraphics[scale=0.38]{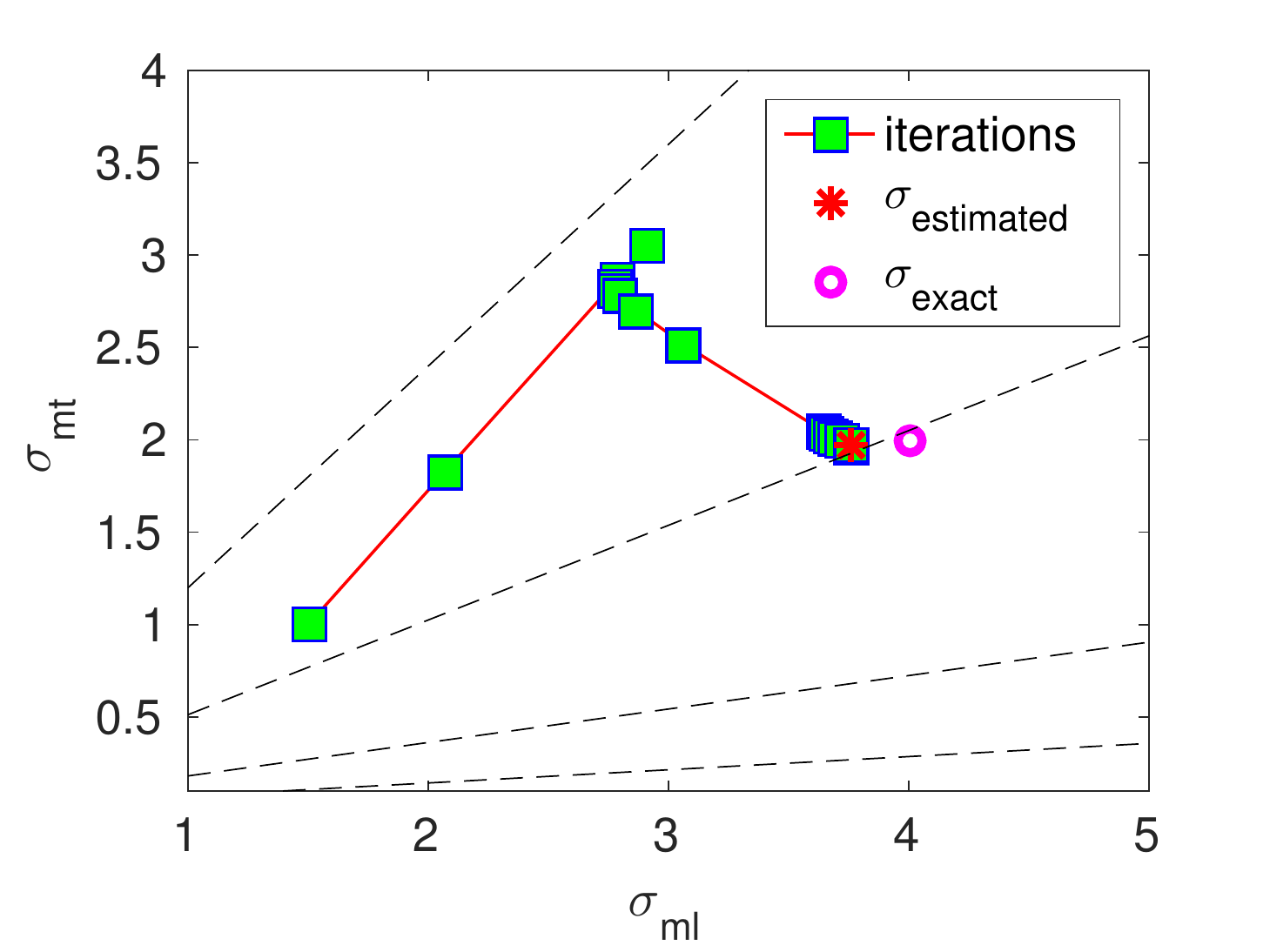}
\end{minipage}

\vspace{1cm}

\begin{minipage}{0.49\textwidth}
\includegraphics[scale=0.38]{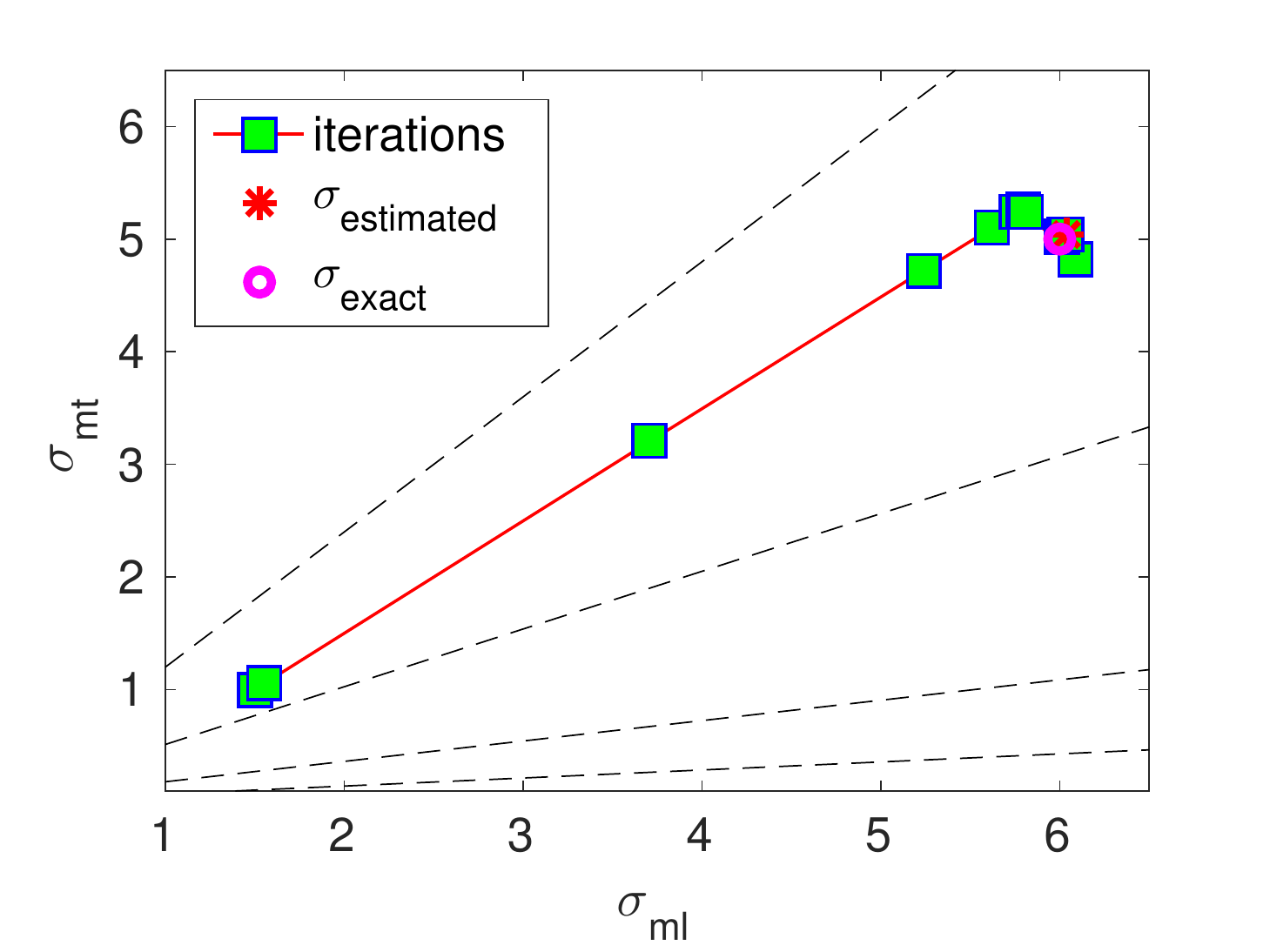}
\end{minipage}
\begin{minipage}{0.49\textwidth}
\includegraphics[scale=0.38]{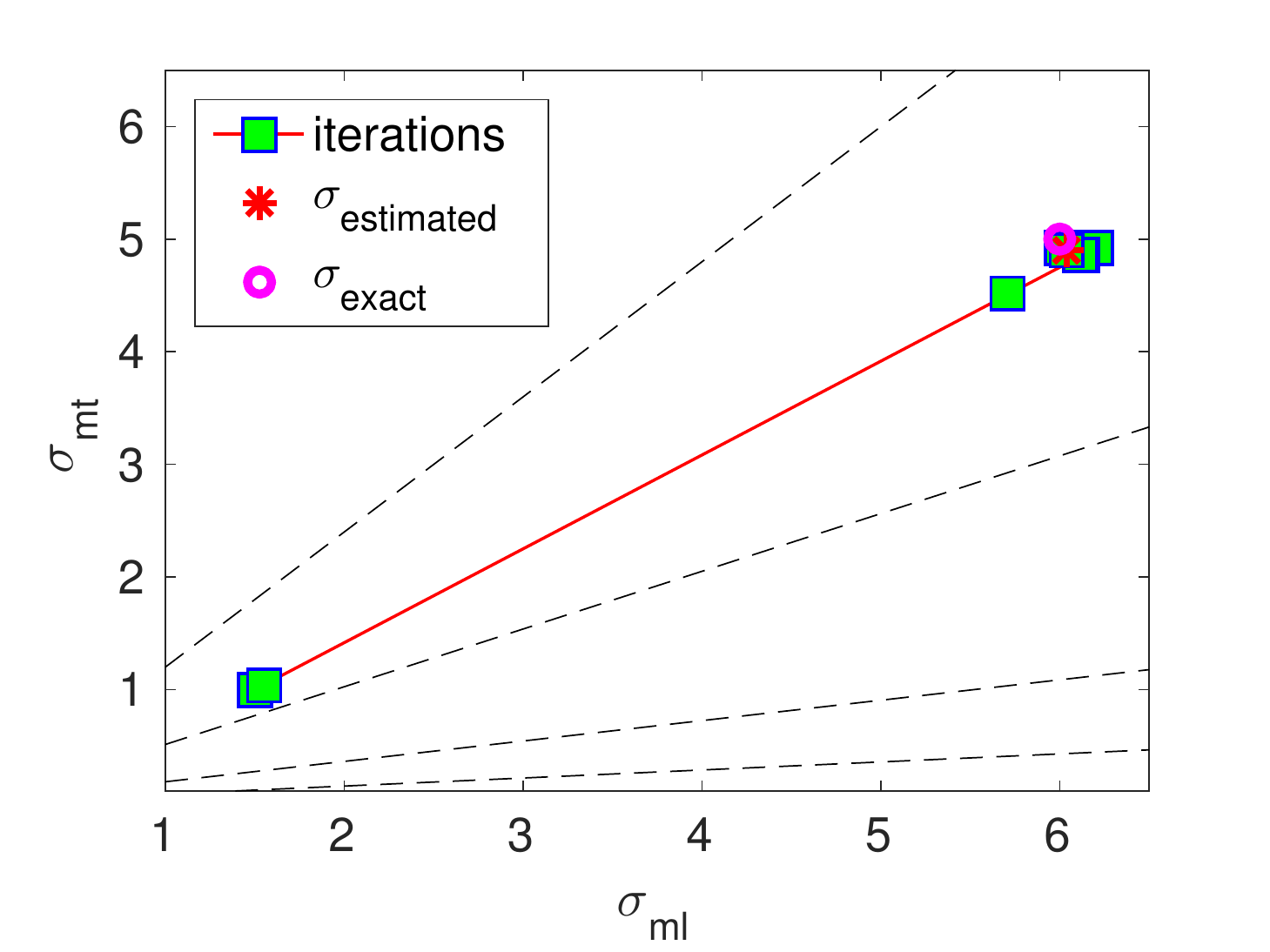}
\end{minipage}
\caption{Optimization trajectory corresponding to Table \ref{76k6testTable} for test points: [4.5, 1], [4, 2], [6, 5]. The left column corresponds to the full-order inverse solver, while the right column corresponds to the reduced inverse solver.}
\label{3testPt-itrPlane}
\end{center}
\end{figure}

\begin{figure}[!ht]
\begin{center}
\begin{minipage}{0.45\textwidth}
\includegraphics[scale=0.32]{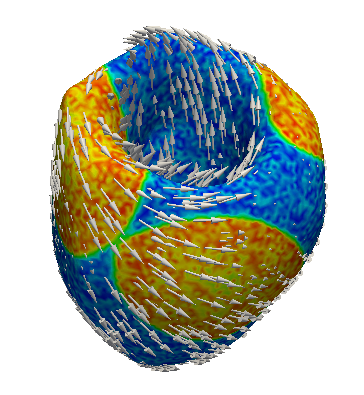}
\end{minipage}
\begin{minipage}{0.45\textwidth}
\includegraphics[scale=0.32]{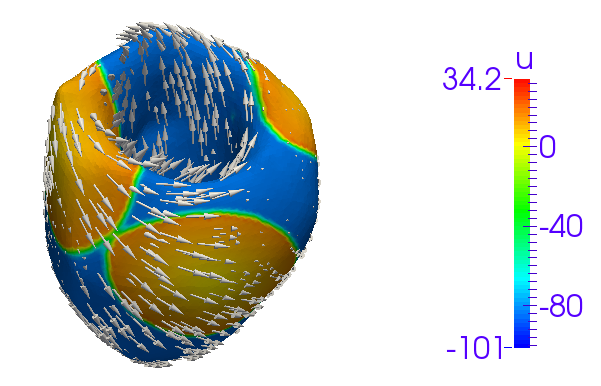}
\end{minipage}
\caption{ Screenshots of the transmembrane potential (in mV) at $t = 26$ ms, computed on a real left ventricular
geometry reconstructed from SPECT images. The white arrows represent myocardial fiber orientation used in the simulation. Left: synthetic measure of the potential created by simulating with $\bsb{\sigma}_{\rm exact} = [3.2, 0.5]$ and adding 15\% uniform noise. Right: reconstruction of the potential computed with the estimated conductivity $\bsb{\sigma}_{\rm estimated} = [3.07, 0.425]$.}
\label{snapOnVentricle}
\end{center}
\end{figure}

The estimated conductivities are listed in Table \ref{76k6testTable}, from which we infer the following conclusion. (a) On average each solve of the reduced Monodomain system (or its dual), including reduced basis importing, takes about only 1 ms as compared with 60 ms for the full-order model. The computational cost reduction on the forward problem is of about two orders of magnitude. (b) For conductivity estimation, the reduced inverse solver returns slightly worse estimates than the full-order inverse solver, due to the loss of accuracy of the ROM. However, the total execution time can be reduced by at least 95\% and the estimation results using ROM are very satisfying with all the test points. Reasonable and similar results were obtained with the coarse 24272-node mesh, which we do not report for the sake of brevity. In summary, the message we retrieve from these results is that even with the 24272-node mesh, notwithstanding the obvious difference between the potential computed with the exact finite element approximation on nonlinearity and the potential solved by component-wise evaluation on nonlinearity, the conductivity estimations using ROM are still vey good.

The iterations corresponding to the 76832-node case are plotted in Fig.~\ref{6testPt-itrFunDJ}. Some of them are also shown in the plane as displayed in Fig.~\ref{3testPt-itrPlane}. We can see that the optimal pathway in the reduced inverse solver is similar to that in the full-order inverse solver.

\subsubsection*{Simulations on the left ventricle} 
To demonstrate the independence of the method on geometries, we also performed simulations on a real left ventricular geometry reconstructed from SPECT images (see Figure \ref{snapOnVentricle}). 
The main features of the cardiac fiber field (shown as white arrows in Figure \ref{snapOnVentricle}) were modeled by an analytical representation of the fiber orientation, which was originally proposed in \cite{colli2004} eqn.~(6.2) for an ellipsoid domain and properly adapted to a real domain retrieved from SPECT images (as done in \cite{gerardo2009model}). Other computational models of cardiac fibers, such as \cite{bayer2012}, are also available.
In this test, we still took measurements every $dt_{\rm snap} = 2$ ms for a duration of $T = 30$ ms.
The synthetic measure of the transmembrane potential at time $t=26$ ms is displayed in the left of Figure \ref{snapOnVentricle}, which was created by simulating with $\bsb{\sigma}_{\rm exact} = [3.2, 0.5]$ and adding 15\% uniform noise. 

With the same sample points of the conductivity parameter mentioned before, we created the corresponding set of reduced bases and employ them in the POD-DEIM procedure. The main difference between current and previous tests lies in the geometric complexity of the fiber orientation. To capture enough variety of the propagation pattern on the ventricle, we need to increase the number of POD modes in the reduced space. A typical example is at the first test point $\boldsymbol{\sigma}_{\rm exact} = [3.2, 0.5]$. When 35 POD modes were taken for $u$ and 80 for $I_{\rm ion}$ we obtain $\boldsymbol{\sigma}_{\rm estimated} = [2.491, 0.6]$. After increasing the number of POD modes to 45 and 90 respectively, the estimation is clearly improved to $\bsb{\sigma}_{\rm estimated} = [3.07, 0.425]$.
For further tests we refer to Table \ref{22kTable}. On the ventricle, the reduced inverse solver lost certain accuracy as compared to the slab tissue, but still returned acceptable estimations of the conductivity. We reconstruct the transmembrane potential corresponding to the first test point and show the snapshot at time $t=26$ ms in the right of Figure \ref{snapOnVentricle}. As demonstrated by the figure, the reconstruction matches the synthetic measure on the wave front. 

\begin{table}[tp]
 \begin{center}
\caption{Conductivity estimation by the reduced inverse solver on a left ventricular mesh: $\boldsymbol{\sigma}_{\rm initial}=[1.5, 1]$.  }
\label{22kTable}
{\footnotesize
\begin{tabular}{c|cccccc}
\toprule
$\boldsymbol{\sigma}_{\rm exact}$ 	& [3.2, 0.5] & [4.5, 1] & [5.5, 3] & [4, 2]  & [3, 2] & [6, 5] \\  %
\midrule
$\boldsymbol{\sigma}_{\rm estimated}$  & [3.07, 0.425] & [4.757, 0.94] & [5.356, 3.151] & [3.933 1.782]  & [2.8, 2.07] & [6.186, 5.488]\\ 
\bottomrule
\end{tabular}
}
\label{22k6testTable}
\end{center}
\end{table}


\subsubsection*{Adaptive POD-DEIM} 
We finally test the viability of adaptation by an online procedure. As described in Sec.~\ref{pod}, POD adaptivity has been studied in many publications, such as \cite{zahr2015progressive}. The main purpose is to let the sampling procedure take into account the optimization trajectory. Since we have observed that merging snapshots from different conductivity values doesn't reduced the necessary number of POD modes for accurate potential approximation, we still used multiple reduced bases in online reduced-order computation.

\begin{figure}[tp]
\begin{center}
\begin{minipage}{0.49\textwidth}
\includegraphics[scale=0.45]{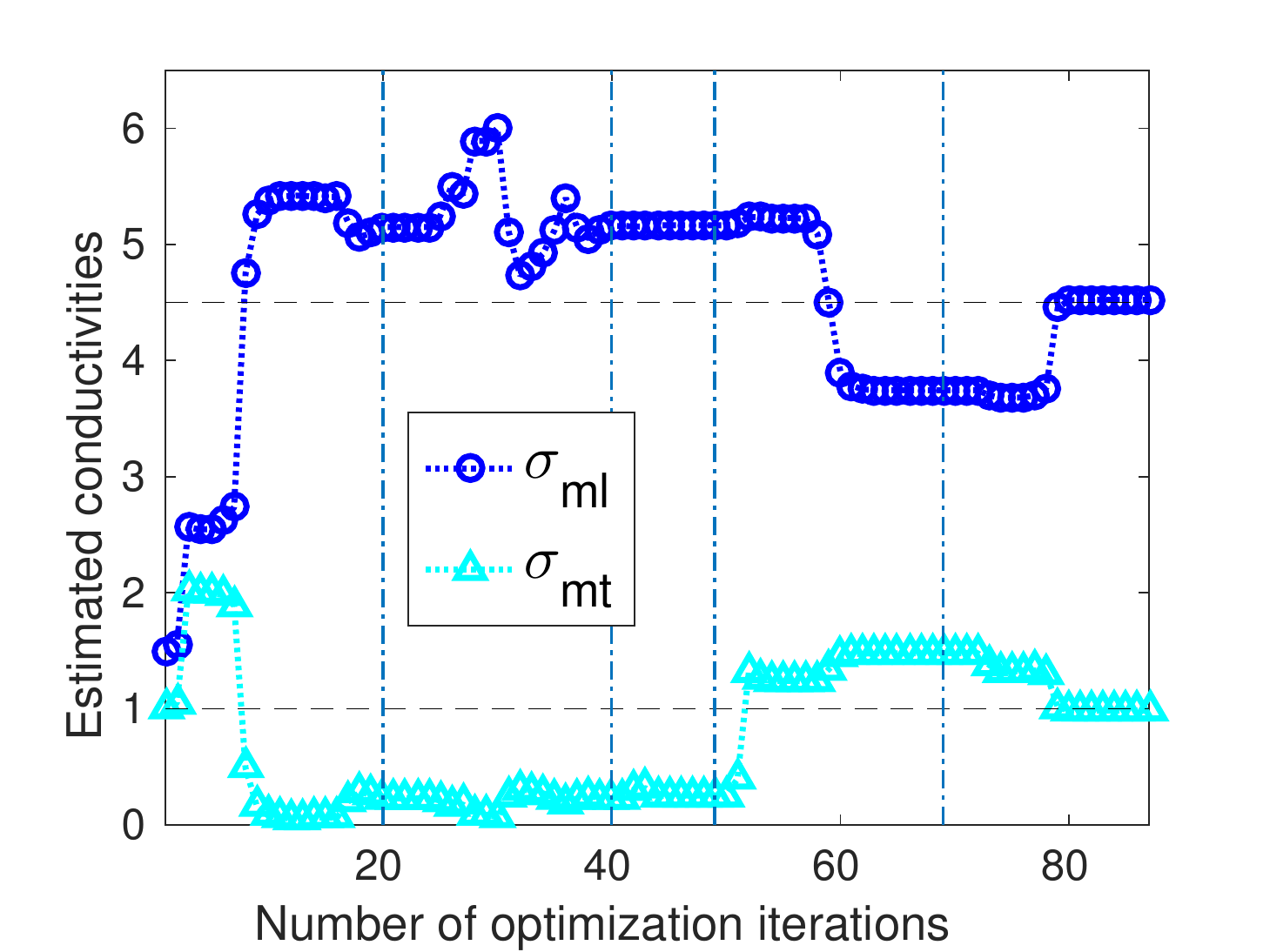}
\end{minipage}
\begin{minipage}{0.49\textwidth}
\includegraphics[scale=0.45]{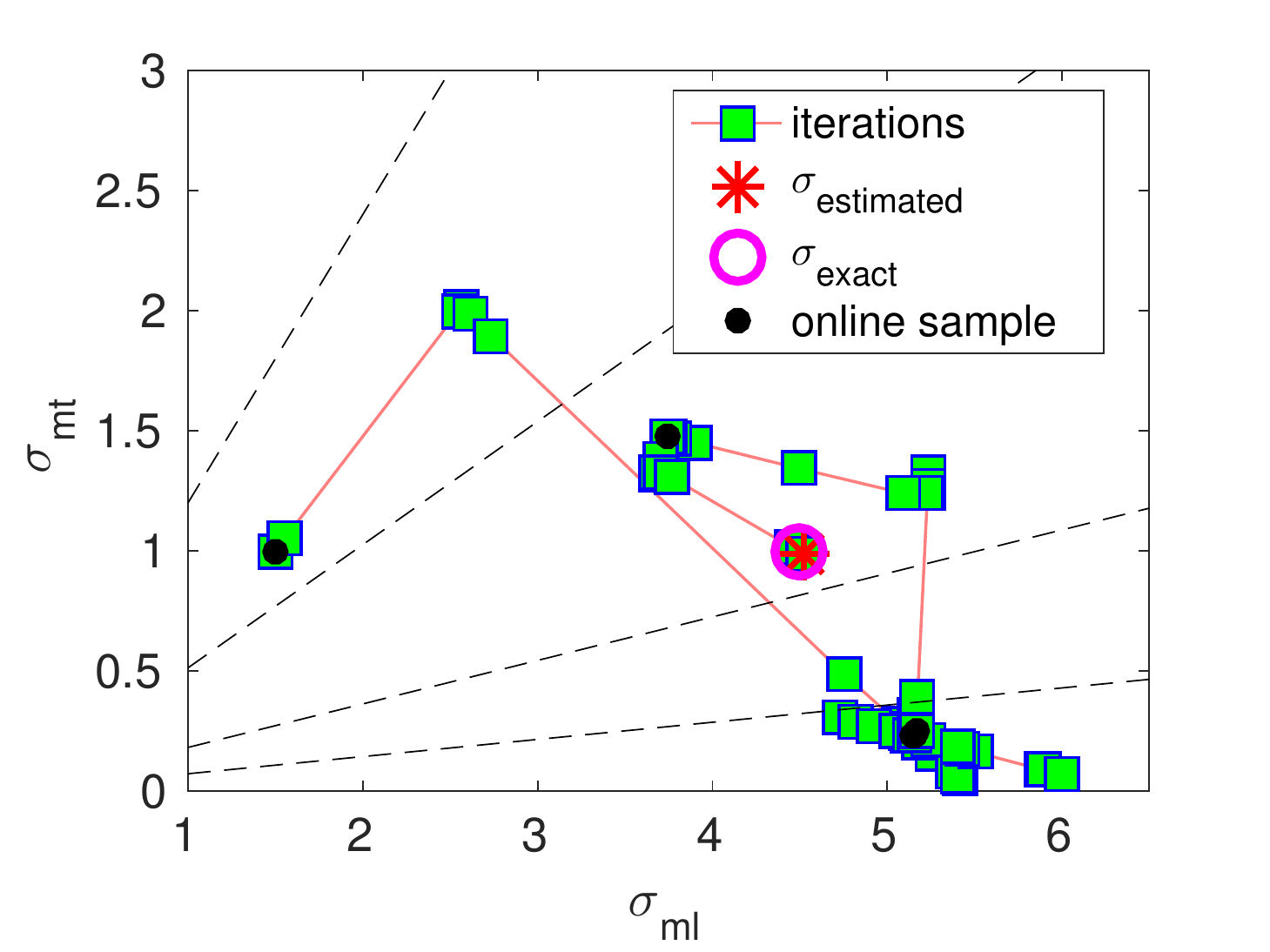}
\end{minipage}
\caption{Optimization with adaptive POD-DEIM corresponding to Table \ref{76k6testOnlineTable} for test point $\bsb{\sigma}_{\rm exact} = [4.5, 1]$. Left: Optimization iterations. The vertical dash lines denote the position where a new POD basis was generated.  Each optimization cycle (five in total) is constrained by a maximum iteration number 20. Right: Optimization trajectory. The points for online reduced basis generation are marked by the black dots.}
\label{onlinePOD4d51}
\end{center}
\end{figure}

\begin{figure}[tp]
\begin{center}
\begin{minipage}{0.49\textwidth}
\includegraphics[scale=0.45]{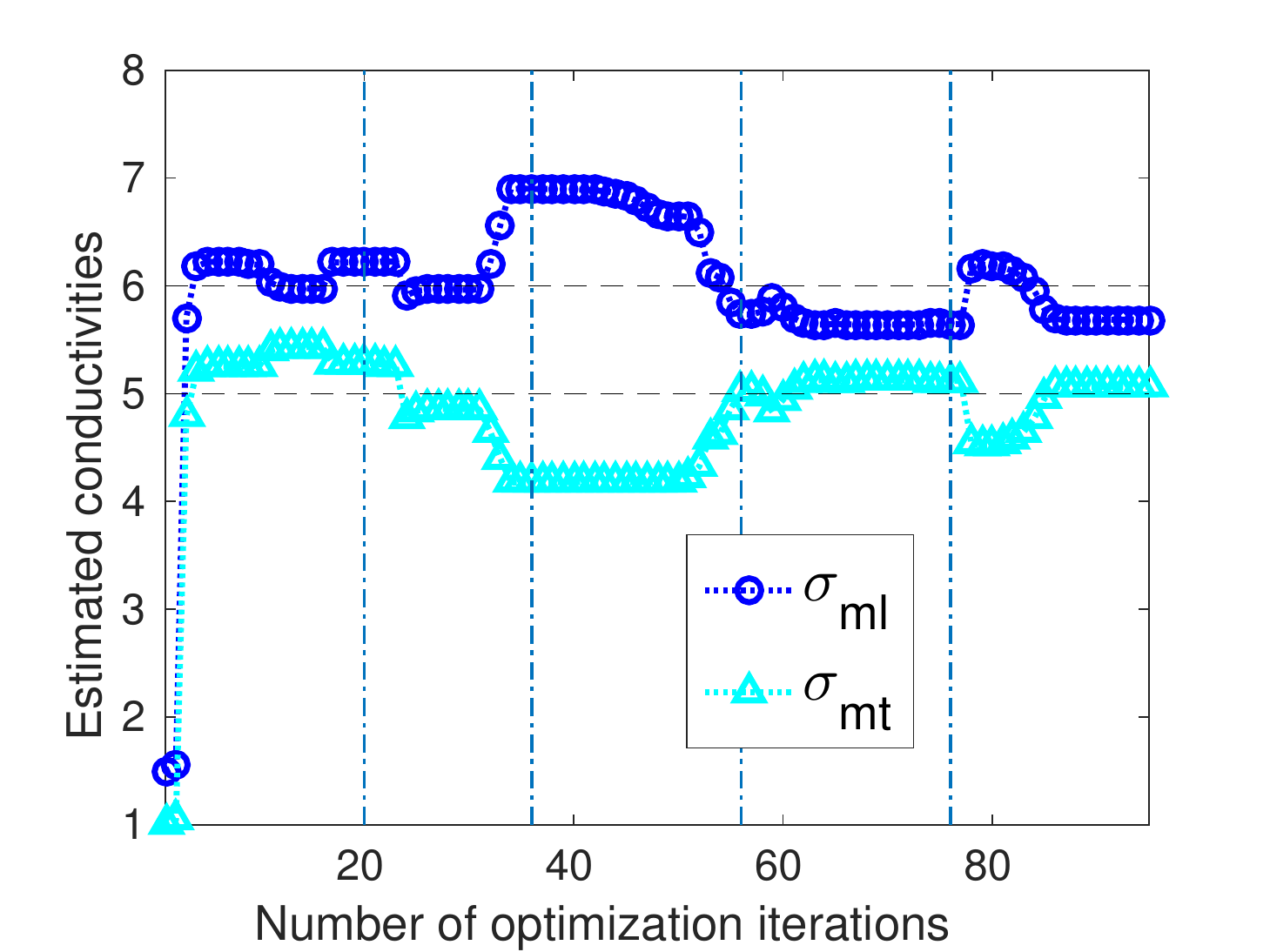}
\end{minipage}
\begin{minipage}{0.49\textwidth}
\includegraphics[scale=0.45]{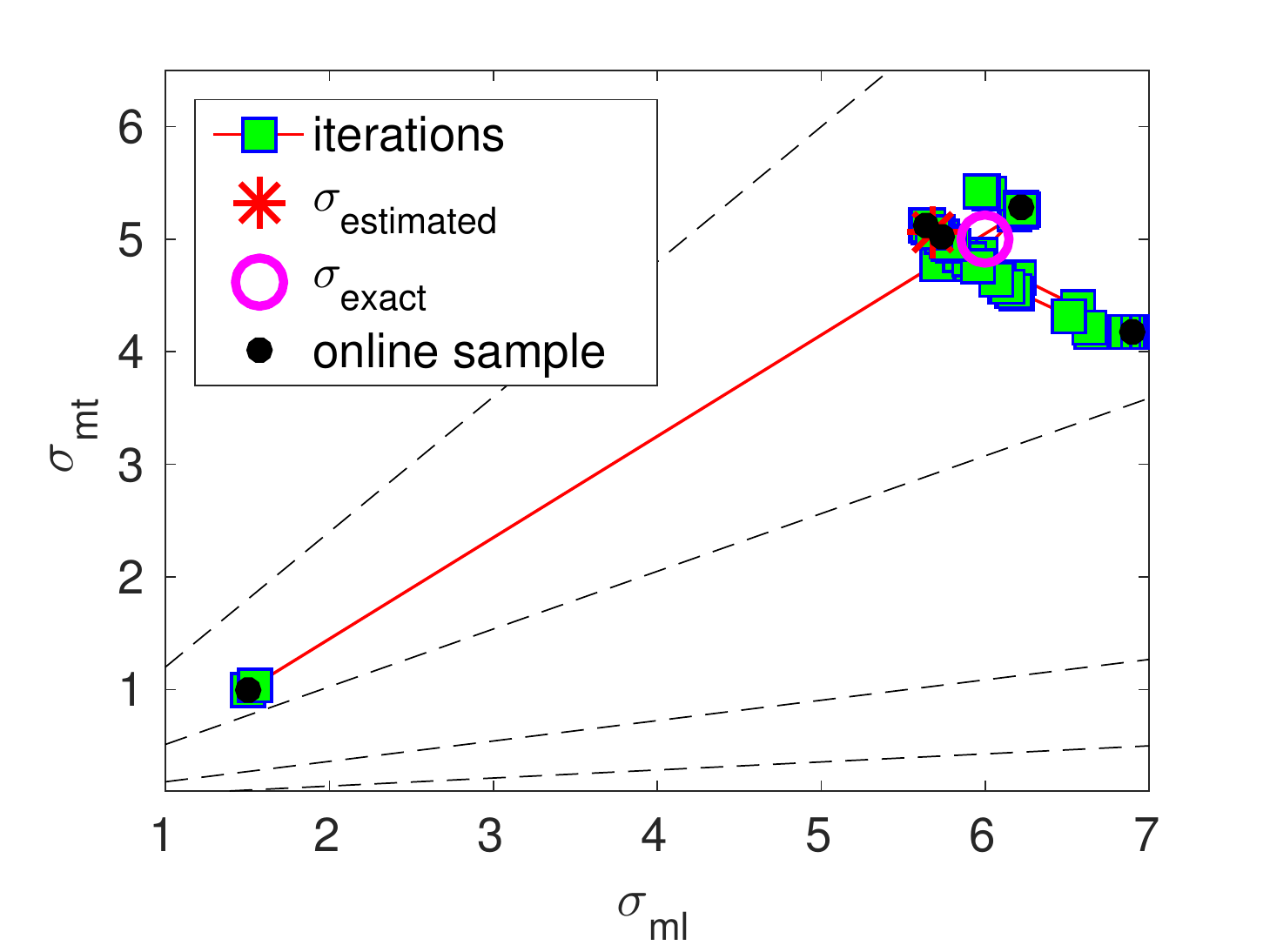}
\end{minipage}
\caption{Optimization with adaptive POD-DEIM corresponding to Table \ref{76k6testOnlineTable} for test point $\bsb{\sigma}_{\rm exact} = [6, 5]$. Left: Optimization iterations. The vertical dash lines denote the position where a new POD basis was generated.  Each optimization cycle (five in total) is constrained by a maximum iteration number 20. Right: Optimization trajectory. The points for online reduced basis generation are marked by the black dots.}
\label{onlinePOD65}
\end{center}
\end{figure}

The following experiments were conducted with the same six test values as before, on the slab mesh having 76832 nodes. The measurement data were still generated synthetically with extra 15\% noise. We execute Algorithm \ref{adaptivePOD-alg} given initial guess $\bsb{\sigma}_{\rm initial} = [1.5, 1]$ and particularly control the reduced-order minimization problem inside the loop by a maximum iteration number 20.  Each single online basis $\mathbb{Z}_{\rm u}$ has size 35 and $\mathbb{Z}_{\rm ion}$ has size 80. Accordingly to our experience, five cycles of reduced-order optimizations are generally enough to obtain stable convergence. We list the estimated conductivities in Table \ref{76k6testOnlineTable}. As one can see, the accuracy of conductivity estimation is promising, and the execution time is reduced to 10.2\%--38.8\%. Although the online computational reduction is not as intense as the offline-online procedure, this is an effective alternative use of model reduction especially when we have no offline reduced basis at hand.
To have an insight on the optimization trajectory, we report in Fig.~\ref{onlinePOD4d51} (resp.~Fig.~\ref{onlinePOD65}) the five cycles of iterations correspondingly to $\bsb{\sigma}_{\rm exact} = [4.5, 1]$ (resp.~$\bsb{\sigma}_{\rm exact} = [6, 5]$).

\begin{table}[tp]
 \begin{center}
\caption{Conductivity estimation by online adaptive POD-DEIM. }
{\footnotesize
\begin{tabular}{|c||c|c|c|c|c|}
\hline\hline
& $\boldsymbol{\sigma}_{\rm exact}$ & $\boldsymbol{\sigma}_{\rm estimated}$ & \# fwd $|$ bwd & Total exe.~time & Time perc\\
\hline
	Full Order & [3.2, 0.5] & [3.237, 0.625] &  51 $|$ 12  & 2720 s & 100\% \\
	Adaptive POD+DEIM  & [3.2, 0.5] & [3.223, 0.6635] &  544 $|$ 89  & 1055  s & 38.8\%\\ 
\hline
	Full Order & [4.5, 1] & [4.535, 1.099] &  58 $|$ 17  & 5409 s & 100\% \\  
	Adaptive POD+DEIM  & [4.5, 1] & [4.522, 0.9839] & 507 $|$ 87 &  733.5 s & 13.6\%\\    
\hline
	Full Order & [5.5, 3] & [5.538, 3.059] &  109 $|$ 36   & 8699 s& 100\% \\   
	Adaptive POD+DEIM  & [5.5, 3] & [5.294, 3.212] & 519 $|$ 93 &  889.4 s & 10.2\%\\     
\hline
	Full Order & [4, 2] & [4.045, 2.076] &  71 $|$ 23   & 5388 s& 100\% \\   
	Adaptive POD+DEIM  & [4, 2] & [3.9, 2.084] & 318 $|$ 81 &  812.6 s & 15.1\%\\     
\hline
	Full Order & [3, 2] & [3.055, 2.077] &  59 $|$ 19   & 5626 s & 100\% \\  
	Adaptive POD+DEIM  & [3, 2] & [3.237, 1.784] & 366 $|$ 80 &  995.9 s & 17.7\%\\ 
\hline
	Full Order & [6, 5] & [6.038, 5.049] &  50 $|$ 16   & 3366 s & 100\%\\   	
	Adaptive POD+DEIM  & [6, 5] & [5.686, 5.064] & 425 $|$ 95 &  709 s & 21\%\\ 
\hline
\end{tabular}
}
\label{76k6testOnlineTable}
\end{center}
\end{table}

\section{Conclusions}
Cardiac conductivity estimation is a critical step for bringing computational electrocardiology in clinical settings.  Currently there are no established procedures, available methodologies are computationally demanding. Model reduction is generally challenged by the mathematical features of the Monodomain problem and this prevents ``classical'' approaches like POD (coupled with DEIM) to be promptly applied. This paper gives a first contribution to model reduction applying  to the inverse conductivity problem. 

As a matter of fact, the main challenge of model reduction lies in the POD basis construction. To quantify this effect, we have introduced the concept of DOE. We have observed that the DOE of a POD basis based on a single generating parameter is narrow, especially when the amplitude of the conductivity value is small. This phenomena has also been shown in \cite{Boulakia2011POD} where POD was applied to estimate ionic model parameters. The situation in our case is even worse, considering the fact that we can not group snapshots for basis construction from different generating parameters, due to the strong sensitivity of the transmembrane potential to the conductivities. In other words, a unique POD basis is inadequate for the inverse conductivity problem, where the simulations are performed with various conductivity parameters. 

Another challenge that deserves to mention is the failure of adding sensitivity snapshots, which is expected to highly improve the effectiveness of the POD basis in other problems \cite{Carlberg2008}. We find that putting extra sensitivity snapshots into the snapshot matrix of state variables does not significantly enlarge the DOE of the constructed RB as much as we desire. Even if it did to some extent, the required number of POD modes is almost tripled, hence it is not appropriate for online application.

Nevertheless, there are still some interesting aspects that finally lead to a success. We detect that the DOE of a reduced basis is confined to an angular region of its generating parameter, and the region would be enlarged if the parameter amplitude gets larger. Based on this, we sample the parameter space utilizing the polar coordinates and the  Gaussian nodes. A sample set of size ten is then obtained.
The usage of multiple POD bases, each generated with a sampled parameter, provides satisfactory results. Overall, by utilizing this POD-DEIM reduced-order model, the computational effort can be reduced by at least 95\% in conductivity estimation.

This work opens several interesting challenges to be investigated in future works. In particular, we would like to theoretically quantify the conductivity estimation error caused by reduced-order modeling. An error analysis for some particular optimal control problems has been studied in \cite{Gubisch2016A-pos-33371,kammann2013posteriori}, but the work on conductivity estimation is still open sine the control parameter appears in the differential core of the system.
We also plan to extend the sampling strategy proposed here to 3D Bidomain inverse conductivity problem.

\section*{Acknowledgements}
This work has been supported by the NSF Project DMS 1412973/1413037 ``Novel Data Assimilation Techniques in Mathematical Cardiology''. We thank the referees for extremely valuable suggestions and remarks in improving this paper.

\bibliographystyle{model1b-num-names}
\bibliography{POD_inverseMono}

\end{document}